\documentclass[12pt]{article}
\pdfoutput=1

\usepackage{amsthm}
\usepackage{mathrsfs}
\usepackage{amsmath, amsfonts, amssymb}
\usepackage{graphicx}
\usepackage{psfrag}
\usepackage[usenames,dvipsnames,svgnames,table]{xcolor}
\usepackage{enumerate}
\usepackage{jheppubcustom}
\usepackage{soul}
\usepackage{subfig}
\usepackage{bbm}

\newcommand{\be}{\begin{equation}}
\newcommand{\ee}{\end{equation}}

\def\be{\begin{equation}}
\def\ee{\end{equation}}
\def\beq{\begin{equation}}
\def\eeq{\end{equation}}

\def\diam{{\rm diam}}
\def\dist{{\rm dist}}

\newtheorem{theorem}{Theorem}[section]
\newtheorem{lemma}[theorem]{Lemma}
\newtheorem{corollary}[theorem]{Corollary}
\newtheorem{proposition}[theorem]{Proposition}
\newtheorem{definition}[theorem]{Definition}

\newtheorem{remark}[theorem]{Remark}
\newtheorem{exercise}[theorem]{Exercise}

\def\Z{\mathbb Z }

\def\la{\langle}
\def\ra{\rangle}

\def\simleq{\; \raise0.3ex\hbox{$<$\kern-0.75em
      \raise-1.1ex\hbox{$\sim$}}\; }
   \def\simgeq{\; \raise0.3ex\hbox{$>$\kern-0.75em
      \raise-1.1ex\hbox{$\sim$}}\; }
      
      \newcommand{\figref}[1]{Fig.~\ref{#1}}

\newcommand{\bs}[1]{\boldsymbol{#1}}

\title{Scaling Limits, Brownian Loops and Conformal Fields}

\author{Federico Camia}

\emailAdd{f.camia@vu.nl, federico.camia@nyu.edu}

\affiliation{\it VU University Amsterdam, the Netherlands}
\affiliation{\it New York University Abu Dhabi, United Arab Emirates}

\begin{document}

\abstract{The main topic of these lecture notes is the continuum scaling limit of planar lattice models. One reason why this topic occupies an important place in the theory of probability and mathematical statistical physics is that scaling limits provide the link between statistical mechanics and Euclidean field theory. In order to explain the main ideas behind the concept of scaling limit, we will focus on a ``toy" model that exhibits the typical behavior of statistical mechanical models at and near the critical point. This model, known as the random walk loop soup, is actually interesting in its own right. It can be described as a Poisson process of lattice loops, or a lattice gas of loops since it fits within the ideal gas framework of statistical mechanics. After introducing the model and discussing some interesting connections with the discrete Gaussian free field, we will present some results concerning its scaling limit, which leads to a Poisson process of continuum loops known as the Brownian loop soup. The latter was introduced by Lawler and Werner and is a very interesting object with connections to the Schramm-Loewner Evolution (SLE) and various models of statistical mechanics. In the second part of these lecture notes, we will use the Brownian loop soup to construct a family of functions that behave like correlation functions of primary fields in Conformal Field Theory (CFT). We will then use these functions and their derivation to introduce the concept of conformal field and to explore the connection between scaling limits and conformal fields. }

\maketitle

\section{Introduction} \label{intro}

\subsection{Critical scaling limits} \label{intro-critical}
One of the main goals of both probability theory and statistical physics
is to understand and describe the behavior of random systems
consisting of a very large number of components (atoms, molecules, pixels, individuals, $\ldots$)
where the effect of each single component is negligible and the behavior of the system as
a whole is determined by the combined effect of all its components.
One usually wishes to study the behavior of such systems via some
\emph{observables}, suitably defined quantities that can be of an
analytic or geometric nature.
The asymptotic (in the number of system components) behavior of these quantities is
often deterministic, but in some interesting cases it turns out to be random.
This type of macroscopic randomness can be observed in \emph{critical systems},
i.e., systems at a continuous phase transition point (the \emph{critical point}).

In the physical theory of critical systems, it is usually assumed
that, when a system approaches the critical point, it is characterized
by a single length scale (the \emph{correlation length}) in terms of
which all other lengths should be measured.
When combined with the experimental observation that the
correlation length diverges at the critical point, this simple
but strong assumption, known as the \emph{scaling hypothesis},
leads to the belief that a critical system has no characteristic
length, and is therefore invariant under scale transformations.
This implies that all thermodynamic functions at criticality
are homogeneous functions, and predicts the appearance of
power laws.


It also suggests that, for models of critical systems realized on a lattice,
one can attempt to take a \emph{continuum scaling limit} in which the mesh
of the lattice is sent to zero while focus is kept on macroscopic
observables that capture the large scale behavior.
In the limit, the discrete model should converge to a continuum one that
encodes the large scale properties of the original model, containing at
the same time more symmetry. In many cases, this allows for the derivation
of additional insight by combining methods of discrete mathematics with
considerations inspired by the continuum limit picture. The simplest
example of such a continuum random model is Brownian motion,
which is the scaling limit of the simple random walk.
In general, though, the complexity of the discrete model makes
it impossible even to guess the nature of the scaling limit,
unless some additional feature can be used to pin down properties of
the continuum limit.  Two-dimensional critical systems belong to the
class of models for which this can be done, and the additional feature
is \emph{conformal invariance}.

Indeed, thanks to the work of Polyakov~\cite{polyakov}
and others~\cite{bpz1,bpz2}, it was understood by physicists
since the early seventies that critical statistical mechanical
models should possess continuum scaling limits with a global
conformal invariance that goes beyond simple scale invariance
(as long as the discrete models have ``enough'' rotation invariance),
but this remains a conjecture for most models describing critical systems.

The appearance of conformal invariance in the scaling limit,
besides being an interesting phenomenon in its own right,
poses constraints on the possible continuum limits.
However, since the conformal group is in general a finite
dimensional Lie group, the resulting constraints are limited
in number and provide only limited information.
The phenomenon becomes particularly interesting in two dimensions,
where every analytic function defines a conformal transformation
(at points where its derivative is non-vanishing), and the conformal
group is infinite-dimensional. 

After this observation was made, a large number of critical
problems in two dimensions were analyzed using conformal
methods, which were applied, among others, to Ising and Potts
models, Brownian motion, Self-Avoiding Walk (SAW), percolation,
and Diffusion Limited Aggregation (DLA). The resulting body of knowledge
and techniques, starting with the work of Belavin, Polyakov and Zamolodchikov~\cite{bpz1,bpz2}
in the early eighties, goes under the name of Conformal Field Theory (CFT).

At the hands of theoretical physicists, CFT has been very
successful in producing many interesting results which until
recently remained beyond any rigorous mathematical justification.
%
%
This has changed in the last fifteen years, with the emergence in the mathematics
literature of important developments in the area of two-dimensional critical systems
which have followed a completely new direction, providing new tools and a new way
of looking at critical systems and related conformal field theories.

These developments came on the heels of interesting
results on the scaling limits of discrete models (see, e.g.,
the work of Aizenman \cite{aizenman1,aizenman2},
Benjamini-Schramm \cite{bs}, Aizenman-Burchard \cite{ab},
Aizenman-Burchard-Newman-Wilson \cite{abnw},
Kenyon \cite{kenyon1,kenyon2} and Aizenman-Duplantier-Aharony \cite{ada}),
but they differ greatly from those because they are based on a radically new
approach whose main tool is the Stochastic Loewner Evolution (SLE),
or Schramm-Loewner Evolution, as it is also known, introduced
by Schramm~\cite{schramm}.
The new approach, which is probabilistic in nature,
focuses directly on non-local structures that characterize
a given system, such as cluster boundaries in Ising, Potts
and percolation models, or loops in the $O(n)$ model.
At criticality, these non-local objects become, in the
continuum limit, random curves whose distributions can be
uniquely identified thanks to their conformal invariance
and a certain Markovian property.
There is a one-parameter family of SLEs, indexed by
a positive real number $\kappa$, and they appear to be
essentially the only possible candidates for the scaling
limits of interfaces of two-dimensional critical systems
that are believed to be conformally invariant.

The identification of the scaling limit of interfaces of critical
lattice models with SLE curves has led to tremendous progress in
recent years. The main power of SLE stems from the fact that it
allows to compute different quantities; for example, percolation
crossing probabilities and various percolation critical exponents.
Therefore, relating the scaling limit of a critical lattice model
to SLE allows for a rigorous determination of some aspects of the
large scale behavior of the lattice model. For the mathematician,
the biggest advantage of SLE over CFT lies maybe in its mathematical
rigor. But many physicists working on critical phenomena and CFT have
promptly recognized the importance of SLE and added it to their toolbox.

In the context of the Ising, Potts and $O(n)$ models, as well as percolation,
an SLE curve is believed to describe the scaling limit of a single interface,
which can be obtained by imposing special boundary conditions.
A single SLE curve is therefore not in itself sufficient to immediately describe
the scaling limit of the unconstrained model without boundary conditions in the
whole plane (or in domains with boundary conditions that do not determine a single
interface), and contains only limited information concerning the connectivity
properties of the model.

A more complete description can be obtained in terms of loops, corresponding
to the scaling limit of cluster boundaries. Such loops should also be random
and have a conformally invariant distribution, closely related to SLE.
This observation led 
to the definition of Conformal Loop Ensembles (CLEs) \cite{werner1,sheffield-cle,sw},
which are, roughly speaking, random collections of fractal loops with a certain
``conformal restriction property." As for SLE, there is a one parameter family of CLEs.


Several interesting models of statistical mechanics, such as percolation and the Ising
and Potts models, can be described in terms of clusters. In two dimensions and at the
critical point, the scaling limit geometry of the boundaries of such clusters is known
(see \cite{smirnov-percolation,cn06,cn07,smirnov-ising-full}) or conjectured (see
\cite{kn}) to be described by some member of the one-parameter family of Schramm-Loewner
Evolutions (SLE$_{\kappa}$ with $\kappa>0$) and related Conformal Loop Ensembles
(CLE$_{\kappa}$ with $8/3<\kappa<8$).
SLEs can be used to describe the scaling limit of single interfaces; CLEs are collections
of loops and are therefore suitable to describe the scaling limit of the collection of all
macroscopic boundaries at once. For example, the scaling limit of the critical percolation
exploration path is SLE$_6$ \cite{smirnov-percolation,cn07}, and the scaling limit of the
collection of all critical percolation interfaces in a bounded domain is CLE$_6$ \cite{cn06,cn08}.
For $8/3 < \kappa \leq 4$, CLE$_\kappa$ can be obtained \cite{sw} from the Brownian
loop soup introduced by Lawler and Werner \cite{lw}.

\subsection{Near-critical scaling limits} \label{intro-near-critical}
A meaningful continuum scaling limit that differs from the critical one can usually be obtained
by considering a system ``\emph{near}'' the critical point, the so-called \emph{off-critical} regime.
This is done by adjusting some parameter of the model and sending it to the critical value at
a specific rate while taking the scaling limit, in such a way that the correlation length, in
macroscopic units, stays bounded away from 0 and $\infty$. Depending on the context, this
situation is described as {\em near-critical}, {\em off-critical} or {\em massive} scaling limit.
The term \emph{massive} refers to the persistence of a macroscopic correlation length, which should
give rise to what is known in the physics literature as a \emph{massive} field theory. In Euclidean
field theory, the term \emph{massless} is used to describe a scale-invariant system, while the
term \emph{massive} refers to a system with exponential decay of correlations.

Near-critical scaling limits are in general not expected to be conformally invariant, since the
persistence of a macroscopic correlation length implies the absence of scale invariance, and
are instead expected to have exponentially decaying correlations, while ``resembling'' their
critical counterparts at distances smaller than the correlation length. (For a rigorous example
of this behavior, see \cite{cjm}.) The lack of conformal invariance implies that the geometry
of near-critical scaling limits cannot be described by an SLE or CLE. Indeed, much less has been
proved about the geometry of off-critical models than about that of critical ones.

\subsection{Random walk loop soups and Brownian loop soups} \label{intro-soups}
Symanzik, in his seminal work on Euclidean quantum field theory~\cite{symanzik}, introduced a
representation of the $\phi^4$ Euclidean field as a ``gas'' of weakly interacting random paths.
The use of random paths in the analysis of Euclidean field theories and statistical mechanical models was
subsequently developed by various authors, most notably 
Brydges, 
Fr\"ohlich, 
Spencer and 
Sokal~\cite{bfsp,bfso}, and 
Aizenman~\cite{aizenman}, proving extremely useful (see~\cite{ffs-book} for a comprehensive account).
The probabilistic analysis of Brownian and random walk paths and associated local times was carried
out by Dynkin~\cite{dynkin1,dynkin2}.
More recently, ``gases'' or ``soups'' (i.e., Poissonian ensembles) of Brownian and random walk
loops have been extensively studied in connection with SLE and the Gaussian free field (see,
e.g., \cite{werner1,lw,werner2,lejan1,lejan2,sznitman-notes}).

In these lecture notes, we provide a prototypical example, amenable to rigorous mathematical
analysis, of critical and off-critical behavior by studying the scaling limit of an ``ideal gas'' of loops,
called \emph{random walk loop soup}. The loops are weighed according to a random walk measure
on the square lattice with killing rates $\{ k_x \}_{x \in {\mathbb Z}^2}$. Lawler and Trujillo Ferreras
\cite{ltf} have shown that, in the absence of killing, the random walk loop soup converges in the
scaling limit to the \emph{Brownian loop soup} introduced by Lawler and Werner \cite{lw}.

The Brownian loop soup is, roughly speaking, a Poisson point process with intensity measure $\lambda\mu$,
where $\lambda$ is a positive constant and $\mu$ is the \emph{Brownian loop measure} studied in \cite{werner3}.
$\mu$ is uniquely determined (up to a multiplicative constant) by its \emph{conformal restriction property},
a combination of conformal invariance and the property that the measure in a subdomain is the original
measure restricted to loops that stay in that subdomain). A realization of the loop soup consists of a countable
collection of loops. Given a bounded domain $D$, there is an infinite number of loops that stay in $D$; however,
the number of loops in $D$ of diameter at least $\varepsilon>0$ is finite. A consequence of conformal invariance
is scale invariance: if $\cal A$ is a realization of the Brownian loop soup and each loop is scaled in space by
$1/N$ and in time by $1/N^2$, the resulting configuration also has the distribution of the Brownian loop soup.


The Brownian loop soup exhibits a connectivity phase transition in the parameter $\lambda>0$ that
multiplies the intensity measure. When $\lambda \leq 1/2$, the loop soup in $D$ is composed of disjoint
clusters of loops \cite{werner1,werner2,sw} (where a cluster is a maximal collection of loops that intersect
each other). When $\lambda>1/2$, there is a unique cluster~\cite{werner1,werner2,sw} and the set of
points not surrounded by a loop is totally disconnected (see~\cite{bc2}). Furthermore, when
$\lambda \leq 1/2$, the outer boundaries of the loop soup clusters are distributed like Conformal Loop
Ensembles (CLE$_{\kappa}$) \cite{werner1,sheffield-cle,sw} with $8/3 < \kappa \leq 4$. As mentioned
earlier, the latter are conjectured to describe the scaling limit of cluster boundaries in various critical models
of statistical mechanics, such as the critical Potts models for $q \in (1,4]$.

More precisely, if $8/3 < \kappa \leq 4$, then $0 < \frac{(3 \kappa -8)(6 - \kappa)}{4 \kappa} \leq 1/2$
and the collection of all outer boundaries of the clusters of the Brownian loop soup with
intensity parameter $\lambda = \frac{(3 \kappa -8)(6 - \kappa)}{4 \kappa}$ is distributed
like $\text{CLE}_{\kappa}$ \cite{sw}. For example, the continuum scaling limit
of the collection of all macroscopic boundaries of critical Ising spin clusters is conjectured
to correspond to $\text{CLE}_3$ and to a Brownian loop soup with $\lambda=1/4$.

We note that there is some confusion in most of the existing literature regarding the critical intensity
corresponding to the connectivity phase transition in the Brownian loop soup, as well as regarding the
relation between $\lambda$ and $\kappa$. This is due to the fact that the Brownian loop measure is
an infinite measure. As a consequence, a choice of normalization is required when defining the Brownian
loop soup. This choice is then reflected in the relation between $\lambda$ and $\kappa$.
Our choice of normalization is made explicit in Section \ref{sec:bls} and coincides with that of \cite{lw}.

The intensity $\lambda$ of the Brownian loop soup is related to the central charge
$c$ of the corresponding statistical mechanical model. A discussion of the central charge of the Brownian
loop soup can be found in Section 6 of \cite{cgk15}. Here we will just mention that, in some vague sense,
$\lambda$ determines how much the system feels a change in the shape of the domain. To understand
this, suppose that ${\cal A}(\lambda,D)$ is a realization of the Brownian loop soup in $D$ with intensity
$\lambda$, and consider a subdomain $D' \subset D$. By removing from ${\cal A}(\lambda,D)$ all loops
that are not contained in $D'$, one obtains the loop soup ${\cal A}(\lambda,D')$ in $D'$ with the same
intensity $\lambda$. (This property of the Brownian loop soup follows from its Poissonian nature and the
conformal restriction property of the Brownian loop measure $\mu$.) The number of loops removed is
stochastically increasing in $\lambda$, so that larger values of $\lambda$ imply that the system is more
sensitive to changes in the shape of the domain. In some sense, the loops can be seen as mediators of
correlations from the boundary of the domain. (A precise formulation of this observation is presented
in Section~\ref{sec:gff}, in the context of the discrete Gaussian free field.) We note that the change
from ${\cal A}(\lambda,D)$ to ${\cal A}(\lambda,D')$ has a nonlocal effect since the loops that are
removed are extended objects, and that even a small local change to the shape of the domain can
have an effect very far away, due to the scale invariance of the Brownian loop soup. This is a
manifestation of the criticality of the system.

If one considers random walk loop soups with killing, the situation is very different. If one
scales space by $1/N$ and time by $1/N^2$ keeping the killing rates constant, the resulting scaling
limit is trivial, in the sense that it does not contain any loops larger than one point. This is so because,
under the random walk loop measure, only loops of duration of order at least $N^2$ have diameter of
order at least $N$ with non-negligible probability as $N \to \infty$, and are therefore ``macroscopic''
in the scaling limit. It is then clear that, in order to obtain a nontrivial scaling limit, the killing rates need
to be rescaled as well. When that is done appropriately, the scaling limit yields a \emph{massive} (non-scale-invariant)
version of the Brownian loop soup \cite{camia15}. This \emph{massive Brownian loop soup} is perhaps
the simplest modification of the massless Brownian loop soup of Lawler and Werner combining several
properties that are considered typical of off-critical systems, including exponential decay of correlations
and a weaker form of conformal symmetry called \emph{conformal covariance}. The mechanism
by which it arises from the Brownian loop soup of Lawler and Werner is analogous to that appearing in the
standard Brownian motion representation of the continuum Gaussian free field when a mass term is present.
Moreover, it is the only possible non-trivial, near-critical scaling limit of a random walk loop soup with killing.

The random walk loop soups studied in these lecture notes are shown to conform to the picture
described in Sections \ref{intro-critical} and \ref{intro-near-critical},
providing explicit examples of critical and near-critical scaling limits satisfying the general features
typically associated with them. Because of their simplicity and amenability to a full analysis, the loop
soups studied in these notes can be considered a prototypical example of a system with critical and
off-critical behavior. Moreover, the Brownian loop soup is deeply related to SLE and the critical scaling
limit of models of statistical mechanics, as explained above, while its massive counterpart could
potentially play a role in the description of various near-critical scaling limits.

\subsection{Conformal correlation functions in the Brownian loop soup}
There is a rich but still not fully understood connection between the conformal field theories studied
by physicists and conformally invariant stochastic models such as SLE and the Brownian loop soup.
Physicists are  often interested in conformal field theories defined by a Lagrangian and focus on the
correlation functions (i.e., expectations of products) of local ``observables'' or ``operators.'' In field
theories related to the scaling limit of lattice models, the local observables typically considered
are sums of (products of)  local (random) variables over some macroscopic region, for example the
local magnetization in the Ising model. In the scaling limit, if the lattice model is critical, such observables
may converge to ``conformal fields.'' Mathematically, these ``fields'' are random generalized
functions (distributions in the sense of Schwartz).

By contrast, conformal stochastic models are often defined and studied using very different methods
and with different goals. In particular, the Brownian loop soup and SLE describe the behavior of
non-local observables, macroscopic random curves such as Ising cluster boundaries.

In the last part of these notes, based on \cite{cgk15}, we use the Brownian loop soup to construct
families of functions that behave like correlation functions of conformal fields, and explore some
mathematical questions related to the scaling limit of critical correlation functions and to conformal
fields. We expect that these correlation functions may define a full-fledged conformal field theory
(one that could be related to the physics of de Sitter spacetime and to eternal inflation).
While this has not been fully established, progress has been made \cite{cl16} in the direction of
constructing random generalized functions with the appropriate correlation functions. The putative field
theory would have several novel features, such as a periodic spectrum of conformal dimensions \cite{cgk15}.

\section*{Acknowledgements}
It is a pleasure to thank the organizers of the program ``Disordered systems, random spatial processes
and some applications'' (IHP, Paris, 5 January -- 3 April, 2015) for their invitation to give a series of lectures,
which resulted in these notes. Part of the content of these notes derives from work conducted with Alberto
Gandolfi and Matthew Kleban, and part of it profited from discussions and previous work with Erik Broman;
I wish to thank all of them for the enjoyable collaborations. Several discussions with Tim van de Brug and
Marcin Lis helped to improve the presentation of the material. My gratitude also goes to Claire Berenger
and all the IHP staff for their friendly efficiency, and to an anonymous referee for a useful suggestion. The
work of the author was supported in part by the Netherlands Organization for Scientific Research (NWO)
through grant Vidi 639.032.916.


\section{Loop soups} \label{sec:definitions}

\subsection{Random walk loop soups} \label{sec:rwls}

Let $k_x \geq 0$ for every $x \in {\mathbb Z}^2$ and define $p_{x,y} = 1/(k_x + 4)$ if $|x-y|=1$ and
$p_{x,y}=0$ otherwise. If $k_x=0$ for all $x$, $\{p_{x,y}\}_{y \in {\mathbb Z}^2}$ is the collection of
transition probabilities for the simple symmetric random walk on ${\mathbb Z}^2$. If $k_x \neq 0$, then
$p_{x,y} = \frac{1}{4}(1+\frac{k_x}{4})^{-1}<\frac{1}{4}$ and one can interpret
$\{p_{x,y}\}_{y \in {\mathbb Z}^2}$ as the collection of transition probabilities for a random walker
``killed'' at $x$ with probability $1-(1+\frac{k_x}{4})^{-1}=\frac{k_x}{k_x+4}$. (Equivalently, one can
introduce a ``cemetery'' state $\Delta$ not in ${\mathbb Z}^2$ to which the random walker jumps from
$x \in {\mathbb Z}^2$ with probability $\frac{k_x}{k_x+4}$, and where it stays forever once it reaches it.)
Because of this interpretation, we will refer to the collection ${\bf k}=\{k_x\}_{x\in{\mathbb Z}^2}$ as
\emph{killing rates}.


Given a $(2n+1)$-tuple $(x_0,x_1,\ldots,x_{2n})$ with $x_0=x_{2n}$ and $|x_i-x_{i-1}|=1$ for $i=1,\ldots,2n$,
we call \emph{rooted lattice loop} the continuous path $\tilde\gamma: [0,2n] \to {\mathbb C}$ with $\tilde\gamma(i)=x_i$
for integer $i=0,\ldots,2n$ and $\tilde\gamma(t)$ obtained by linear interpolation for other $t$. We call $x_0$ the
\emph{root} of the loop and denote by $|\tilde\gamma|=2n$ the \emph{length} or \emph{duration} of the loop.

Now let $D$ denote either ${\mathbb C}$ or a connected subset of ${\mathbb C}$. 
Following Lawler and Trujillo Ferreras \cite{ltf}, but within the more general framework of the previous paragraph,
we introduce the \emph{rooted random walk loop measure} $\nu^{r,{\bf k}}_D$ which assigns the loop $\tilde\gamma$ of
length $|\tilde\gamma|$, with root $x_0$, weight
\be \nonumber
|\tilde\gamma|^{-1} p_{x_0,x_1} p_{x_1,x_2} \ldots p_{x_{|\tilde\gamma|-1},x_0}
= |\tilde\gamma|^{-1} \prod_{i=0}^{|\tilde\gamma|-1} \left(k_{x_i}+4 \right)
= |\tilde\gamma|^{-1} \prod_{x \in \tilde\gamma} \left(k_{x}+4 \right)^{-n(x,\tilde\gamma)}
\ee
if $x_0,\ldots,x_{|\tilde\gamma|-1} \in D$ and 0 otherwise, where $x \in \tilde\gamma$ means that $x$ is visited by
$\tilde\gamma$ and $n(x,\tilde\gamma)$ denotes the number of times $\tilde\gamma$ makes a jump from $x$.
(Note that $n(x,\tilde\gamma)$ equals the number of visits of $\tilde\gamma$ to $x$ for all vertices visited by $\tilde\gamma$
except the root $x_0$. For the root $x_0$, $n(x,\tilde\gamma)$ is the number of visits minus one since the walk ends
at $x_0$, so the last visit doesn't count.)

The \emph{unrooted random walk loop measure} $\nu^{u,{\bf k}}_D$ is obtained from the rooted one by
``forgetting the root.'' More precisely, if $\tilde\gamma$ is a rooted lattice loop and $j$ a positive integer,
$\theta_j \tilde\gamma : t \mapsto \tilde\gamma(j+t \mod |\tilde\gamma|)$ is again a rooted loop.
This defines an equivalence relation between rooted loops; an \emph{unrooted lattice loop} is an equivalence
class of rooted lattice loops under that relation. By a slight abuse of notation, in the rest of the paper we will
use $\tilde\gamma$ to denote unrooted lattice loops and $\tilde\gamma(\cdot)$ to denote any rooted
lattice loop in the equivalence class of $\tilde\gamma$.
The $\nu^{u,{\bf k}}_D$-measure of the unrooted loop $\tilde\gamma$ is the sum of the $\nu^{r,{\bf k}}_D$-measures
of the rooted loops in the equivalence class of $\tilde\gamma$. 
The \emph{length} or \emph{duration}, $|\tilde\gamma|$, of an unrooted loop $\tilde\gamma$ is the length
of any one of the rooted loops in the equivalence class $\tilde\gamma$.


\begin{definition} \label{def:rwls}
A \emph{random walk loop soup} in $D$ with intensity $\lambda$ is a Poisson realization from $\lambda\nu^{u,{\bf k}}_D$.
\end{definition}

A realization of the random walk loop soup in $D$ is a multiset (i.e., a set whose elements can occur multiple times)
of unrooted loops. If we denote by $N_{\tilde\gamma}$ the multiplicity of $\tilde\gamma$ in a loop soup with
intensity $\lambda$, $\{N_{\tilde\gamma}\}$ is a collection of independent Poisson random variables with
parameters $\lambda\nu^{u,{\bf k}}_D(\tilde\gamma)$. Therefore, the probability that a realization of the
random walk loop soup in $D$ with intensity $\lambda$ contains each loop $\tilde\gamma$ in $D$ with
multiplicity $n_{\tilde\gamma}\geq0$ is equal to
\begin{equation} \label{rwls}
\prod_{\tilde\gamma} \exp{\left(-\lambda\nu^{u,{\bf k}}_D(\tilde\gamma)\right)}
\frac{1}{n_{\tilde\gamma}!} \left(\lambda \nu^{u,{\bf k}}_D(\tilde\gamma)\right)^{n_{\tilde\gamma}}
= \frac{1}{{\cal Z}_D^{\lambda,{\bf k}}} \prod_{\tilde\gamma} \frac{1}{n_{\tilde\gamma}!} \left(\lambda \nu^{u,{\bf k}}_D(\tilde\gamma)\right)^{n_{\tilde\gamma}} \, ,
\end{equation}
where the product $\prod_{\tilde\gamma}$ is over all unrooted lattice loops in $D$ and
\begin{equation} \label{partition-function}
{\cal Z}_D^{\lambda,{\bf k}} := \exp{\left(\lambda \sum_{\tilde\gamma} \nu^{u,{\bf k}}_D(\tilde\gamma)\right)}
= 1 + \sum_{n=1}^{\infty} \frac{1}{n!} \sum_{(\tilde\gamma_1,\ldots,\tilde\gamma_n)} \prod_{i=1}^n \lambda \nu^{u,{\bf k}}_D(\tilde\gamma_i) \, ,
\end{equation}
where the sum over $(\tilde\gamma_1,\ldots,\tilde\gamma_n)$ is over all ordered configurations of $n$ loops,
not necessarily distinct.
From a statistical mechanical viewpoint, ${\cal Z}_D^{\lambda,{\bf k}}$ can be interpreted as the grand canonical
partition function of a ``gas'' of loops, and one can think of the random walk loop soup as describing a grand
canonical ensemble of noninteracting loops (an ``ideal gas'') with the killing rates $\{k_x\}$ and the intensity
$\lambda$ as free ``parameters.''
(For more on the statistical mechanical interpretation of the model, see Section~6.4 of~\cite{bb}.)

\begin{exercise} \label{exercise:generating-function}
Compute the probability $p_{\ell}(A)$ that a random walk loop soup configuration contains $\ell$ loops from a given set $A$
and use the resut to compute the probability generating function $\sum_{\ell=0}^{\infty} p_{\ell} x^{\ell}$.
\end{exercise}

When $k_x=0 \; \forall x \in D \cap {\mathbb Z}^2$, we use $\nu^u_D$ to denote the unrooted random walk
loop measure in $D$; for reasons that will be clear when we talk about scaling limits, later in this section,
a random walk loop soup obtained using such a measure will be called \emph{critical}.

Now let $m: {\mathbb C} \to {\mathbb R}$ be a nonnegative function; we say that a random walk loop soup
has \emph{mass (function)} $m$ if $k_x = 4(e^{m^2(x)}-1)$ for all $x \in D \cap {\mathbb Z}^2$, and call
\emph{massive} a random walk loop soup with mass $m$ that is not identically zero on $D \cap {\mathbb Z}^2$.
For a massive random walk loop soup in $D$ with intensity $\lambda$ and mass $m$ we use the notation
$\tilde{\cal A}(\lambda,m,D)$.

The next proposition gives a construction for generating a massive random walk loop soup from a critical one,
establishing a useful probabilistic coupling between the two (i.e., a way to construct the two loop soups on the
same probability space).
\begin{proposition} \label{coupling2}
A random walk loop soup in $D$ with intensity $\lambda$ and mass function $m$ can be realized in the following way.
\begin{enumerate}
\item Take a realization of the \emph{critical} random walk loop soup in $D$ with intensity $\lambda$.
\item Assign to each loop $\tilde\gamma$ an independent, mean-one, exponential random variable $T_{\tilde\gamma}$.
\item Remove from the soup the loop $\tilde\gamma$ of length $|\tilde\gamma|$ if
\begin{equation} \label{eq:loop-killing-2}
\sum_{i=0}^{|\tilde\gamma|-1} m^2(\tilde\gamma(i)) > T_{\tilde\gamma} \, .
\end{equation}
\end{enumerate}
\end{proposition}

\begin{remark}
Note that equation~\eqref{eq:loop-killing-2} requires choosing a rooted loop from the equivalence class
$\tilde\gamma$ but is independent of the choice.
\end{remark}

\noindent{\bf Proof.}
The proof is analogous to that of Proposition~\ref{coupling1}, so we leave the details to the reader.
Below we just compare the expected number of loops in $D$ generated by the construction of
Proposition~\ref{coupling2} with that of the massive random walk loop soup, to verify the relation
between killing rates $\{k_x\}$ and mass function $m$. Since different loops are independent, we
can compare the expected numbers of individual loops.

Writing $p_{x,y} = \frac{1}{k_x+4} = \frac{1}{4} \frac{4}{k_x+4}$ when $|x-y|=1$,
for the massive random walk soup and loop $\tilde\gamma$, we have
\begin{equation} \label{mean1}
\lambda \nu^{u,{\bf k}}_D(\tilde\gamma) =
\lambda \nu^u_D(\tilde\gamma) \prod_{i=0}^{|\tilde\gamma|-1} \frac{4}{k_{x_i}+4} \, ,
\end{equation}
where, in the right hand side of the equation, we have chosen a representative for $\tilde\gamma$
such that $\tilde\gamma(i)=x_i$ (but the equation is independent of the choice of representative).

The expected number of loops $\tilde\gamma$ resulting from the construction of Proposition~\ref{coupling2} is
\be \label{mean2}
\lambda \nu^{u}_D(\tilde\gamma) \int_0^{\infty} e^{-t} \mathbbm{1}_{\{\sum_{i=0}^{|\tilde\gamma|-1} m^2(\tilde\gamma(i))<t\}} dt
= \lambda \nu^{u}_D(\tilde\gamma) \prod_{i=0}^{|\tilde\gamma|-1}e^{-m^2(x_i)} \, ,
\ee
where we have chosen the same representative with $\tilde\gamma(i)=x_i$ as before and the result is again
independent of the choice. Comparing equations \eqref{mean1} and \eqref{mean2}, we see that the two expected
numbers are indeed the same when $e^{-m^2(x)}=\frac{4}{k_x+4}=(1+k_x/4)^{-1}$ or $k_x=4(e^{m^2(x)}-1)$. \fbox{}

\subsection{Boundary correlations in the discrete Gaussian free field} \label{sec:gff}

In this section we discuss some interesting relations between the random walk loop soups of the
previous section and the discrete Gaussian free field. We will use the setup of the previous section,
but we need some additional notation and definitions.

Let $D$ be a bounded subset of $\mathbb C$, define $D^{\#} := D \cap {\mathbb Z}^2$ and let
$\Phi^{\bf k}_D=\{\phi_x\}_{x \in D^{\#}}$ denote a collection of mean-zero Gaussian random variables with
covariance ${\mathbb E}^{\bf k}_D(\phi_ x \phi_y) = G^{\bf k}_D(x,y)$, where $G^{\bf k}_D(x,y)$ denotes the
Green function of the random walk introduced at the beginning of Section~\ref{sec:scal-lim}, with killing rates
${\bf k}=\{k_x\}_{x \in D^{\#}}$ and killed upon exiting the domain $D$ (i.e., if the random walker attempts to
leave $D$, it is sent to the cemetery $\Delta$, where it stays forever). The lattice field $\Phi^{\bf k}_D$ is the
\emph{discrete Gaussian free field} in $D$ with zero (Dirichlet) boundary condition.
If $k_x=0 \; \forall x \in D^{\#}$, the field is called \emph{massless}, otherwise we will call it \emph{massive}.
(If the nature of the field is not specified, it means that it can be either massless or massive.)

The distribution of $\Phi^{\bf k}_D$ has density with respect to the Lebesgue measure on ${\mathbb R}^{D^{\#}}$
given by
\begin{equation} \nonumber
\frac{1}{Z^{\bf k}_D} \exp\left(-H^{\bf k}_D(\varphi)\right) \, ,
\end{equation}
where 
\begin{equation} \nonumber
Z^{\bf k}_D = \int_{{\mathbb R}^{D^{\#}}} \exp\left(-H^{\bf k}_D(\varphi)\right) \prod_{x \in D^{\#}} d\varphi_x
\end{equation}
is a normalizing constant (the \emph{partition function} of the model) and the \emph{Hamiltonian} $H^{\bf k}_D$
is defined as follows:
\begin{eqnarray*}
H^{\bf k}_D(\varphi) & := & \frac{1}{4} \sum_{x,y \in D^{\#} : x \sim y} (\varphi_y-\varphi_x)^2 + \frac{1}{2} \sum_{x \in D^{\#}} k_x \varphi^2_x
+ \frac{1}{2} \sum_{x \in \partial D^{\#}} \sum_{y \notin D^{\#} : x \sim y} \varphi^2_x \\
& = & -\frac{1}{2} \sum_{x,y \in D^{\#} : x \sim y} \varphi_x \varphi_y + \frac{1}{2} \sum_{x \in D^{\#}} (k_x+4) \varphi^2_x \, ,
\end{eqnarray*}
where the first sum is over all ordered pairs $x,y \in D^{\#}$ such that $|x-y|=1$ (denoted by $x \sim y$),
and $\partial D^{\#}$ is the set $\{x \in D^{\#} : \exists y \notin D^{\#} \text{ such that } |x-y|=1\}$.

In the first expression for $H^{\bf k}_D$, the second sum accounts for the massive nature of the field, while the
third sum accounts for the Dirichlet boundary condition. (To understand the third sum, note that one
can extend the field $\Phi^{\bf k}_D$ on $D^{\#}$ to a field $\Phi^{\bf k}=\{\phi_x\}_{x \in {\mathbb Z}^2}$
on ${\mathbb Z}^2$ by setting $\phi_x=0 \; \forall x \notin D^{\#}$.)

\begin{exercise} \label{exercise:partition-functions}
Use Lemma~1.2 of~\cite{bfsp} and Gaussian integration to show that the partition function of the
discrete Gaussian free field in $D$ can be expressed as
\be \nonumber
Z^{\bf k}_D = \left(\prod_{x \in D^{\#}}\frac{2\pi}{k_x+4}\right)^{1/2} {\cal Z}^{1/2,{\bf k}}_{D} .
\ee
\end{exercise}

The next result shows that the probability that the value of the field at a point inside the domain is affected
by a change in the shape of the domain can be computed using
the random walk loop soup with intensity $\lambda=1/2$.
\begin{proposition} \label{gff-loops}
Let $m$ be a nonnegative function (possibly identically zero), $D$ and $D'$ be bounded subsets of $\mathbb C$
containing $x_0 \in {\mathbb Z}^2$, with $D' \subset D$, and ${\bf k}$ denote the collection
$\{4(e^{m^2(x)}-1)\}_{x \in D \cap {\mathbb Z}^2}$. There exist versions of
$\Phi^{\bf k}_D = \{ \phi_x \}_{x \in D\cap {\mathbb Z}^2}$ and
$\Phi^{\bf k}_{D'} = \{ \phi'_x \}_{x \in D' \cap {\mathbb Z}^2}$, defined on the same probability space, such that
\begin{equation} \nonumber
P(\phi_{x_0} \neq \phi'_{x_0}) = P_{1/2,m} (\text{there is a loop through } x_0 \text{ that intersects } D \setminus D') \, ,
\end{equation}
where $P$ denotes the joint probability distribution of $\Phi^{\bf k}_D$ and $\Phi^{\bf k}_{D'}$, and $P_{1/2,m}$
is the law of the random walk loop soup in $D$ with intensity $\lambda=1/2$ and mass function $m$.
\end{proposition}

The proof of Proposition~\ref{gff-loops} will follow from a probabilistic coupling that allows us to define the random
walk loop soup in $D$ with intensity $1/2$ and the Gaussian free field in $D$ with Dirichlet boundary condition on
the same probability space. The coupling is given in Lemma~\ref{coupling3} below, but first we need some
additional notation.

We say that $x \in {\mathbb Z}^2$ is \emph{touched} by the unrooted loop $\tilde\gamma$, and write
$x \in \tilde\gamma$, if $\tilde\gamma(i)=x$ for some $i \in \{0,\ldots,|\tilde\gamma|-1\}$ and some
representative $\tilde\gamma(\cdot)$ of $\tilde\gamma$. If $\tilde\gamma(\cdot)$ is any rooted version
of $\tilde\gamma$, the number of indices in $\{0,\ldots,|\tilde\gamma|-1\}$ such that $\tilde\gamma(i)=x$
is denoted by $n(x,\tilde\gamma)$ (note that the notation makes sense because $n(x,\tilde\gamma)$ is
independent of the choice of representative $\tilde\gamma(\cdot)$).
To each $x \in {\mathbb Z}^2$ touched by $\tilde\gamma$, we associate $n(x,\tilde\gamma)$ independent,
exponentially distributed random variables with mean one, denoted by $\{\tau_x^i(\tilde\gamma)\}_{i=1,\ldots,n(x,\tilde\gamma)}$. 
We call the quantity
\begin{equation} \nonumber
T_x(\tilde\gamma) := \left\{ \begin{array}{ll}
\sum_{i=1}^{n(x,\tilde\gamma)} \frac{\tau_x^i(\tilde\gamma)}{k_x+4} & \mbox{ if } x \in \tilde\gamma \\
0 & \mbox{ if } x \notin \tilde\gamma
                     \end{array} \right.
\end {equation}
the \emph{occupation time} at $x$ associated to $\tilde\gamma$.
%

If $\tilde{\cal A}_{\lambda,m}$ is a realization of the random walk loop soup,
we define the \emph{occupation field} at $x$ associated to $\tilde{\cal A}_{\lambda,m}$ as
\begin{equation} \nonumber
L_x(\tilde{\cal A}_{\lambda,m}) 
:= \sum_{\tilde\gamma \in \tilde{\cal A}_{\lambda,m}} T_x(\tilde\gamma) + \frac{\tau_x^0/2}{k_x+4} \, ,
\end{equation}
where $\{\tau_x^0\}_{x \in {\mathbb Z}^2}$ is an additional collection of independent, exponential
random variables with mean one. We denote by ${\bf L}^{\bf k}_D$ the collection $\{ L_x \}_{x \in D^{\#}}$.

The next result provides the coupling needed to prove Proposition~\ref{gff-loops}. It says that, if $\{S_x\}$ are
random variables with the Ising-type distribution~\eqref{ising-type} below, where the $\{L_x\}$ are distributed
like the components of the occupation field of the random walk loop soup in $D$ with intensity $1/2$ and mass
function $m$, then the random variables $\psi_x = \sqrt{2 L_x} S_x$ are equidistributed with the components of the
discrete Gaussian free field in $D$ with $k_x = 4(e^{m^2(x)}-1)$ and Dirichlet boundary condition. More formally,
one has the following proposition, where $\hat P_{1/2,m}$ denotes the joint distribution of the random walk loop
soup in $D$ with intensity $1/2$ and mass function $m$, and the collection of all exponential random variables
needed to define the occupation field. The proof of the proposition follows easily from Theorem \ref{gff-of} in
Appendix~\ref{appendix}, which is a version of a recent result of Le Jan~\cite{lejan1} (see also \cite{lejan2}
and Theorem~4.5 of \cite{sznitman-notes}).
\begin{lemma} \label{coupling3}
Let ${\bf L}^{\bf k}_D = \{ L_x \}_{x \in D^{\#}}$, denote the occupation field of the random walk loop soup 
in $D$ with intensity $1/2$ and mass function $m$. Let ${\bf S}=\{ S_x \}_{x \in D^{\#}}$ be $(\pm 1)$-valued
random variables with (random) distribution
\begin{eqnarray}
\lefteqn{P_{{\bf L}^{\bf k}_{D}} (S_x=\sigma_x \; \forall x \in D^{\#})} \nonumber \\
& = & \frac{1}{Z} \exp\left( \sum_{x,y \in D^{\#} : x \sim y} \sqrt{L_x L_y} \sigma_x \sigma_y \right)
\, , \label{ising-type}
\end{eqnarray}
where $\sigma_x= 1$ or $-1$ and $Z$ is a normalization constant. Let $\psi_x = \sqrt{2 L_x} S_x$ for all $x \in D^{\#}$;
then under $\hat P_{1/2,m} \otimes P_{{\bf L}^{\bf k}_{D}}$, $\{ \psi_x \}_{x \in D^{\#}}$ is distributed like the
discrete Gaussian free field $\Phi^{\bf k}_D = \{ \phi_x \}_{x \in D^{\#}}$ in $D$ with $k_x = 4(e^{m^2(x)}-1)$
and Dirichlet boundary condition.
\end{lemma}

\noindent{\bf Proof.}
According to Theorem~\ref{gff-of}, $\phi^2_x$ and $\psi_x^2=2L_x$ have the same distribution for every $x \in D^{\#}$.
The lemma follows immediately from this fact and the observation that, letting $\sigma_x=\text{sgn}(\varphi_x)$,
one can write
\begin{equation} \nonumber
H^{\bf k}_D(\varphi) = -\frac{1}{2} \sum_{x,y \in D^{\#} : x \sim y} \sqrt{\varphi^2_x \varphi^2_y} \sigma_x \sigma_y
+ \frac{1}{2} \sum_{x \in D^{\#}} (k_x+4) \varphi^2_x \, . \;\; \fbox{}
\end{equation} 

\noindent{\bf Proof of Proposition~\ref{gff-loops}.}
Let $\tilde{\cal A}_{1/2,m}$ be a realization of the random walk loop soup in $D$ with intensity $\lambda=1/2$
and mass function $m$, and let $\tilde{\cal A}'_{1/2,m}$ denote the collection of loops obtained by removing from
$\tilde{\cal A}_{1/2,m}$ all loops that intersect $D \setminus D'$. $\tilde{\cal A}'_{1/2,m}$ is a random walk loop
soup in $D'$ with intensity $1/2$ and mass function $m$.

If $x_0$ is touched by a loop from $\tilde{\cal A}_{1/2,m}$ that intersects $D \setminus D'$, use $\tilde{\cal A}_{1/2,m}$
and Lemma~\ref{coupling3} to generate a collection $\{\psi_x\}_{x \in D \cap {\mathbb Z}^2}$ distributed like the
discrete Gaussian free field in $D$ with Dirichlet boundary condition, and $\tilde{\cal A}'_{1/2,m}$ and Lemma~\ref{coupling3}
to generate a collection $\{\psi'_x\}_{x \in D' \cap {\mathbb Z}^2}$ distributed like the free field in $D'$.
Then, with probability one,
\begin{equation} \nonumber
\psi^2_{x_0} = 2L_{x_0}(\tilde{\cal A}_{1/2,m}) > 2L_{x_0}(\tilde{\cal A}'_{1/2,m}) = (\psi'_{x_0})^2 \, .
\end{equation}

If $x_0$ is not touched by any loop from $\tilde{\cal A}_{1/2,m}$ that intersects $D \setminus D'$, use
$\tilde{\cal A}_{1/2,m}$ and Lemma~\ref{coupling3} to generate a $\{\psi_x\}_{x \in D \cap {\mathbb Z}^2}$
distributed like the free field in $D$. Because of symmetry, $\psi_{x_0}>0$ with probability $1/2$.
\begin{itemize}
\item If $\psi_{x_0}>0$, use $\tilde{\cal A}'_{1/2,m}$ and Lemma~\ref{coupling3} to generate a collection
$\{\psi'_x\}_{x \in D' \cap {\mathbb Z}^2}$, conditioned on $\psi'_{x_0}>0$.
\item If $\psi_{x_0}<0$, use $\tilde{\cal A}'_{1/2,m}$ and Lemma~\ref{coupling3} to generate a collection
$\{\psi'_x\}_{x \in D' \cap {\mathbb Z}^2}$, conditioned on $\psi'_{x_0}<0$.
\end{itemize}
Because of the $\pm$ symmetry of the Gaussian free field, it follows immediately that
$\{\psi'_x\}_{x \in D' \cap {\mathbb Z}^2}$ is distributed like a Gaussian free field in
$D'$, and that $\psi'_{x_0}=\psi_{x_0}$. \fbox{}

\subsection{Brownian loop soups} \label{sec:bls}


A \emph{rooted loop} $\gamma:[0,t_{\gamma}] \to {\mathbb C}$ is a continuous function with
$\gamma(0)=\gamma(t_{\gamma})$. We will consider only loops with $t_{\gamma} \in (0,\infty)$.
The \emph{Brownian bridge measure} $\mu^{br}$ is the probability measure on rooted loops of duration 1
with $\gamma(0)=0$ induced by the Brownian bridge $B_t := W_t - t W_1$, $t \in [0,1]$, where $W_t$
is standard, two-dimensional Brownian motion. A measure $\mu^{br}_{z,t}$ on loops rooted at
$z \in {\mathbb C}$ (i.e., with $\gamma(0)=z$) of duration $t$ is obtained from $\mu^{br}$ by
Brownian scaling, using the map
\begin{equation} \nonumber
(\gamma,z,t) \mapsto z + t^{1/2} \gamma(s/t) , \; s \in [0,t] \, .
\end{equation}
More precisely, we let
\begin{equation} \label{eq:brownian-bridge}
\mu^{br}_{z,t}(\cdot) := \mu^{br}(\Phi^{-1}_{z,t}(\cdot)) \, ,
\end{equation}
where
\begin{equation} \label{eq:map}
\Phi_{z,t}: \gamma(s), s \in [0,1] \mapsto z + t^{1/2} \gamma(s/t) , s \in [0,t] \, .
\end{equation}

The \emph{rooted Brownian loop measure} is defined as
\begin{equation} \label{eq:rooted-brownian-loop-measure}
\mu_r := \int_{\mathbb C} \int_0^{\infty} \frac{1}{2 \pi t^2} \, \mu^{br}_{z,t} \, dt \, d{\bf A}(z) \, ,
\end{equation}
where $\bf A$ denotes area.

The (\emph{unrooted}) \emph{Brownian loop measure} $\mu$ is obtained from the rooted one by
``forgetting the root.'' More precisely, if $\gamma$ is a rooted loop,
$\theta_u \gamma: t \mapsto \gamma(u+t \mod t_{\gamma})$ is again a rooted loop. This defines
an equivalence relation between rooted loops, whose equivalence classes we refer to as (\emph{unrooted})
\emph{loops}; $\mu(\gamma)$ is the $\mu_r$-measure of the equivalence class $\gamma$.
With a slight abuse of notation, in the rest of the paper we will use $\gamma$ to denote an unrooted loop
and $\gamma(\cdot)$ to denote any representative of the equivalence class $\gamma$.

%

The \emph{massive} (\emph{unrooted}) \emph{Brownian loop measure} $\mu^m$ is defined
by the relation
\begin{equation} \label{eq:massive-brrownian-loop-measure}
d\mu^m(\gamma) = \exp{(-R_m(\gamma))} \, d\mu(\gamma) \, ,
\end{equation}
where $m:{\mathbb C} \to {\mathbb R}$ is a nonnegative \emph{mass function} and
\begin{equation} \nonumber
R_m(\gamma) := \int_0^{t_{\gamma}} m^2(\gamma(t)) dt
\end{equation}
for any rooted loop $\gamma(t)$ in the equivalence class of the unrooted loop $\gamma$.
(Analogously, one can also define a massive \emph{rooted} Brownian loop measure:
$d\mu^m_r(\gamma) := \exp{(-R_m(\gamma))} \, d\mu_r(\gamma)$.)

If $D$ is a subset of $\mathbb C$, we let $\mu_D$ (respectively, $\mu^{m}_D$) denote $\mu$
(resp., $\mu^{m}$) restricted to loops that lie in $D$. The family of measures $\{ \mu_D \}_D$
(resp., $\{ \mu^{m}_D \}_D$), indexed by $D \subset {\mathbb C}$, satisfies the
\emph{restriction property}, i.e., if $D' \subset D$, then $\mu_{D'}$ (resp., $\mu^{m}_{D'}$) is
$\mu_D$ (resp., $\mu^{m}_D$) restricted to loops lying in $D'$.

An equivalent characterization of the Brownian loop measure $\mu$ is as follows (see~\cite{werner3}).
Given a conformal map $f:D \to D'$, let $f\circ\gamma(s)$ denote the loop $f(\gamma(t))$ in $D'$
with parametrization
\begin{equation} \label{eq:time-parametrization}
s=s(t) = \int_0^t |f'(\gamma(u))|^2 du \, .
\end{equation}
Given a subset $A$ of the space of loops in $D$, let
$f \circ A = \{ \hat\gamma = f \circ \gamma \text{ with } \gamma \in A \}$.
Up to a multiplicative constant,
$\mu_D$ is the unique measure satisfying the following two properties,
collectively known as \emph{conformal restriction property}.
\begin{itemize}
\item 
For any conformal map $f:D \to D'$,
\begin{equation} \label{eq:conf-inv}
\mu_{D'}(f \circ A) = \mu_D(A) \, .
\end{equation}
\item If $D' \subset D$, $\mu_{D'}$ is $\mu_D$ restricted to loops that stay in $D'$.
\end{itemize}
The conformal invariance of the Brownian loop measure and the result just mentioned are discussed
in Appendix \ref{BrownianLoopMeasure}.

As a consequence of the conformal invariance of $\{ \mu_D \}_D$, the family of massive
measures $\{ \mu^{m}_D \}_D$ satisfies a property called \emph{conformal covariance},
defined below.
Given a conformal map $f:D \to D'$, let $\tilde m$ be defined by the map
\begin{equation} \label{eq:mass}
m(z) \stackrel{f}{\mapsto} \tilde m(w) = \left| f'(f^{-1}(w)) \right|^{-1} m(f^{-1}(w)) \, ,
\end{equation}
where $w = f(z)$. This definition, combined with \eqref{eq:time-parametrization} and
Brownian scaling, implies that $m^2 dt = \tilde m^2 ds$. From this and \eqref{eq:conf-inv},
it follows that
\begin{equation} \label{eq:conf-cov}
\mu^{\tilde m}_{D'}(f \circ A) = \mu^m_D(A) \, ,
\end{equation}
where $A$ and $f \circ A$ have the same meaning as in equation~(\ref{eq:conf-inv}).
We call this property conformal covariance and say that the massive Brownian
loop measure $\mu^m$ is \emph{conformally covariant}.

\begin{definition} \label{def:Bls}
A \emph{Brownian loop soup} 
in $D$ with intensity $\lambda$ is a Poissonian realization from $\lambda\mu_D$. 
A \emph{massive Brownian loop soup} 
in $D$ with intensity $\lambda$ and mass function $m$ is a Poissonian realization from $\lambda\mu^{m}_D$.
\end{definition}

The (massive) Brownian loop soup ``inherits'' the property of conformal invariance (covariance) from the
(massive) Brownian loop measure. Note that in a homogeneous massive Brownian loop soup, that is, if
$m$ is constant, loops are exponentially suppressed at a rate proportional to their time duration. We will
sometimes call the conformally invariant Brownian loop soup introduced by Lawler and Werner \emph{critical},
to distinguish it from the \emph{massive} Brownian loop soup defined above.

The definition of the massive Brownian loop soup has a nice interpretation in terms of ``killed'' Brownian motion.
For a given function $f$ on the space of loops, one can write
\begin{equation} \nonumber
\int f(\gamma) e^{-R_m(\gamma)} d\mu(\gamma) 
= \int {\mathbb E}_{T_{\gamma}} f(\gamma) \mathbbm{1}_{\{R_m(\gamma)<T_{\gamma}\}} d\mu(\gamma) \, , 
\end{equation}
where ${\mathbb E}_{T_{\gamma}}$ denotes expectation with respect to the law of the mean-one,
exponential random variable $T_{\gamma}$, and $\mathbbm{1}_{\{\cdot\}}$ denotes the indicator
function. In view of this, one can think of the Brownian loop $\gamma$ under the measure $\mu^m$
defined in \eqref{eq:massive-brrownian-loop-measure} as being ``killed'' at rate $m^2(\gamma(t))$.
More precisely, one has the following alternative and useful characterization.
\begin{proposition} \label{coupling1}
A \emph{massive} Brownian loop soup in $D$ with intensity $\lambda$ and mass function $m$ can
be realized in the following way.
\begin{enumerate}
\item Take a realization of the \emph{critical} Brownian loop soup in $D$ with intensity $\lambda$.
\item Assign to each loop $\gamma$ of duration $t_{\gamma}$ an independent, mean-one,
exponential random variable, $T_{\gamma}$.
\item 
Remove from the soup all loops $\gamma$ such that
\begin{equation} \label{eq:loop-killing-1}
\int_0^{t_{\gamma}} m^2(\gamma(t)) dt > T_{\gamma} \, .
\end{equation}
\end{enumerate}

\end{proposition}
\begin{remark}
Note that Eq.~\eqref{eq:loop-killing-1} requires choosing a time parametrization
for the loop $\gamma$ but is independent of the choice.
\end{remark}

\noindent{\bf Proof.}
Let ${\cal L}^D$ denote the set of loops contained in $D$ and define
${\cal L}^D_{>\varepsilon} := \{\gamma \in {\cal L}^D : \text{diam}(\gamma)>\varepsilon \}$ and
${\cal L}^D_{>\varepsilon,r} := \{\gamma \in {\cal L}^D_{>\varepsilon} : R_m(\gamma)=r \}$.
For a subset $A$ of ${\cal L}^D_{>\varepsilon}$, let $A_r = A \cap {\cal L}^D_{>\varepsilon,r}$.
For every $\varepsilon>0$, the restriction to loops of diameter larger
than $\varepsilon$ of the massive Brownian loop soup in $D$ with mass function $m$ is a Poisson point process
on ${\cal L}^D_{>\varepsilon}$ such that the expected number of loops in $A \subset {\cal L}^D_{>\varepsilon}$
at level $\lambda>0$ is
\begin{eqnarray*}
\lambda \mu^m(A) & = & \lambda \int_A e^{-R_m(\gamma)} d\mu(\gamma) \\
& = & \lambda \int_0^{\infty} \int_{A_r} e^{-r} d\mu(\gamma) dr \\
& = & \lambda \int_0^{\infty} e^{-r} \mu(A_r) dr \, .
\end{eqnarray*}

We will now show that, when attention is restricted to loops of diameter larger than $\varepsilon$, the
construction of Proposition~\ref{coupling1} produces a Poisson point process on ${\cal L}^D_{>\varepsilon}$
with the same expected number of loops at level $\lambda>0$.

Let $N_{\lambda}(A)$ denote the number of loops in $A$ obtained from the construction of Proposition~\ref{coupling1}.
Because the Brownian loop soup is a Poisson point process and loops are removed independently, for every
$A \subset {\cal L}^D_{>\varepsilon}$ we have that
\begin{enumerate}
\item[(i)] $N_{0}(A)=0$,
\item[(ii)] $\forall \lambda,\delta>0$ and $0 \leq \ell \leq \lambda$, $N_{\lambda+\delta}(A)-N_{\lambda}(A)$ is independent of $N_{\ell}(A)$,
\item[(iii)] $\forall \lambda,\delta>0$, $\Pr(N_{\lambda+\delta}(A)-N_{\lambda}(A) \geq 2) = o(\delta)$,
\item[(iiii)] $\forall \lambda,\delta>0$, $\Pr(N_{\lambda+\delta}(A)-N_{\lambda}(A)=1) = \mu^{m}(A) \delta + o(\delta)$,
\end{enumerate}
where
(iiii) follows from the fact that, conditioned on the event $N_{\lambda+\delta}(A)-N_{\lambda}(A)=1$,
the additional point (i.e., loop) that appears going from $\lambda$ to $\lambda+\delta$ is distributed
according to the density $\frac{\mu(A_r) dr}{\mu(A)}$ on $A$. Conditions (i)-(iiii) ensure that the point
process is Poisson.

In order to identify the Poisson point process generated by the construction of Proposition~\ref{coupling1}
with the massive Brownian loop soup, it remains to compute the expected number of loops in $A$ at
level $\lambda$. For every $\varepsilon>0$ and $A \subset {\cal L}^D_{>\varepsilon}$, this is given by
\begin{equation} \nonumber
\int_0^{\infty} e^{-r} \lambda\mu(A_r) dr = \lambda\mu^m(A),
\end{equation}
which concludes the proof. \fbox{}

\subsection{Some properties of the massive Brownian loop soup} \label{sec:properties}

Let ${\cal A}(\lambda,m,D)$ denote a massive Brownian loop soup in $D \subset {\mathbb C}$
with mass function $m$ and intensity $\lambda$.
We say that two loops are \emph{adjacent} if they intersect; this adjacency relation defines \emph{clusters}
of loops, denoted by $\cal C$. (Note that clusters can be nested.) For each cluster $\cal C$, we write
$\overline{\cal C}$ for the closure of the union of all the loops in $\cal C$; furthermore, we write $\hat{\cal C}$
for the \emph{filling} of $\cal C$, i.e., the complement of the unbounded connected component of
${\mathbb C} \setminus \overline{\cal C}$.
With a slight abuse of notation, we call $\hat{\cal C}$ a \emph{cluster} and denote by $\hat{\cal C}_z$
the cluster containing $z$. We set $\hat{\cal C}_z = \emptyset$ if $z$ is not contained in any cluster $\hat{\cal C}$,
and call the set $\{ z \in D : \hat{\cal C}_z = \emptyset \}$ the \emph{carpet} (or \emph{gasket}).
(Informally, the carpet is the complement of the ``filled-in'' clusters.)

It is shown in~\cite{sw} that, in the \emph{critical} case ($m=0$), if $D$ is bounded, the set of outer boundaries
of the clusters $\hat{\cal C}$ that are not surrounded by other outer boundaries are distributed like a Conformal
Loop Ensemble in $D$, as explained in the introduction.

\begin{theorem} \label{thm:phase-trans}
Let ${\cal A}(\lambda,m,D)$ be a massive Brownian loop soup in $D$ with intensity $\lambda$ and mass function
$m$, and denote by ${\mathbb P}_{\lambda,m}$ the distribution of ${\cal A}(\lambda,m,{\mathbb C})$.
\begin{itemize}
\item If $\lambda>1$, $m$ is bounded and $D$ is bounded, with probability one the vacant set of
${\cal A}(\lambda,m,D)$ is totally disconnected.
\item If $\lambda \leq 1$ and $m$ is bounded away from zero, the vacant set of
${\cal A}(\lambda,m,{\mathbb C})$ contains a unique infinite connected component. Moreover, there is a
$\xi<\infty$ such that, for any $z \in {\mathbb C}$ and all $L>0$,
\begin{equation} \label{eq:exp-decay}
{\mathbb P}_{\lambda,m}(\text{\emph{diam}}(\hat{\cal C}_z) \geq L) \leq e^{-L/\xi} \; .
\end{equation}
\end{itemize}
\end{theorem}

Note that, although in a massive loop soup individual large loops are exponentially suppressed, 
Eq.~(\ref{eq:exp-decay}) is far from obvious, and in fact false when $\lambda>1$, since in that
equation the exponential decay refers to clusters of loops.


\begin{theorem} \label{thm:gasket}
For any bounded domain $D \subset {\mathbb C}$ and any $m: D \to {\mathbb R}$ nonnegative
and bounded, the carpet of the \emph{massive} Brownian loop soup in $D$ with mass function
$m$ and intensity $\lambda$ has the same Hausdorff dimension as the carpet of the \emph{critical}
Brownian loop soup in $D$ with the same intensity. 
\end{theorem}

It is expected that certain features of a near-critical scaling limit be the same as for the critical scaling limit.
One of these features is the Hausdorff dimension of certain geometric objects. For instance, it is proved in
\cite{nw} that the almost sure Hausdorff dimension of near-critical percolation interfaces in the scaling limit is
$7/4$, exactly as in the critical case. In view of the results in Sect.~\ref{sec:scal-lim}, Theorem~\ref{thm:gasket}
can be interpreted in the same spirit.



In the rest of the section, we present the proofs of the two theorems. To prove Theorem~\ref{thm:phase-trans},
we will use the following lemma, where, according to the notation of the theorem, ${\mathbb P}_{\lambda,m}$
denotes the distribution of the massive Brownian loop soup in $\mathbb  C$ with intensity $\lambda$ and
mass function $m$.

\begin{lemma} \label{no-loop}
Let $D \subset {\mathbb C}$ be a bounded domain with $\diam(D)>1$ and $m$ a positive
function bounded away from zero. There exist constants $c<\infty$ and
$m_0>0$, independent of $D$, such that, for every $\lambda>0$ and $\ell_0>1$,
\begin{equation} \nonumber
{\mathbb P}_{\lambda,m}(\nexists \gamma \text{ with } \diam(\gamma) > \ell_0 \text{ and } \gamma \cap D \neq \emptyset)
\geq 1 - \lambda c \, (\diam(D)+2\ell_0)^2 \, e^{-m_0 \, \ell_0} \, .
\end{equation}
\end{lemma}

\noindent{\bf Proof.}
Let $D_{\ell} := \cup_{z \in D} B_{\ell}(z)$, where $B_{\ell}(z)$ denotes the disc of radius $\ell$
centered at $z$. We define several sets of loops, namely,
\begin{eqnarray*}
{\cal L}_{\ell} & := & \{ \text{loops } \gamma \text{ with } \diam(\gamma) = \ell \} \\
{\cal L}'_{\ell} & := & \{ \text{loops } \gamma \text{ with } \diam(\gamma) = \ell \text{ and duration } t_{\gamma} \geq \ell \} \\
{\cal L}''_{\ell} & := & \{ \text{loops } \gamma \text{ with } \diam(\gamma) = \ell \text{ and duration } t_{\gamma} < \ell \} \\
{\cal L}^D_{\ell} & := & \{ \text{loops } \gamma \text{ with } \diam(\gamma) = \ell \text{ and } \gamma \cap D \neq \emptyset \} \\
{\cal L}^D_{>\ell_0} & := & \cup_{\ell>\ell_0} {\cal L}^D_{\ell}
= \{ \text{loops } \gamma \text{ with } \diam(\gamma) > \ell_0 \text{ and } \gamma \cap D \neq \emptyset \} \, .
\end{eqnarray*}

We note that
\begin{eqnarray}
{\mathbb P}_{\lambda,m}(\nexists \gamma \text{ with diam}(\gamma) > \ell_0 \text{ and } \gamma \cap D \neq \emptyset)
& = & \exp[- \lambda \mu^m({\cal L}^D_{>\ell_0})] \nonumber \\
& \geq & 1 - \lambda \mu^m({\cal L}^D_{>\ell_0}) \label{eq:lower-bound}
\end{eqnarray}
where $\mu^m$ denotes the massive Brownian loop measure with mass function $m$.
Thus, in order to prove the lemma, we look for an upper bound for $\mu^m({\cal L}^D_{>\ell_0})$.

Denoting by  $\mu^m_{D_{\ell}}$ the restriction of $\mu^m$ to $D_{\ell}$, we can write
\begin{equation} \label{eq:upper-bound}
\mu^m({\cal L}^D_{>\ell_0}) \leq \int_{\ell>\ell_0} d\mu^m_{D_{\ell}}({\cal L}_{\ell})
= \int_{\ell>\ell_0} d\mu^m_{D_{\ell}}({\cal L}'_{\ell}) + \int_{\ell>\ell_0} d\mu^m_{D_{\ell}}({\cal L}''_{\ell}) \, .
\end{equation}


From the definition of the massive Brownian loop measure (see Eqs.~\eqref{eq:massive-brrownian-loop-measure}
and~\eqref{eq:rooted-brownian-loop-measure}), we have that
\begin{eqnarray} \label{bound-first-term}
\int_{\ell>\ell_0} d\mu^m_{D_{\ell}}({\cal L}'_{\ell}) &\leq &
\int_{\ell_0}^{\infty} \frac{\pi}{4} \frac{\diam^2(D_{\ell})}{2 \pi \ell^2} \exp\left(- \ell \,\inf_{z \in D_{\ell}}m(z) \right) d\ell \nonumber \\
& < & \frac{1}{8} \int_{\ell_0}^{\infty} \diam^2(D_{\ell}) \exp\left(- \ell \, \inf_{z \in D_{\ell}}m(z) \right) d\ell \, .
\end{eqnarray}

To bound 
the second term of the r.h.s. of~\eqref{eq:upper-bound}, we observe that, if a loop $\gamma$ of duration
$t_{\gamma}$ has $\diam(\gamma) \geq \ell$, for any time parametrization of the loop, there exists a
$t_0 \in (0,t_{\gamma})$ such that $|\gamma(t_0)-\gamma(0)| \geq \ell/2$. The image $\hat\gamma$ of
$\gamma$ under $\Phi^{-1}$ must then satisfy $|\hat\gamma(t_0/t_{\gamma})| \geq \ell/(2 \sqrt{t_{\gamma}})$
(see equation~\eqref{eq:map}).

Let $W_s$ and $B_s := W_s - s W_1$ denote standard, two-dimensional Brownian motion and Brownian bridge,
respectively, with $s \in [0,1]$. Let $W^1_s$ denote standard, one-dimensional Brownian motion.
Noting that $|B_s| \leq |W_s|+|W_1| \leq 2 |W_s| + |W_1-W_s|$, and using the reflection principle and known
properties of the complementary error function (erfc), the above observation gives the bound
\begin{eqnarray} \label{brownian-bound}
\lefteqn{\mu^{br}_{z,t}(\gamma: \diam(\gamma) \geq \ell)} \nonumber \\
& \leq & \mu^{br}_{z,t}(\exists s \in (0,t) : |\gamma(s)-z| \geq \ell/2) \nonumber \\
& = & \mu^{br}\left(\exists s \in (0,1) : |\hat\gamma(s)| \geq \frac{\ell}{2\sqrt{t}}\right) \nonumber  \\
& \leq & \Pr\left(\exists s \in (0,1) : |W_s|+|W_1-W_s| \geq \frac{\ell}{4\sqrt{t}}\right) \nonumber \\
& \leq & \Pr\left(\exists s \in (0,1) : |W_s| \geq \frac{\ell}{12\sqrt{t}}\right) \nonumber \\
& \leq & 4 \Pr\left(\sup_{0 \leq s \leq 1} W^1_s \geq \frac{\ell}{12\sqrt{t}}\right) \nonumber \\
& = & 8 \Pr\left( W^1_1 \geq \frac{\ell}{12\sqrt{t}} \right) \nonumber \\
& = & 4 \, \text{erfc}\left(\frac{\ell}{12\sqrt{t}}\right) \leq 4 e^{-\ell^2/288t} \, .
\end{eqnarray}

Using \eqref{brownian-bound},
we have that, for $\ell_0>1$,
\begin{eqnarray*}
\int_{\ell>\ell_0} d\mu^m_{D_{\ell}}({\cal L}''_{\ell})
& \leq & \int_{\ell_0}^{\infty} \frac{\pi}{4} \, \diam^2(D_{\ell}) \, \int_0^{\ell} \frac{1}{2\pi t^2}
\, \mu^{br}_{0,t}(\gamma : \diam(\gamma) \geq \ell) \, dt \, d\ell \\
& \leq & \int_{\ell_0}^{\infty} \frac{\pi}{4} \, \diam^2(D_{\ell}) \, \int_0^{\ell} \frac{2}{\pi t^2} e^{-\ell^2/288t} \, dt \, d\ell \\
& \leq & \tilde c \int_{\ell_0}^{\infty} \diam^2(D_{\ell}) \, e^{-\ell/288} \, d\ell \, ,
\end{eqnarray*}
where $\tilde c < \infty$ is a suitably chosen constant that does not depend on $D$,  $\ell_0$ or $m$.

Combining this bound with the bound~\eqref{bound-first-term}, we can write
\begin{equation} \nonumber
\mu^m({\cal L}^D_{>\ell_0}) \leq \hat c \, \int_{\ell_0}^{\infty} \diam^2(D_{\ell}) \, e^{-m_0 \, \ell} \, d\ell \, ,
\end{equation}
where $m_0>0$ is any positive number smaller than $\min(1/288,\inf_{z \in {\mathbb C}} m(z))$ and
$\hat c < \infty$ is a suitably chosen constant independent of $D$, $\ell_0$ and $m$.
From this, a simple calculation leads to
\begin{equation} \nonumber
\mu^m({\cal L}^D_{>\ell_0}) \leq \hat c \int_{\ell_0}^{\infty} (\diam(D) + 2 \ell)^2 e^{-m_0 \, \ell} \, d\ell
\leq c \, (\diam(D)+2\ell_0)^2 \, e^{-m_0 \, \ell_0} \, ,
\end{equation}
where $c<\infty$ is a constant that does not depend on $D$ and $\ell_0$.
The proof is concluded using inequality~\eqref{eq:lower-bound}. \fbox \\

\medskip
\noindent {\bf Proof of Theorem~\ref{thm:phase-trans}.}
We first prove the statement in the first bullet.
Let $\overline m := \sup_{z \in {\mathbb C}} m(z)$; since $m$ is bounded, $\overline m < \infty$.
Let $\tau>0$ be so small that $\lambda' := e^{- \overline m^2 \tau} \lambda > 1$ and
denote by ${\cal A}(\lambda, 0, D)$ a critical loop soup in $D$ with intensity $\lambda$
obtained from the full-plane loop soup ${\cal A}(\lambda,0,{\mathbb C})$.
Let ${\cal A}(\lambda, m, D)$ be a massive soup obtained from ${\cal A}(\lambda, 0, D)$
via the construction of Proposition~\ref{coupling1}, and let $\overline{\cal A}(\lambda', \tau, D)$
denote the collection of loops obtained from ${\cal A}(\lambda, 0, D)$ by removing all
loops of duration $>\tau$ with probability one, and other loops with probability
$1 - e^{- \overline m^2 \tau}$, so that $\overline{\cal A}(\lambda', \tau, D)$ is a loop soup with
intensity $\lambda' = e^{- \overline m^2 \tau} \lambda$, restricted to loops of duration $\leq \tau$.
If we use the same exponential random variables to generate ${\cal A}(\lambda,m,D)$ and
$\overline{\cal A}(\lambda', \tau, D)$ from ${\cal A}(\lambda,0,D)$, the resulting loop soups are
coupled in such a way that ${\cal A}(\lambda, m, D)$ contains all the loops that are contained
in $\overline{\cal A}(\lambda', \tau, D)$, and the vacant set of ${\cal A}(\lambda, m, D)$ is a
subset of the vacant set of $\overline{\cal A}(\lambda', \tau, D)$.

We will now show that, for every $S \subset D$ at positive distance from the boundary of $D$, the
intersection with $S$ of the vacant set of $\overline{\cal A}(\lambda', \tau, D)$ is totally disconnected.
The same statement is then true for the intersection with $S$ of the vacant set of ${\cal A}(\lambda, m, D)$,
which concludes the proof of the first bullet, since the presence of a connected component larger
than one point in the vacant set of ${\cal A}(\lambda,m,D)$ would lead to a contradiction.

For a given $S \subset D$ and any $\varepsilon>0$, take $\tau=\tau(S,\varepsilon)$ so small that
$e^{- \overline m^2 \tau} \lambda>1$ and the probability that a loop from ${\cal A}(\lambda,0, {\mathbb C})$
of duration $\leq \tau$ intersects both $S$ and the complement of $D$ is less than $\varepsilon$.
(This is possible because the $\mu$-measure of the set of loops that intersect both $S$ and the
complement of $D$ is finite.)
If that event does not happen, the intersection between $S$ and the vacant set of
$\overline{\cal A}(\lambda',\tau, D)$ coincides with the intersection between $S$ and the
vacant set of the full-plane loop soup $\overline{\cal A}(\lambda',\tau, {\mathbb C})$ with cutoff
$\tau$ on the duration of loops. The latter intersection is a totally disconnected set with probability
one by an application of Theorem 2.5 of~\cite{bc2}. (Note that the result does not follow directly
from Theorem 2.5 of~\cite{bc2}, which deals with full-space soups with a cutoff on the \emph{diameter}
of loops, but can be easily obtained from it, for example with a coupling between
$\overline{\cal A}(\lambda',\tau, {\mathbb C})$ and a full-plane soup with a cutoff on the diameter of
loops chosen to be much smaller than $\sqrt{\tau}$, and using arguments along the lines of those in
the proof of Lemma~\ref{no-loop}. We leave the details to the interested reader.) Since $\tau$ can be
chosen arbitrarily small, this shows that the intersection with any $S$ of the vacant set of
$\overline{\cal A}(\lambda',\tau,D)$ is totally disconnected with probability one.

To prove the statement in the second bullet we need some definitions.
Let $R_l := [0,3 l] \times [0,l]$ and denote by $A_l$ the event that the vacant set
of a loop soup contains a crossing of $R_l$ in the long direction, i.e., that it contains a connected
component which stays in $R_l$ and intersects both $\{0\} \times [0,l]$ and $\{3 l\} \times [0,l]$.
Furthermore, let
$E^{\ell}_l := \{ \nexists \gamma \text{ with diam}(\gamma) > \ell \text{ and } \gamma \cap R_l \neq \emptyset \}$
and denote by 
$\overline{\mathbb P}_{\lambda,\ell}$ the distribution of the critical Brownian loop soup in $\mathbb C$
with intensity $\lambda$ and cutoff $\ell$ on the diameter of the loops (i.e., with all the loops of diameter
$>\ell$ removed).

We have that, for any $\ell_0>0$ and $n \in {\mathbb N}$,
\begin{eqnarray*}
{\mathbb P}_{\lambda,m}(A_{3^n}) & \geq & {\mathbb P}_{\lambda,m}(A_{3^n} \cap E^{n \, \ell_0}_{3^n}) \\
 						& = & {\mathbb P}_{\lambda,m}(A_{3^n} | E^{n \, \ell_0}_{3^n}) {\mathbb P}_{\lambda,m}(E^{n \, \ell_0}_{3^n}) \\
						& \geq & \overline{\mathbb P}_{\lambda, n \, \ell_0}(A_{3^n}) {\mathbb P}_{\lambda,m}(E^{n \, \ell_0}_{3^n}) \\
						& = & \overline{\mathbb P}_{\lambda, 3^{-n} n \, \ell_0}(A_1) {\mathbb P}_{\lambda,m}(E^{n \, \ell_0}_{3^n}) \, ,
\end{eqnarray*}
where we have used the Poissonian nature of the loop soup in the second inequality, and the last equality
follows from scale invariance.

Now consider a sequence $\{ \overline{\cal A}(\lambda, 3^{-n} n \, \ell_0, {\mathbb C}) \}_{n \geq 1}$
of full-plane soups with cutoffs $\{ 3^{-n} n \, \ell_0 \}_{n \geq 1}$, obtained from the same critical Brownian
loop soup ${\cal A}(\lambda,0,{\mathbb C})$ by removing all loops of diameter larger than the cutoff.
The soups are then coupled in such a way that their vacant sets form an increasing (in the sense of
inclusion of sets) sequence of sets. Therefore, by Kolmogorov's zero-one law,
$\lim_{n \to \infty} \overline{\mathbb P}_{\lambda, 3^{-n} n \, \ell_0}(A_1)$
is either 0 or 1. (Note that this limit can be seen as the probability of the union over $n \geq 1$ of the events that
the rectangle $R_1$ is crossed in the long direction by the vacant set of the soup with cutoff $3^{-n} \, n \, \ell_0$.)
Since $\overline{\mathbb P}_{\lambda, 3^{-n} n \, \ell_0}(A_1)$ is strictly positive for $n=1$ (see Section~3 of~\cite{bc2})
and clearly increasing in $n$,
we conclude that $\lim_{n \to \infty} \overline{\mathbb P}_{\lambda, 3^{-n} n \, \ell_0}(A_1) = 1$.

Moreover, it follows from Lemma~\ref{no-loop} that
\begin{equation} \nonumber
{\mathbb P}_{\lambda,m}(E^{n \, \ell_0}_{3^n}) \geq 1- \lambda c \, (\sqrt{10}+2\ell_0)^2 \, e^{(2 \log 3 - m_0 \ell_0) n} \, , 
\end{equation}
for some constants $c<\infty$ and $m_0>0$ independent of $\ell_0$ and $n$. Choosing $\ell_0 > (2 \log 3)/m_0$,
we have that ${\mathbb P}_{\lambda,m}(E^{n \, \ell_0}_{3^n}) \stackrel{n \to \infty}{\longrightarrow} 1$, which implies
that $\lim_{n \to \infty} {\mathbb P}_{\lambda,m}(A_{3^n}) = 1$. Note that the result does not depend on the position
and orientation of the rectangles $R_l$ chosen to define the event $A_l$.

Crossing events for the vacant set are decreasing in $\lambda$ and are therefore positively correlated (see, e.g.,
Lemma 2.2 of \cite{janson}). (An event $A$ is decreasing if ${\cal A} \notin A$ implies ${\cal A}' \notin A$ whenever
${\cal A}$ and ${\cal A}'$ are two soup realizations such that  ${\cal A}'$ contains all the loops contained in ${\cal A}$.)
Let ${\mathbb A}_n := [-3^{n+1}/2, 3^{n+1}/2] \times [-3^{n+1}/2, 3^{n+1}/2] \setminus [-3^n/2, 3^n/2] \times [-3^n/2, 3^n/2]$
and $C_n$ denote the event that a connected component of the vacant set makes a circuit inside ${\mathbb A}_n$
surrounding $[-3^n/2, 3^n/2] \times [-3^n/2, 3^n/2]$.
Using the positive correlation of crossing events, and the fact that the circuit described above can be obtained
by ``pasting'' together four crossings of rectangles, we conclude that $\lim_{n \to \infty} {\mathbb P}_{\lambda,m}(C_n) = 1$.

The existence of a unique unbounded component in the vacant set now follows from standard arguments
(see, e.g., the proof of Theorem~3.2 of~\cite{bc2}).

The exponential decay of loop soup clusters also follows immediately,
since the occurrence of $C_n$ prevents the cluster of the origin from extending
beyond the square $[-3^{n+1}/2, 3^{n+1}/2] \times [-3^{n+1}/2, 3^{n+1}/2]$. \fbox \\


\medskip
\noindent {\bf Proof of Theorem~\ref{thm:gasket}.}
Let $\overline m_D := \sup_{z \in D} m(z)$; since $m$ is bounded, $\overline m_D < \infty$.
Fix $\tau \in (0,\infty)$ and define $\lambda' := e^{- \overline m_D^2 \tau} \lambda < \lambda$.
Denote by ${\cal A}(\lambda, 0, D)$ a critical loop soup in $D$ with intensity $\lambda$.
Let ${\cal A}(\lambda, m, D)$ be a massive loop soup obtained from ${\cal A}(\lambda, 0, D)$
via the construction of Prop.~\ref{coupling1}, and let $\overline{\cal A}(\lambda', \tau, D)$
denote the set of loops obtained from ${\cal A}(\lambda, 0, D)$ by removing all loops of duration
$>\tau$ with probability one, and other loops with probability $1 - e^{- \overline m_D^2 \tau}$,
so that $\overline{\cal A}(\lambda', \tau, D)$ is a loop soup with intensity
$\lambda' = e^{- \overline m_D^2 \tau} \lambda$, restricted to loops of duration $\leq \tau$.
If we use the same exponential random variables to generate ${\cal A}(\lambda,m,D)$ and
$\overline{\cal A}(\lambda', \tau, D)$ from ${\cal A}(\lambda,0,D)$, the resulting loop soups  are
coupled in such a way that ${\cal A}(\lambda, m, D)$ contains all the loops that are contained
in $\overline{\cal A}(\lambda', \tau, D)$, and the vacant set of ${\cal A}(\lambda, m, D)$ is a
subset of the vacant set of $\overline{\cal A}(\lambda', \tau, D)$. Let $\mathbb P$ denote the
probability distribution corresponding to the coupling between soups described above.

Note that, if we denote carpets by $G[\cdot]$, we have that
\begin{equation} \label{eq:inclusion}
G[\overline{\cal A}(\lambda',\tau,D)] \supset G[{\cal A}(\lambda,m,D)] \supset G[{\cal A}(\lambda,0,D)] \, .
\end{equation}
Moreover, letting 
\begin{equation} \nonumber
h(\ell) := \frac{187 - 7 \ell + \sqrt{25 + \ell^2 - 26 \ell}}{96} \, ,
\end{equation}
denoting the Hausdorff dimension of a set $S$ by $H(S)$, and combining the computation of the
expectation dimension for carpets/gaskets of Conformal Loop Ensembles~\cite{SSW} with the results of~\cite{nawe}
(see, in particular, Section 4.5 of~\cite{nawe}), we have that, for any $\ell \in [0,1]$, with probability one,
\begin{equation} \nonumber
H(G({\cal A}(\ell,0,D))) = h(\ell) \, .
\end{equation}

We will now show that the almost sure Hausdorff dimension of $\overline{\cal A}(\lambda', \tau, D)$
equals $h(\lambda')$. To do this, consider the event $E$ that ${\cal A}(\lambda,0,D)$ contains no
loop of duration $>\tau$. Note that ${\mathbb P}(E)>0$, since the set of loops of duration $>\tau$ that
stay in $D$ has finite mass for the Brownian loop measure $\mu$. Note also that, on the event $E$,
$\overline{\cal A}(\lambda', \tau, D)$ coincides with ${\cal A}(\lambda',0,D)$. Thus, since the sets
of loops of duration $>\tau$ and $\leq \tau$ are disjoint, the Poissonian nature of the loop soups
implies that
\begin{equation} \nonumber
{\mathbb P}(H(G[\overline{\cal A}(\lambda', \tau, D)]) = h(\lambda')) =
{\mathbb P}(H(G[{\cal A}(\lambda', 0, D)]) = h(\lambda') | E) = 1 \, .
\end{equation}

From this and~\eqref{eq:inclusion}, it follows that, with probability one,
\begin{equation} \nonumber
h(\lambda') = H(G[\overline{\cal A}(\lambda', \tau, D)])
\geq H(G[{\cal A}(\lambda, m, D)])
\geq H(G[{\cal A}(\lambda, 0, D)]) = h(\lambda) \, .
\end{equation}
Since $h$ is continuous, letting $\tau \to 0$ (so that $\lambda' \to \lambda$)
concludes the proof. \fbox{}

\section{Scaling limits} \label{sec:scal-lim}

We are now going to consider scaling limits for the random walk loop soup defined in the previous section.

\subsection{The critical case} \label{sec:crit-scal-lim}

Consider a critical, full-plane, random walk loop soup
$\tilde{\cal A}_{\lambda} \equiv \tilde{\cal A}(\lambda,0,{\mathbb Z}^2)$.
Following~\cite{ltf}, for each integer $N \geq 2$, we define the \emph{rescaled random walk loop soup}
\begin{equation} \label{rescaled-soup}
\tilde{\cal A}^N_{\lambda} := \{ \tilde\Phi_N \tilde\gamma : \tilde\gamma \in \tilde{\cal A}_{\lambda}\}
\text{ with } \tilde\Phi_N \tilde\gamma(t) := N^{-1} \tilde\gamma(2N^2t) \, .
\end{equation}
$\tilde\Phi_N \tilde\gamma$ is a lattice loop of duration $t_{\tilde\gamma}:=|\tilde\gamma|/(2N^2)$
on the rescaled lattice $\frac{1}{N}{\mathbb Z}^2$ and so
$\tilde{\cal A}^N_{\lambda}$ is a 
random walk loop soup on $\frac{1}{N}{\mathbb Z}^2$, with rescaled time.

For each positive integer $N$ we also define the \emph{rescaled Brownian loop soup}
\begin{equation} \label{rescaled-soup}
{\cal A}^N_{\lambda} := \{ \Phi_N \gamma : \gamma \in {\cal A}_{\lambda}\}
\text{ with } \Phi_N \gamma(t) := N^{-1} \gamma(N^2t) \, .
\end{equation}
$\Phi_N \gamma$ is a lattice loop of duration $t_{\gamma}/(N^2)$
on the rescaled lattice $\frac{1}{N}{\mathbb Z}^2$ and so
${\cal A}^N_{\lambda}$ is a 
random walk loop soup on $\frac{1}{N}{\mathbb Z}^2$, with rescaled time.

It is shown in~\cite{ltf} that, as $N \to \infty$, $\tilde{\cal A}^N_{\lambda}$
converges to the Brownian loop soup of~\cite{lw} in an appropriate sense.
In this section we give a very brief sketch of the proof of this convergence result, which implies
that the critical random walk loop soup has a conformally invariant scaling limit (the Brownian loop soup)
and explains our use of the term \emph{critical}.

\begin{theorem} {\bf (\cite{ltf})} \label{prop:critical-soups}
There exist two sequences $\{{\cal A}^N_{\lambda}\}_{N\geq2}$ and $\{\tilde{\cal A}^N_{\lambda}\}_{N\geq2}$
of loop soups, defined on the same probability space, such that the following holds.
\begin{itemize}
\item For each $\lambda>0$, ${\cal A}^N_{\lambda}$ is a 
Brownian loop soup in $\mathbb C$ with intensity $\lambda$; the realizations
of the loop soup are increasing in $\lambda$.
\item For each $\lambda>0$, $\tilde{\cal A}^N_{\lambda}$ is a 
random walk loop soup on $\frac{1}{N}{\mathbb Z}^2$ with intensity $\lambda$ and time scaled
as in~\eqref{rescaled-soup}; the realizations of the loop soup are increasing in $\lambda$.
\item For every bounded $D \subset {\mathbb C}$, with probability going to one as $N \to \infty$,
loops from ${\cal A}^N_{\lambda}$ and $\tilde{\cal A}^N2_{\lambda}$ that are contained in $D$ and
have duration at least $N^{-1/6}$ can be put in a one-to-one correspondence with the following property. If
$\gamma \in {\cal A}^N_{\lambda}$ and $\tilde\gamma \in \tilde{\cal A}^N_{\lambda}$ are paired
in that correspondence and $t_{\gamma}$ and $t_{\tilde\gamma}$ denote their respective durations, then
\begin{eqnarray*}
|t_{\gamma} - t_{\tilde\gamma}| \leq \frac{5}{8} N^{-2} \\
\sup_{0 \leq s \leq 1} | \gamma(s t_{\gamma}) - \tilde\gamma(s t_{\tilde\gamma})| \leq c N^{-1} \log N
\end{eqnarray*}
for some constant $c<\infty$.
\end{itemize}
\end{theorem}

\noindent{\bf Sketch of the proof.} For simplicity we consider the \emph{rooted} random walk loop soup;
this suffices since one can obtain an unrooted loop soup by starting with a rooted one and forgetting the root.
Let ${\cal L}^z_n$ be the set of loops of length $2n$ rooted at $z$. If $\tilde\gamma$ has length $2n$,
then $\nu^r(\tilde\gamma) = \frac{1}{2n}(\frac{1}{4})^{2n}$, which implies that
\begin{eqnarray*}
\tilde q_n := \lambda\nu^r({\cal L}^z_n) & = & \lambda\nu^r({\cal L}^0_n) = \frac{\lambda}{2n}\left[\left(\frac{1}{4}\right)^n {2n \choose n}\right]^2 \\
& = & \frac{\lambda}{2n}\left[ \frac{1}{\pi n} - \frac{1}{4\pi n^2} + O\left(\frac{1}{n^3}\right) \right] \, ,
\end{eqnarray*}
where we have used the fact that $n! = \sqrt{2\pi} n^{n+1/2} e^{-n} \left[1 + \frac{1}{12n} + O\left(\frac{1}{n^2}\right)\right]$.
Thus, the number of loops rooted at $z$ with length $2n$ is a Poisson random variable with parameter
$\tilde q_n = \frac{\lambda}{2\pi n^2} - \frac{\lambda}{8\pi n^3} + O\left(\frac{1}{n^4}\right)$.

A realization of the random walk loop soup can be obtained by drawing a Poisson random variable $\tilde N^z_n$
with parameter $\tilde q_n =\frac{\lambda}{2\pi n^2} - \frac{\lambda}{8\pi n^3} + O\left(\frac{1}{n^4}\right)$
for each $n \in {\mathbb N}$ and each vertex $z \in {\mathbb Z}^2$, and rooting $\tilde N^z_n$ loops of length
$2n$ at $z$, with the loops chosen independently according to the uniform probability measure on ${\cal L}^z_n$.

If we now consider the \emph{rooted} Brownian loop soup with intensity measure
\begin{equation*}
\lambda\mu_r = \lambda\int_{\mathbb C} \int_0^{\infty} \frac{1}{2\pi t^2} \, \mu^{br}_{z,t} \, dt \, d{\bf A}(z) \, ,
\end{equation*}
we see that the number $N^z_n$ of rooted Brownian loops with root in the unit square centered at $z \in {\mathbb Z}^2$
with duration between $n - \frac{3}{8}$ and $n + \frac{5}{8}$ is a Poisson random variable with parameter
\begin{equation*}
q_n := \lambda\int_{n-3/8}^{n+5/8} \frac{dt}{2\pi t^2} = \frac{\lambda}{2\pi [(n+5/8)-(n-3/8)]}
= \frac{\lambda}{2\pi n^2} - \frac{\lambda}{8\pi n^3} + O\left(\frac{1}{n^4}\right) \, ,
\end{equation*}
so that $q_n - \tilde q_n = O\left(\frac{1}{n^4}\right)$.

A realization of the Brownian loop soup can be obtained by drawing a Poisson random variable $N^z_n$
with parameter $q_n =\frac{\lambda}{2\pi n^2} - \frac{\lambda}{8\pi n^3} + O\left(\frac{1}{n^4}\right)$ for
each $n \in {\mathbb N}$ and each vertex $z \in {\mathbb Z}^2$, and rooting $N^z_n$ loops of duration between
$n-3/8$ and $n+5/8$ uniformly in the unit square centered at $z$, with the loops chosen independently according
to the Brownian bridge measure of duration $T^z_n(i), i=1, \ldots, N^z_n$, where the $T^z_n(i)$'s are independent
random variables with density
\begin{equation} \nonumber
\frac{(n+\frac{5}{8})(n-\frac{3}{8})}{s^2} \, , \; n-\frac{3}{8} \leq s \leq n+\frac{5}{8} \, .
\end{equation}
Finally, in order to get a full soup, one needs to add loops of duration less than $5/8$. We will not attempt to couple
these short loops with random walk loops.

Now let $N(n,z;t)$, $n \in {\mathbb N}$, $z \in {\mathbb Z}^2$, be independent Poisson processes (in the time variable
$t$) with parameter $1$, and let $\tilde N^z_n = N(n,z;\tilde q_n)$ and $N^z_n = N(n,z;q_n)$. In order to conclude the
proof, Lawler and Trujillo Ferreras use a version of a coupling lemma due to Koml\'os, Major and Tusn\'ady \cite{kmt75}
which shows that the random walk bridge and the Brownian bridge can be coupled very closely. Using this lemma to draw
coupled versions of the random walk loops and Brownian loops needed in the constructions of the random walk and Brownian
loop soups described above, one obtains coupled versions of the soups that are very close, except possibly for small loops. \fbox{}

\subsection{The near-critical case} \label{sec:near-crit-scal-lim}

If we rescale a massive random walk loop soup with constant mass function $m>0$ in the same way as
discussed in Section~\ref{sec:crit-scal-lim},
the resulting scaling limit is trivial, in the sense that it does not contain any loops larger than one point.
This is so because, under the random walk loop measure, only loops of duration of order at least $N^2$
have diameter of order at least $N$ with non-negligible probability as $N \to \infty$, and are therefore
``macroscopic'' in the scaling limit. It is then clear that, in order to obtain a nontrivial scaling limit, the
mass function needs to be rescaled while taking the scaling limit.

Suppose, for simplicity, that the mass function $m$ is constant, and let $m_N$ denote the rescaled
mass function. When $m_N$ tends to zero, $k_x \approx 4 m_N^2$ and one has the following dichotomy.
\begin{itemize}
\item If $\lim_{N \to \infty}N m_N = 0$, loops with a number of steps of the order of $N^2$ or smaller
are not affected by the killing in the scaling limit and one recovers the critical Brownian loop soup.
\item If $\lim_{N \to \infty}N m_N = \infty$, all loops with a number of steps of the order of $N^2$ or
more are removed from the soup in the scaling limit and no ``macroscopic'' loop (larger than one point)
is left.
\end{itemize}

In view of this observation, a \emph{near-critical} scaling limit, that is, a nontrivial scaling limit that
differs from the critical one, can only be obtained if the mass function $m$ is rescaled by $O(1/N)$.
This leads us to considering the 
loop soup $\tilde{\cal A}^N_{\lambda,m}$ defined as a random walk loop soup on the rescaled lattice
$\frac{1}{N}{\mathbb Z}^2$ with mass function $m/(\sqrt{2}N)$ and rescaled time as in~\eqref{rescaled-soup}.
Such a soup can be obtained from $\tilde{\cal A}^N_{\lambda}$ using the construction in
Proposition~\ref{coupling2}, replacing $m^2(\tilde\gamma(i))$ with $\frac{1}{2N^2}m^2(\tilde\gamma(i)/N)$
in equation~(\ref{eq:loop-killing-2}).

\begin{theorem} \label{prop:massive-soups}
Let $m$ be a nonnegative function such that $m^2$ is Lipschitz continuous.
There exist two sequences $\{{\cal A}^N_{\lambda,m}\}_{N\geq2}$ and $\{\tilde{\cal A}^N_{\lambda,m}\}_{N\geq2}$
of loop soups, defined on the same probability space, such that the following holds.
\begin{itemize}
\item For each $\lambda>0$, ${\cal A}^N_{\lambda,m}$ is a 
massive Brownian loop soup in $\mathbb C$ with intensity $\lambda$ and mass $m$; the realizations
of the loop soup are increasing in $\lambda$.
\item For each $\lambda>0$, $\tilde{\cal A}^N_{\lambda,m}$ is a 
massive random walk loop soup on $\frac{1}{N}{\mathbb Z}^2$ with intensity $\lambda$, mass
$m/(\sqrt{2}N)$ and time scaled as in~\eqref{rescaled-soup}; the realizations of the loop soup are
increasing in $\lambda$.
\item For every bounded $D \subset {\mathbb C}$, with probability going to one as $N \to \infty$,
loops from ${\cal A}^N_{\lambda,m}$ and $\tilde{\cal A}^N_{\lambda,m}$ that are contained in $D$ and
have duration at least $N^{-1/6}$ can be put in a one-to-one correspondence with the following property. If
$\gamma \in {\cal A}^N_{\lambda,m}$ and $\tilde\gamma \in \tilde{\cal A}^N_{\lambda,m}$ are paired
in that correspondence and $t_{\gamma}$ and $t_{\tilde\gamma}$ denote their respective durations, then
\begin{eqnarray*}
|t_{\gamma} - t_{\tilde\gamma}| \leq \frac{5}{8} N^{-2} \\
\sup_{0 \leq s \leq 1} | \gamma(s t_{\gamma}) - \tilde\gamma(s t_{\tilde\gamma})| \leq c_1 N^{-1} \log N
\end{eqnarray*}
for some constant $c_1<\infty$.
\end{itemize}

\end{theorem}

\noindent{\bf Proof.} 
Let ${\cal A}_{\lambda}$ be a critical Brownian loop soup in $\mathbb C$ with intensity $\lambda$
and $\tilde{\cal A}_{\lambda}$ a critical random walk loop soup on ${\mathbb Z}^2$ with intensity
$\lambda$, coupled as in Theorem~\ref{prop:critical-soups}. Consider the scaled loop soups ${\cal A}^N_{\lambda}$
and $\tilde{\cal A}^N_{\lambda}$, where $\tilde{\cal A}^N_{\lambda}$ is defined in~\eqref{rescaled-soup}
and ${\cal A}^N_{\lambda} := \{\Phi_N \gamma : \gamma \in {\cal A}_{\lambda}\}$ with
$\Phi_N \gamma(t) = N^{-1} \gamma(N^2t)$ for $0 \leq t \leq t_{\gamma}/N^2$.
Note that, because of scale invariance, ${\cal A}^N_{\lambda}$ is a critical Brownian loop soup
in $\mathbb C$ with parameter $\lambda$.

It follows from Theorem~\ref{prop:critical-soups} that, if one considers only loops of duration
greater than $N^{-1/6}$, loops from ${\cal A}^N_{\lambda}$ and $\tilde{\cal A}^N_{\lambda}$ contained
in $D$ can be put in a one-to-one correspondence with the properties described in Theorem~\ref{prop:massive-soups},
except perhaps on an event of probability going to zero as $N \to \infty$. For simplicity, in the rest of the proof we will
call {\it macroscopic} the loops of duration greater than $N^{-1/6}$.

On the event that such a one-to-one correspondence between macroscopic loops in $D$ exists, we
construct the massive loop soups ${\cal A}^N_{\lambda,m}$ and $\tilde{\cal A}^N_{\lambda,m}$ in the
following way. To each pair of macroscopic loops $\gamma \in {\cal A}^N_{\lambda}$ and
$\tilde\gamma \in \tilde{\cal A}^N_{\lambda}$, paired in the correspondence of Theorem~\ref{prop:critical-soups},
we assign an independent, mean-one, exponential random variable $T_{\gamma}$. We let $t_{\gamma}$
denote the (rescaled) duration of $\gamma$ and $t_{\tilde\gamma}$ the (rescaled) duration of $\tilde\gamma$,
and let $M=2N^2 t_{\tilde\gamma}$ denote the number of steps of the lattice loop $\tilde\gamma$.
As in the constructions described in Propositions~\ref{coupling1} and \ref{coupling2}, we remove $\gamma$
from the Brownian loop soup if $\int_0^{t_{\gamma}} m^2(\gamma(s)) ds > T_{\gamma}$ and remove
$\tilde\gamma$ from the random walk loop soup if
$\frac{1}{2N^2}\sum_{k=0}^{M-1} m^2(\tilde\gamma(\frac{k}{2N^2})) > T_{\gamma}$.
The resulting loop soups, ${\cal A}^N_{\lambda,m}$ and $\tilde{\cal A}^N_{\lambda,m}$, are defined on the same
probability space and are distributed like a massive Brownian loop soup with mass function $m$ and a random
walk loop soup with mass function $m/(\sqrt{2}N)$, respectively. We use $\mathbb P$ to denote the joint distribution
of ${\cal A}^N_{\lambda,m}$, $\tilde{\cal A}^N_{\lambda,m}$ and the collection $\{T_{\gamma}\}$.

For loops that are not macroscopic, the removal of loops is done independently for the Brownian loop soup and the
random walk loop soup. If there is no one-to-one correspondence between macroscopic loops in $D$, the removal
is done independently for all loops, including the macroscopic ones.

We want to show that, on the event that there is a one-to-one correspondence between macroscopic loops in $D$,
the one-to-one correspondence survives the removal of loops described above with probability going to one as
$N \to \infty$. For that purpose, we need to compare $\int_0^{t_{\gamma}} m^2(\gamma(s)) ds$ and
$\frac{1}{2N^2}\sum_{k=0}^{M-1} m^2(\tilde\gamma(\frac{k}{2N^2}))$ for loops $\gamma$ and $\tilde\gamma$
paired in the above correspondence. In order to do that, we write
\begin{eqnarray*}
\int_0^{t_{\gamma}} m^2(\gamma(s)) ds & = & t_{\gamma} \int_0^1 m^2(\gamma(t_{\gamma}u)) du \\
& = & \lim_{n \to \infty} \frac{t_{\gamma}}{n t_{\tilde\gamma}} \sum_{i=0}^{\lfloor n t_{\tilde\gamma}\rfloor} m^2\left(\gamma\left(\frac{t_{\gamma} i}{n t_{\tilde\gamma}} \right)\right) \\
& = & \lim_{q \to \infty} \frac{t_{\gamma}/t_{\tilde\gamma}}{4qN^2} \sum_{i=0}^{2qM-1} m^2\left(\gamma\left(\frac{i}{2qM} t_{\gamma} \right)\right) \, ,
\end{eqnarray*}
where $t_{\tilde\gamma} = \frac{M}{2N^2}$ and the last expression is obtained by letting $n=4qN^2$, with $q \in {\mathbb N}$.
Thus, for fixed $N$ and $\gamma$, the quantity
\begin{equation} \nonumber
\Omega(N,q;\gamma) := \left| \int_0^{t_{\gamma}} m^2(\gamma(s)) ds - \frac{t_{\gamma}/t_{\tilde\gamma}}{4qN^2}
\sum_{i=0}^{2qM-1} m^2 \left(\gamma\left(\frac{i}{2qM} t_{\gamma} \right)\right) \right|
\end{equation}
can be made arbitrarily small by choosing $q$ sufficiently large.

Define the sets of indexes $I_0 = \{ i: 0 \leq i < q \} \cup \{ i: (2M-1)q \leq i < 2qM \}$ and
$I_k = \{ i: (2k-1)q \leq i < (2k+1)q \}$ for $1 \leq k \leq M-1$.
For $i \in I_k$, $0 \leq k \leq M-1$, we have that
\begin{eqnarray*}
\lefteqn{\left|\gamma\left(\frac{i}{2qM} t_{\gamma}\right) - \tilde\gamma\left(\frac{k}{M}t_{\tilde\gamma}\right)\right| } \\
& & \leq \left|\gamma\left(\frac{i}{2qM} t_{\gamma}\right) - \tilde\gamma\left(\frac{i}{2qM}t_{\tilde\gamma}\right)\right|
+ \left|\tilde\gamma\left(\frac{i}{2qM}t_{\tilde\gamma}\right) - \tilde\gamma\left(\frac{k}{M}t_{\tilde\gamma}\right)\right| \\
& & \leq \frac{c_1 \log N}{N} + \frac{\sqrt 2}{N}
\end{eqnarray*}
for some constant $c_1$, where the first term in the last line comes from Theorem~\ref{prop:critical-soups} and the second term
comes from the fact that $\tilde\gamma(s)$ is defined by interpolation and that
\begin{itemize}
\item if $i \in I_0$, either $0 \leq \frac{i}{2qM}t_{\tilde\gamma} < \frac{1}{2M}t_{\tilde\gamma}$ so that $\tilde\gamma(\frac{i}{2qM}t_{\tilde\gamma})$
falls on the edge of $\frac{1}{N} {\mathbb Z}^2$ between $\tilde\gamma(0)$ and $\tilde\gamma(\frac{t_{\tilde\gamma}}{M}) = \tilde\gamma(\frac{1}{2N^2})$,
or $(1 - \frac{1}{2M})t_{\tilde\gamma} \leq \frac{i}{2qM}t_{\tilde\gamma} < t_{\tilde\gamma}$
so that $\tilde\gamma(\frac{i}{2qM}t_{\tilde\gamma})$ falls on the edge between
$\tilde\gamma(t_{\tilde\gamma} - \frac{t_{\tilde\gamma}}{M}) = \tilde\gamma(t_{\tilde\gamma}-\frac{1}{2N^2})$ and $\tilde\gamma(t_{\tilde\gamma}) = \tilde\gamma(0)$,
\item if $i \in I_k$ with $1 \leq k \leq M-1$,
$\frac{k}{M}t_{\tilde\gamma} - \frac{t_{\tilde\gamma}}{2M} \leq \frac{i}{2qM}t_{\tilde\gamma} < \frac{k}{M}t_{\tilde\gamma} + \frac{t_{\tilde\gamma}}{2M}$
so that $\tilde\gamma(\frac{i}{2qM}t_{\tilde\gamma})$ falls either on the edge of $\frac{1}{N} {\mathbb Z}^2$ between
$\tilde\gamma(\frac{k-1}{M}t_{\tilde\gamma}) = \tilde\gamma(\frac{k-1}{2N^2})$ and $\tilde\gamma(\frac{k}{2N^2})$,
or on the edge between $\tilde\gamma(\frac{k}{2N^2})$ and $\tilde\gamma(\frac{k+1}{M}t_{\tilde\gamma}) = \tilde\gamma(\frac{k+1}{2N^2})$.
\end{itemize}

Since $m^2$ is Lipschitz continuous, for each $i \in I_k, \, 0 \leq k \leq M-1$, we have that
\begin{equation} \nonumber
\left|m^2\left(\gamma\left(\frac{i}{2qM}t_{\gamma}\right)\right) - m^2\left(\tilde\gamma\left(\frac{k}{M}t_{\tilde\gamma}\right)\right)\right| \leq \frac{c_2\log N}{N}
\end{equation}
for some constant $c_2<\infty$ and all $N\geq2$.
We let $\overline{m}^2_D:=\sup_{x \in D}m^2(x)$ and observe that, since $t_{\tilde\gamma} \geq N^{-1/6}$, the inequality
$|t_{\gamma}-t_{\tilde\gamma}|<\frac{5}{8}N^{-2}$ from Theorem~\ref{prop:critical-soups} implies that
$t_{\gamma}/t_{\tilde\gamma} < 1 + \frac{5}{8} N^{-11/6}$ and that $M=2N^2 t_{\tilde\gamma} < 2N^2 t_{\gamma} + \frac{5}{4}$.
It then follows that
\begin{eqnarray*} 
\lefteqn{\left| \frac{1}{2N^2} \sum_{k=0}^{M-1} m^2\left(\tilde\gamma\left(\frac{k}{2N^2}\right)\right)
- \frac{t_{\gamma}/t_{\tilde\gamma}}{4qN^2} \sum_{i=0}^{2qM-1} m^2\left(\gamma\left(\frac{i}{2qM} t_{\gamma} \right)\right) \right| } \nonumber \\
& & \leq \frac{1}{2N^2} \sum_{k=0}^{M-1} \left| m^2\left(\tilde\gamma\left(\frac{k}{2N^2}\right)\right)
- \frac{t_{\gamma}/t_{\tilde\gamma}}{2q} \sum_{i \in I_k} m^2\left(\gamma\left(\frac{i}{2qM} t_{\gamma} \right)\right) \right| \nonumber \\
& & \leq \frac{1}{2N^2} \sum_{k=0}^{M-1} \frac{1}{2q} \sum_{i \in I_k} \left| m^2\left(\tilde\gamma\left(\frac{k}{2N^2}\right)\right)
- \frac{t_{\gamma}}{t_{\tilde\gamma}} m^2\left(\gamma\left(\frac{i}{2qM} t_{\gamma} \right)\right) \right| \nonumber \\
& & \leq  \frac{|1-t_{\gamma}/t_{\tilde\gamma}|}{2N^2} \sum_{k=0}^{M-1} m^2\left(\tilde\gamma\left(\frac{k}{2N^2}\right)\right) \nonumber \\
& & + \frac{t_{\gamma}/t_{\tilde\gamma}}{4qN^2} \sum_{k=0}^{M-1} \sum_{i \in I_k} \left| m^2\left(\tilde\gamma\left(\frac{k}{2N^2}\right)\right)
- m^2\left(\gamma\left(\frac{i}{2qM} t_{\gamma}\right)\right) \right| \nonumber \\
& & \leq  \frac{|1-t_{\gamma}/t_{\tilde\gamma}|}{2N^2} \sum_{k=0}^{M-1} m^2\left(\tilde\gamma\left(\frac{k}{2N^2}\right)\right) + t_{\gamma} \frac{c_2\log N}{N} \nonumber \\
& & \leq  \frac{|1-t_{\gamma}/t_{\tilde\gamma}|}{2N^2} M \overline{m}^2_D + t_{\gamma} \frac{c_2\log N}{N}
< t_{\gamma} \frac{c'_3\log N}{N} \, ,
\end{eqnarray*}
for some positive constant $c'_3=c'_3(D,m)<\infty$ independent of $\gamma$ and $\tilde\gamma$.

Therefore, for fixed $N$ and pair of macroscopic loops, $\gamma$ and $\tilde\gamma$, and for any $q \in {\mathbb N}$,
\begin{equation} \nonumber
\left| \int_0^{t_{\gamma}} m^2(\gamma(s)) ds -
\frac{1}{2N^2}\sum_{k=0}^{M-1} m^2\left(\tilde\gamma\left(\frac{k}{2N^2}\right)\right) \right|
< t_{\gamma} \frac{c'_3\log N}{N} + \Omega(N,q;\gamma) \, .
\end{equation}
For fixed $N$ and $\gamma$, one can choose $q^*$ so large that
\begin{equation} \nonumber
\Omega(N,q^*;\gamma) < t_{\gamma} \frac{c'_3 \log N}{N} \, .
\end{equation}
Hence, there is a positive constant $c_3=2c'_3$ such that, for every $N \geq 2$ and every pair of
macroscopic loops, $\gamma$ and $\tilde\gamma$, paired in the correspondence of Theorem~\ref{prop:critical-soups},
\begin{equation} \nonumber
\left| \int_0^{t_{\gamma}} m^2(\gamma(s)) ds -
\frac{1}{2N^2}\sum_{k=0}^{M-1} m^2\left(\tilde\gamma\left(\frac{k}{2N^2}\right)\right) \right|
< t_{\gamma} \frac{c_3\log N}{N} \, .
\end{equation}

We now need to estimate the number of macroscopic loops contained in $D$. For that purpose,
we note that, using the rooted Brownian loop measure~\eqref{eq:rooted-brownian-loop-measure},
the mean number, $\cal M$, of macroscopic loops contained in $D$ can be bounded above by
\begin{equation} \label{eq:mean}
{\cal M} = \lambda \int_D \int_{N^{-1/6}}^{\infty} \frac{1}{2\pi t^2} \, \mu^{br}_{z,t}(\gamma: \gamma \subset D) \, dt \, dA(z) \leq \frac{\lambda \, \diam^2(D)}{8} N^{1/6} \, .
\end{equation}

Let ${\cal A}^N(\lambda,m;D)$ (respectively, $\tilde{\cal A}^N(\lambda,m;D)$) denote the massive Brownian
(resp., random walk) loop soup in $D$, i.e., the set of loops from ${\cal A}^N_{\lambda,m}$
(respectively, $\tilde{\cal A}^N_{\lambda,m}$) contained in $D$. For the critical soups, we use the
same notation omitting the $m$.

Let $A_N$ denote the event that there is a one-to-one correspondence between macroscopic loops from
${\cal A}^N(\lambda;D)$ and $\tilde{\cal A}^N(\lambda;D)$, and let $A^m_N$ denote the event that
there is a one-to-one correspondence between macroscopic loops from ${\cal A}^N(\lambda,m;D)$ and
$\tilde{\cal A}^N(\lambda,m;D)$. Furthermore, we denote by $X$ the number of macroscopic loops in
${\cal A}^N(\lambda;D)$, and by $T$ a mean-one exponential random variable. We have that, for any
$c_4,\theta>0$ and for all $N$ sufficiently large,

\begin{eqnarray*}
{\mathbb P}(A^m_N) & \geq & {\mathbb P}(A^m_N \cap A_N \cap \{X \leq c_4 N^{1/6}\}
\cap \{ \nexists \gamma \in {\cal A}^N(\lambda;D) : t_{\gamma} \geq \theta \}) \\
& = & {\mathbb P}(A^m_N | A_N \cap \{X \leq c_4 N^{1/6}\}
\cap \{ \nexists \gamma \in {\cal A}^N(\lambda;D) : t_{\gamma} \geq \theta \}) \\
& & {\mathbb P}(A_N \cap \{X \leq c_4 N^{1/6}\}
\cap \{ \nexists \gamma \in {\cal A}^N(\lambda;D) : t_{\gamma} \geq \theta \}) \\
& \geq & \left[ 1 - \sup_{x \geq 0} \Pr\left(x \leq T \leq x + \frac{c_3 \theta \log N}{N}\right)
\right]^{c_4 N^{1/6}} \\
& & {\mathbb P}(A_N \cap \{X \leq c_4 N^{1/6}\}
\cap \{ \nexists \gamma \in {\cal A}^N(\lambda;D) : t_{\gamma} \geq \theta \}) \\
& = & \exp\left(-\frac{c_5 \theta \log N}{N^{5/6}}\right) \\
& & {\mathbb P}(A_N \cap \{X \leq c_4 N^{1/6}\}
\cap \{ \nexists \gamma \in {\cal A}^N(\lambda;D) : t_{\gamma} \geq \theta \}) \, ,
\end{eqnarray*}
where $c_5=c_3 c_4$.

Since $\exp\left(-\frac{c_5 \theta \log N}{N^{5/6}}\right) \to 1$ as $N \to \infty$ for any $c_5,\theta>0$,
in order to conclude the proof, it suffices to show that ${\mathbb P}(A_N \cap \{ X \leq c_4 N^{1/6} \} \cap
\{ \nexists \gamma \in {\cal A}^N(\lambda;D) : t_{\gamma} \geq \theta \})$ can be made arbitrarily close
to one for some choice of $c_4$ and $\theta$, and $N$ sufficiently large.
But by Theorem~\ref{prop:critical-soups}, ${\mathbb P}(A_N) \geq 1 - c (\lambda+1) \diam^2(D) N^{-7/2} \to 1$ as $N \to \infty$;
moreover, if $c_4>\frac{\lambda \, \diam^2(D)}{8}$, by equation~\eqref{eq:mean}, $c_4 N^{1/6}$ is larger than
the mean number of macroscopic loops in $D$. Since $X$ is a Poisson random variable with parameter equal to
the mean number $\cal M$ of macroscopic loops in $D$, the latter fact (together with a Chernoff bound argument)
implies that
\begin{equation} \nonumber
{\mathbb P}(X > c_4 N^{1/6}) \leq \frac{e^{-{\cal M}} (e {\cal M})^{c_4 N^{1/6}}}{(c_4 N^{1/6})^{c_4 N^{1/6}}}
\leq \left(\frac{e \, \lambda \, \diam^2(D)}{8 c_4}\right)^{c_4 N^{1/6}} \, .
\end{equation}
This shows that, if $c_4 > e \, \lambda \, \diam^2(D)/8$, ${\mathbb P}(X \leq c_4 N^{1/6}) \to 1$ as $N \to \infty$.

To find a lower bound for ${\mathbb P}(\nexists \gamma \in {\cal A}^N(\lambda;D) : t_{\gamma} \geq \theta)$, we define
\begin{equation} \nonumber
{\cal L}_{\theta,D} := \{ \text{loops } \gamma \text{ with } t_{\gamma} \geq \theta \text{ that stay in } D \} \, .
\end{equation}
We then have
\begin{eqnarray} \label{lower-bound}
{\mathbb P}(\nexists \gamma \in {\cal A}^N(\lambda;D) : t_{\gamma} \geq \theta) & = &
\exp\left[ -\lambda \mu_D({\cal L}_{\theta,D}) \right] \nonumber \\
& \geq & 1 - \lambda \mu_D({\cal L}_{\theta,D}) \nonumber \\
& \geq & 1 - \frac{\lambda \, \diam^2(D)}{\theta} \, ,
\end{eqnarray}
where the last line follows from the bound
\begin{equation} \nonumber
\mu_D({\cal L}_{\theta,D})
= \int_D \int_{\theta}^{\infty} \frac{1}{2 \pi t^2} \, \mu^{br}_{z,t}(\gamma : \gamma \text{ stays in } D) \, dt \, d{\bf A}(z)
\leq \frac{\diam^2(D)}{\theta} \, .
\end{equation}
The lower bound~\eqref{lower-bound}, together with the previous observations, shows that ${\mathbb P}(A_N^m)$
can be made arbitrarily close to one by choosing $c_4 > e \, \lambda \, \diam^2(D)/8$, $\theta$ sufficiently large,
depending on $D$, and then $N$ sufficiently large, depending on the values of $c_4$ and $\theta$. \fbox{}

\section{Conformal correlation functions in the Brownian loop soup}

In the rest of these lecture notes, we will focus on the \emph{critical} (\emph{massless}) Brownian loop soup and analyze certain
correlation functions that characterize aspects of its distribution. (One could perform a similar analysis using the \emph{massive}
Brownian loop soup, but we will not do that in these notes.)
We will show that these correlation functions need to be defined via a regularization procedure that requires introducing a cutoff,
which is then removed via a limiting procedure. The limiting procedure produces a nontrivial limit only if the correlation functions
with a cutoff are scaled by an appropriate power of the cutoff. In that case, the limiting functions have multiplicative scaling behavior
under conformal transformations, a type of behavior that characterizes correlation functions of \emph{primary fields} in conformal
field theory. The models described in this section are interesting in their own right and are studied in~\cite{cgk15}. (One of them
was inspired by~\cite{fk}.) Here, we use them as ``toy'' models for the behavior of lattice correlation functions in the scaling limit.

\subsection{Motivation and summary of results} \label{sum}

The study of lattice models involves local lattice ``observables.'' A typical example of such local observables are the spin variables in the Ising model, or,
more generally, sums of (products of) spin variables. At the critical point, the scaling limit of such observables can lead to conformal fields, which are not
defined pointwise (they are not functions), but can be defined as generalized functions (distributions in the sense of Schwartz). For example, the scaling limit
of the critical Ising magnetization in two dimensions leads to a conformally covariant magnetization field \cite{cgn1}.

\emph{Correlation functions} (or \emph{$n$-point functions}), that is, expectations of products of local lattice observables, are important tools in analyzing
the behavior of both lattice models and field theories. In general, lattice correlation functions do not possess a finite, nontrivial scaling limit. However, the
theory of \emph{renormalization} (based on the study of exactly solved models, as well as the perturbative analysis of cutoff quantum field theories) suggests
that certain correlation functions are \emph{multiplicatively renormalizable}; that is, they possess a scaling limit when multiplied by appropriate powers of the
lattice spacing. When that happens, the limiting functions are expected to have multiplicative scaling behavior under conformal transformations.

From the point of view of field theory, a lattice can be seen as a way to avoid divergencies, with the lattice spacing playing the role of an \emph{ultraviolet cutoff},
preventing infinities due to small scales. Below, we will define and study correlation functions from the Brownian loop soup. We will see that, in the absence of a
lattice, we will need to impose a different ultraviolet cutoff $\delta>0$. The behavior of the correlation functions described below is conceptually the same
as that of lattice correlation functions, with the cutoff $\delta$ playing the role of the lattice spacing. We will analyze the correlation functions for two different
types of ``fields,'' which we now describe informally.

\begin{figure}
\centering
\includegraphics[scale=0.5]{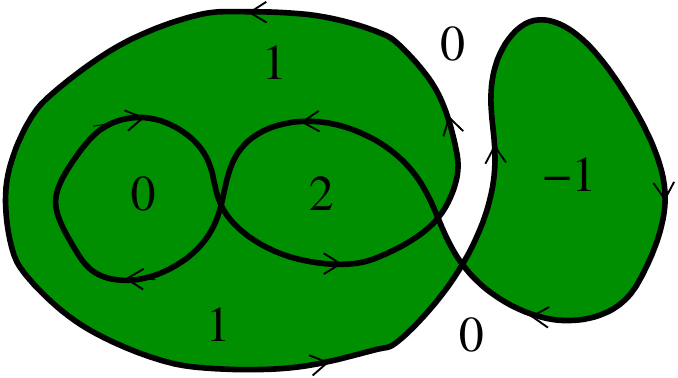}
\caption{A stylized Brownian loop. The numbers indicate the winding numbers of the loop that contribute additively to $N_{w}$,
while the green shaded region is the interior of the loop (the set of points disconnected from infinity by the loop) that contributes
$\pm 1$ (where the sign is a Boolean variable assigned randomly to each loop) to the layering number $N_\ell$.}
\label{winds}
\end{figure}

\paragraph{Winding operator.} Let $N_{w}(z)$ count the total winding number of loops around a point $z$.
We call $N_{w}$ the \emph{winding operator}.
Since the loops in a loop soup have an orientation, winding numbers can be positive or negative (see \figref{winds}).
We are interested in the $n$-point correlation functions of the ``field'' $N_w$. However, since the Brownian loop soup
is scale invariant, any given $z$ is almost surely ``surrounded'' by infinitely many loops and $N_{w}(z)$ is almost surely
\emph{infinite}, so the ``field'' $N_w$ is not well defined. We will need to restrict the soup to a bounded domain and to
introduce a cutoff $\delta>0$ and count only the winding of loops of diameter at least $\delta$.

\begin{figure}
\centering
\includegraphics[scale=.5]{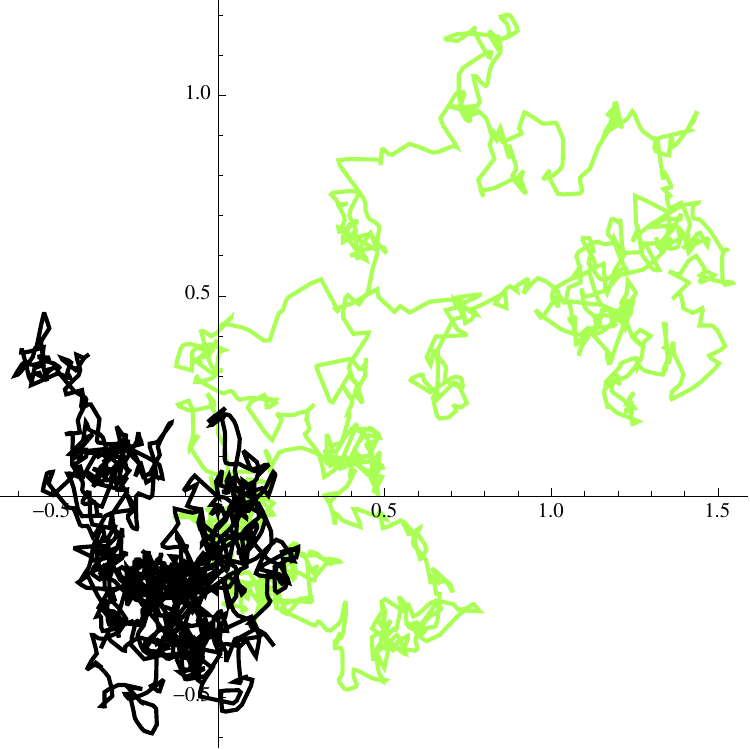}
\caption{Two simulated Brownian loops of time length $t=1$.}
\label{loops}
\end{figure}

\paragraph{Layering operator.} For each loop, define the interior as the set of points inside the outermost edge of the loop.
Points in this set are ``covered'' by the ``filled-in'' loop. Furthermore, declare each loop to be of type 1 or type 2 with equal probability,
independently of all other loops. For a point $z$, let $N_j(z)$ be the number of distinct loops of type $j=1,2$ that cover $z$, and let
$N_{\ell}(z) = N_1(z)- N_2(z)$. We'll call $N_{\ell}(z)$ the ``layering number'' at $z$, and $N_{\ell}$ the \emph{layering operator}.
Because the Brownian loop soup is scale invariant, any given point of the plane is covered by infinitely many loops with probability one,
so the layering ``field'' $N_{\ell}$. is not well defined. As before, we will need to restrict to bounded domains, as well as introduce a
cutoff $\delta>0$ and count only the contribution from loops of diameter at least $\delta$. \\

The first model described above appears more natural; the advantage of the second model is that it is simpler to analyze and it can be
easily defined on the full plane, as we will see later. Both models present logarithmic divergences (like massless fields in two dimensions)
as one tries to remove the cutoff by letting $\delta \to 0$. For this reason, we will focus on exponentials of the winding and layering
numbers times imaginary coefficients. We will denote these \emph{exponential operators} as $$V_{\beta}(z) = e^{i \beta N(z)},$$
where $N$ can be either a layering or winding number, and will prove the following results.

\begin{itemize}

\item For both models in finite domains $D$, \emph{correlators} of $n \in {\mathbb N}$ exponential operators 
\be \label{cor} 
\left\la \Pi_{j=1}^n V_{\beta_{j}}(z_{j}) \right\ra_{\delta,D} = \left\la e^{i \sum_{j=1}^n \beta_{j} N(z_{j}) } \right\ra_{\delta,D} 
\ee
exist as long as a cutoff $\delta>0$ on the the diameters of the loops is imposed, and
\begin{equation} \nonumber
\lim_{\delta \to 0} 
\frac{\left\langle \prod_{j=1}^n V_{\beta_j}(z_{j})\right\rangle_{\delta, D}}{\prod_{j=1}^n \delta^{2\Delta(\beta_j)}} =:
\phi_D(z_1, \dots, z_n; \beta_1, \dots, \beta_n) \equiv \phi_D (\boldsymbol{z}; \boldsymbol{\beta}) 
\end{equation}
exists and is finite. Moreover, if $D'$ is another finite domain and $f:D \to D'$ is a conformal map such that
$z'_1=f(z_1),\ldots,z'_n=f(z_n)$, then
\begin{equation} \nonumber
\phi_{D'}(\boldsymbol{z'};\boldsymbol{\beta}) =
\prod_{j=1}^n \left| f'(z_j) \right|^{-2\Delta(\beta_j)} \phi_D(\bs{z};\boldsymbol{\beta}) \, ,
\end{equation}
where the $\Delta(\beta)$ is defined below. This is the behavior expected for a conformal primary operator/field.
(In the conformal field theory literature, the terms field and operator are sometimes used interchangeably.)

\item For both versions in infinite volume, correlators of $n$ exponential operators 
\be \label{cor} 
\left\la \Pi_{j} V_{\beta_{j}}(z_{j}) \right\ra_{\delta} = \left\la e^{i \sum_{j} \beta_{j} N(z_{j}) } \right\ra_{\delta} 
\ee
vanish.
However, in the case of the layering model, one can let $\delta \to 0$ and still obtain a nontrivial limit by imposing
the following condition, satisfied mod $2 \pi$,
\be \label{cc}
\sum_{j} \beta_{j} = 2 \pi k, \,  \,  \,  k \in {\mathbb Z}.
\ee

\item The correlators \eqref{cor} of layering model in the plane are finite and non-zero when \eqref{cc}
is satisfied, so long as the loop soup is cut off at short scales (no other cutoff is necessary). 

\item In the case of 2 points, assuming \eqref{cc}, the $\delta \to 0$ limit of the renormalized correlators of
the layering model in the plane \eqref{cor} can be explicitly computed up to an overall multiplicative constant.
The result is
\be \nonumber
 \phi_{\mathbb C}(z_1,z_2;\beta_1,\beta_2) = C_{2} \left| \left( {1 \over z_1
     - z_2} \right)^{ \Delta(\beta_1)+\Delta(\beta_2)}  \right|^{2},
 \ee
where $C_{2}$ is a constant (see Section~\ref{2-point-function-layering-model}).

\item The $\delta \to 0$ limit of the renormalized 3-point function for the layering model in the plane,
assuming \eqref{cc}, is
\begin{eqnarray*}
\lefteqn{ \phi_{\mathbb C}(z_1,z_2,z_3;\beta_1,\beta_2,\beta_3) =} \\
& & C_{3} \left| \left( {1 \over |z_1 - z_2|} \right)^{\Delta(\beta_1) + \Delta(\beta_2) - \Delta(\beta_3)}  \left( {1 \over |z_1 - z_3|} \right)^{\Delta(\beta_1) + \Delta(\beta_3) - \Delta(\beta_2)} \left( {1 \over |z_2 - z_3|} \right)^{\Delta(\beta_2) + \Delta(\beta_3) - \Delta(\beta_1)} \right|^2 ,
\end{eqnarray*}
where $C_{3}$ is a constant (Corollary~\ref{3-point-function}).

\item The exponents $\Delta(\beta)$ are called \emph{conformal dimensions} and differ for the two models.
For the layering model,
$$
\Delta_{l}(\beta) = {\lambda \over 10} (1-\cos \beta).
$$
For the winding model, 
$$
\Delta_{w}(\beta) = \lambda \beta (2 \pi-\beta)/8 \pi^2,
$$
where this formula applies for $0\leq \beta < 2 \pi$, and $\Delta_{w}(\beta)$ is periodic under $\beta \to \beta+2 \pi$.  

\end{itemize}

\subsection{Correlation functions of the layering and winding models}

\subsubsection{Correlation functions of layering model}

As already discussed, the layering model is  defined by randomly assigning a Boolean variable to each loop in the Brownian loop soup.
Alternatively, one can think of this as two independent Brownian loop soups, each with a Poisson distribution $P_{\lambda_{+(-)}, \mu}$
with intensity measure $\lambda_{+(-)}\mu$, where we take $\lambda_+ = \lambda_- = \lambda/2$.
(This follows from the fact that the collection of all loops from a Brownian loop soup of intensity $\lambda_{+}$ and an independent
one of intensity $\lambda_{-}$ is distributed like a Brownian loop soup with intensity $\lambda_{+} + \lambda_{-}$.)

Denote by $N_{+(-)}(z)$ the number of loops $\gamma$ in the first (respectively, second) class such that the
point $z \in \mathbb C$ is separated from infinity by the image of $\gamma$ in $\mathbb C$.  If $\bar \gamma$ is the ``filled in"
loop $\gamma$, then this condition becomes $z \in \bar \gamma$, or $z$ is covered by $\gamma$. We are interested in the layering field $N_{\ell}$,
with $N_{\ell}(z) = N_{+}(z)-N_{-}(z)$. This is purely formal as  both $N_{+(-)}(z)$ are infinite with probability one for any $z$. They are  infinite
for two reasons: both because there are infinitely many large loops surrounding $z$ (infrared, or IR, divergence), and because there are 
infinitely many small loops around $z$ (ultraviolet, or UV, divergence). 

We will consider correlators of the exponential operator $V_{\beta} = e^{i \beta N_{\ell}(z)}$, and show that there are
choices of $\beta$ that remove the IR divergence and a normalization which removes the UV divergence. Specifically, we are interested in 
the correlators $V_{\boldsymbol{\beta}}(z_1, \dots, z_n) := \prod_{j=1}^n V_{\beta_j}(z_j) = e^{i \sum_{j=1}^n \beta_j N_{\ell}(z_j)}$
and their moments
\be \nonumber
\left\la V_{\boldsymbol{\beta}}(z_1, \dots, z_n) \right\ra := {\mathbb E}_{\lambda}(V_{\boldsymbol{\beta}}(z_1, \dots, z_n))
\ee
where $z_j \in \mathbb C$, $\boldsymbol {\beta}=(\beta_1, \dots, \beta_n) \in \mathbb R^n$, and the expected value
${\mathbb E}_{\lambda}$ is  taken with respect to two independent copies of the Brownian Loop Soup with equal
intensities $\lambda/2$. 

As the field $N_{\ell}$ has both IR and UV divergences, we get a meaningful definition by introducing cutoffs which
restrict the loops to have diameter
within some $\delta$ and $R \in \mathbb R^+$, $\delta < R$: let
$\mu_{\delta, R}(\cdot)=\mu(\cdot \cap \{\gamma:\delta \leq \diam(\gamma) < R\})$
and consider the correlators
\be \nonumber
\left\la V_{\boldsymbol{\beta}}(z_1, \dots, z_n) \right\ra_{\delta, R} := {\mathbb E}_{\lambda,\delta,R}\left(e^{i \sum_{j=1}^n \beta_j N_{\ell}(z_j)}\right),
\ee
where the expectation ${\mathbb E}_{\lambda,\delta,R}$ is with respect to the distribution $P_{\lambda,\mu_{\delta,R}} \otimes P_{1/2}$,
where $P_{\lambda,\mu_{\delta,R}}$ is the Poisson distribution with intensity measure $\lambda\mu_{\delta,R}$ and $P_{1/2}$  is the distribution
of a countable sequence of independent Bernoulli random variables with parameter $1/2$ (remember that each loop belongs to one of two classes
with equal probability). It is for these correlators that a suitable choice of $\boldsymbol {\beta}$ and a suitable normalization will allow us to remove
both IR and UV cutoffs.

\subsubsection{The $1$-point function in the layering model}
In this section we explicitly compute the $1$-point function in the presence of IR and UV cutoffs.
Replacing the area of a filled Brownian loop of time length $1$ with that of a disk of radius $1$,
the result reproduces the $1$-point function in the disk model of~\cite{fk}.
\begin{lemma}
For all  $z \in \mathbb C$, we have that
\be \nonumber
\left\la V_{\beta}(z) \right\ra_{\delta, R} = \left(\frac{R}{\delta}\right)^{- \frac{\lambda}{5}(1-\cos\beta)}.
\ee
\end{lemma}

\noindent{\bf Proof.}
With IR and UV cutoffs in place, the field $N_{\ell}(z)$ can be realized as follows. Let $\eta$ be a realization of loops,
and let $\{X_{\gamma}\}_{\gamma \in \eta}$ be a collection of Bernoulli symmetric random variables taking values
in $\{-1,1\}$. The quantity
\be \nonumber
N_\ell(z) = \sum_{\gamma \in \eta, z \in \bar \gamma, \delta \leq \diam(\gamma) <R} X_{\gamma} =: {\sum}^* X_{\gamma}
\ee
is finite $P_{\lambda, \mu}$ almost surely, since 
$\mu\{\gamma: z \in \bar \gamma, \delta \leq \diam(\gamma) <R\}
=\mu_{\delta, R}\{\gamma: z \in \bar \gamma \}  < \infty$
(see Appendix~\ref{BrownianLoopMeasure} and \cite{werner3}). Now,
\begin{eqnarray*}
\la V_{\beta}(z) \ra_{\delta, R} &=& {\mathbb E}_{\lambda,\delta, R}\left(e^{i \beta N_{\ell}(z)}\right) \\
&=& \sum_{k=0}^{\infty} {\mathbb E}_{\lambda,\delta, R}\left(e^{i \beta N_{\ell}(z)} | \mathcal L_k \right) 
P_{\lambda, \mu_{\delta, R}}(\mathcal L_k),
\end{eqnarray*}
where  $\mathcal L_k= \{\eta: |\{\gamma \in \eta: z \in \bar \gamma,\delta \leq \diam(\gamma) <R\}|=k \}$.
If $X$ denotes a $(\pm 1)$-valued symmetric random variable,
\be \nonumber
{\mathbb E}_{\lambda,\delta, R}\left(e^{i \beta \sum^* X_{\gamma}} | \mathcal L_k \right) =
\left(E\left(e^{i \beta X}\right)\right)^k = (\cos\beta)^k .
\ee

Therefore, for $\alpha_{z, \delta, R}= \mu_{\delta, R} (\gamma: z \in \bar \gamma)$, we have that
\begin{eqnarray*}
\la V_{\beta}(z) \ra_{\delta, R} &=& \sum_{k=0}^{\infty} (\cos\beta)^k
\frac{(\lambda \alpha_{z, \delta, R})^k}{k!} e^{-\lambda \alpha_{z, \delta, R}}
 \\
&=& e^{-\lambda \alpha_{z, \delta, R}(1-\cos\beta)}.
\end{eqnarray*}
Moreover, by Lemma \ref {FirstLemma},
\be \nonumber
\alpha_{z, \delta, R} = \frac{1}{5} \log\frac{R}{\delta} ,
\ee
which implies
\be \nonumber
\la V_{\beta}(z) \ra_{\delta, R} = \left(\frac{R}{\delta}\right)^{- \frac{\lambda}{5}(1-\cos(\beta))} ,
\ee
as claimed. \fbox{}

\subsubsection{The $1$-point function in the winding model} \label{Sec:WindingOperator}

To define the second model, let $N_{w}(z)$ denote the total winding number  about the  point $z$ of all
loops in a Brownian loop soup; as for the layering operators, this is a formal definition as, in general,
$N_{w}(z)$ is infinite. Consider again the correlators
$V_{\boldsymbol{\beta}}(z_1, \dots, z_n) = e^{i \sum_{j=1}^n \beta_j N_{w}(z_j)}$ and their moments
$
\la V_{\boldsymbol{\beta}}(z_1, \dots, z_n) \ra = E(V_{\boldsymbol{\beta}}(z_1, \dots, z_n))
$
where $z_j \in \mathbb C$, $\boldsymbol {\beta}=(\beta_1, \dots, \beta_n) \in \mathbb R^n$, and the expected
value is taken with respect to the Brownian loop soup distribution. Denoting by $P_{\lambda, \mu}$ the Poisson distribution
with intensity measure $\lambda\mu$, and restricting the loops to have diameter between some
$\delta$ and $R \in \mathbb R^+$, with $\delta < R$, we let
$\mu_{\delta, R}(\cdot)=\mu(\cdot \cap \{\gamma:\delta \leq \diam(\gamma) < R\})$
and consider the correlators
\be \nonumber
\la V_{\boldsymbol{\beta}}(z_1, \dots, z_n) \ra_{\delta, R} := {\mathbb E}_{\lambda,\delta, R}\left(e^{i \sum_{j=1}^n \beta_j N_{w}(z_j)}\right).
\ee

We now explicitly compute the $1$-point function in the presence of IR and UV cutoffs
for the winding model, using the following result, which we leave as an exercise.

\begin{exercise} \label{exercise:sum}
Show that
\be \nonumber
\sum_{m=1}^{\infty}\frac{1}{m^2} \left(1-\cos(m \beta)\right) = \frac{1}{4}\beta(2 \pi - \beta) ,
\ee
where, on the right hand side of the equation, $\beta$ is to be interpreted modulo $2\pi$. (Hint: the left hand
side is a convergent series representing a periodic function and can therefore be written as a Fourier series.)
\end{exercise}

\begin{lemma}
For all  $z \in \mathbb C$, we have that
\be \label{winddim}
\la V_{\beta}(z) \ra_{\delta, R} = \left(\frac{R}{\delta}\right)^{- \lambda \frac{\beta(2 \pi - \beta)}{4 \pi^2} },
\ee
where the formula is valid for $\beta \in [0,2\pi)$, and for $\beta \not\in [0,2\pi)$, in the right hand side, $\beta$
should be replaced by $(\beta \mod 2\pi)$.
\end{lemma}

\noindent{\bf Proof.}
For a point $z$ and a loop $\gamma$, let
$\theta_{\gamma}(z)$ indicate the winding number of $\gamma$ around $z$.
Moreover, for $\boldsymbol{k}  \in (\mathbb N \cup \{0\})^{\mathbb N}$   let 
\be \nonumber
\mathcal L_{\boldsymbol{k} }= \{\eta: |\{\gamma \in \eta: 
z \in \bar \gamma,\delta \leq {\rm diam}(\gamma) <R,
|\theta_{\gamma}(z)| = m \}|=k_m \text{ for all } m \in \mathbb N \}.
\ee
If a loop $\gamma$ has $|\theta_{\gamma}(z)| = m $, then the winding number
$\theta_{\gamma}(z)$ is $\pm m$ with equal probability under $P_{\lambda, \mu}$.
Finally, for $m \geq 1$ we have
\begin{eqnarray*}
\alpha_{z, \delta, R,m}& := &  \mu_{\delta, R}(\gamma: z \in \bar {\gamma}, |\theta_{\gamma}(z)| = m) \\
&=& \mu(\gamma \in \eta: z \in \bar \gamma, \delta \leq {\rm diam}(\gamma) <R, |\theta_{\gamma}(z)| = m) \\
&=& \mu(\gamma: z \in \bar \gamma,\gamma \not\subset B_{z, \delta}, \gamma \subset B_{z,R}, |\theta_{\gamma}(z)| = m) \\
&=& 
\frac{1}{\pi^2 m^2}\log\frac{R}{\delta} ,
\end{eqnarray*}
where, in the last equality, we have used Lemma~\ref{SecondLemma}.

For all $\boldsymbol{k}  \in (\mathbb N \cup \{0\})^{\mathbb N}$ we have 
\be \nonumber
P_{\lambda, \mu_{\delta, R}}( \mathcal L_{\boldsymbol{k} })
= \prod_{m=1}^{\infty} 
\frac{(\lambda \alpha_{z, \delta, R,m})^{k_m}}{{k_m}!} e^{-\lambda \alpha_{z, \delta, R,m}} ,
\ee
as for different $m$'s the sets of loops with those winding numbers are disjoint.
Hence with IR and UV cutoffs in place, denoting by $E_{\lambda,\delta,R}$ the expectation
with respect to the Poisson distribution $P_{\lambda,\mu_{\delta,R}}$ with intensity measure
$\mu_{\delta,R}$, we have, for all $z$,
\begin{eqnarray*}
\left\la V_{\beta}(z) \right\ra_{\delta, R} &=& E_{\lambda,\delta,R}\left(e^{i \beta N_{w}(z)}\right) \\
&=& \sum_{\boldsymbol{k}  \in (\mathbb N \cup \{0\})^{\mathbb N}} E_{\lambda,\delta,R}
\left(e^{i \beta N_{w}(z)} | \mathcal L_{\boldsymbol{k} }\right) 
P_{\lambda, \mu_{\delta, R}}(\mathcal L_{\boldsymbol{k} })\\
&=&  \sum_{\boldsymbol{k}  \in (\mathbb N \cup \{0\})^{\mathbb N}}
\prod_{m=1}^{\infty} (\cos(m \beta))^{k_m}
\frac{(\lambda \alpha_{z, \delta, R,m})^{k_m}}{{k_m}!} e^{-\lambda \alpha_{z, \delta, R,m}}
\\
&=& \prod_{m=1}^{\infty}e^{-\lambda \alpha_{z, \delta, R,m}(1-\cos(m \beta))}
\\
&=&\left(\frac{R}{\delta}\right)^{-\lambda  \sum_{m=1}^{\infty}\frac{1}{\pi^2 m^2}(1-\cos(m \beta))}
= \left(\frac{R}{\delta}\right)^{- \lambda \frac{\beta(2 \pi - \beta)}{4 \pi^2} } ,
\end{eqnarray*}
where the $\beta$ on the r.h.s. of the last equality is to be interpreted modulo $2\pi$. \fbox{}

\subsubsection{The $2$-point function in the layering model} \label{2-point-function-layering-model}

We now analyze the $2$-point function when the IR cutoff is removed by the charge conservation condition \eqref{cc}.

\begin{theorem} \label{thm-two-point-function}
If $\beta_{1}+\beta_{2}= 2 k \pi$ with $k \in {\mathbb Z}$, there is a positive constant $C_2 < \infty$ such that, for all $z_1 \neq z_2$,
\be \nonumber
\lim_{R \to \infty} \left\la V_{\beta_1}(z_1) V_{\beta_2}(z_2) \right\ra_{\delta, R}
= C_2 \left(\frac{|z_1-z_2|}{\delta} \right)^{-\frac{\lambda}{5}(2-\cos \beta_{1}- \cos \beta_{2})} .
\ee
As a consequence,
\be \nonumber
\lim_{\delta \to 0} \lim_{R \to \infty}
\frac{\left\la V_{\beta_1}(z_1) V_{\beta_2}(z_2) \right\ra_{\delta, R}}{\delta^{\frac{\lambda}{5}(2-\cos \beta_{1}- \cos \beta_{2})}}
= C_2 |z_1-z_2|^{-\frac{\lambda}{5}(2-\cos \beta_{1}- \cos \beta_{2})} .
\ee
\end{theorem}

\noindent{\bf Proof.}
Letting $d := |z_{1}-z_{1}|$, for given $\beta_1$ and $\beta_2$, and $d \geq \delta$, we have that
\begin{eqnarray*}
\left\la V_{\beta_1}(z_1) V_{\beta_2}(z_2) \right\ra_{\delta, R} 
& = & \left\la e^{i (\beta_1 N_{\ell}(z_1)+ \beta_2 N_{\ell}(z_2))} \right\ra_{\delta, R} \\
& = & \left\la e^{i (\beta_1 + \beta_2 ) N_{12}} \right\ra_{d, R}
\left\la e^{i \beta_1  N_{1}} \right\ra_{\delta, R} \left\la e^{i  \beta_2  N_{2}} \right\ra_{\delta, R} \, , 
\end{eqnarray*}
where $N_{12}$ is the number of loops that cover both
$z_1$ and $z_2$,  and $N_1$ ($N_2$) is the number of loops that cover 
$z_1$ but not  $z_2$ ($z_2$ but not  $z_1$, resp.). The two-point function
factorizes because the sets of loops contributing to $N_{1 2}$,
$N_1$ and $N_2$ are disjoint; the $\delta$ is replaced by $d$ in the first factor
in the second line because a loop covering both $z_1$ and $z_2$ must have
diameter at least $d$.

As in the $1$-point function calculation, we can write
\be \nonumber
\left\la e^{i (\beta_1 + \beta_2 ) N_{12}} \right\ra_{d, R} =
\sum_{n=0}^{\infty} \left(\cos(\beta_1 + \beta_2)\right)^n P_{\lambda, \mu^{loop}_{d, R}}(N_{1 2}=n)
= e^{-\lambda \alpha_{d, R}(z_1,z_2)\left( 1 - \cos(\beta_1 + \beta_2) \right)},
\ee
where $\alpha_{d, R}(z_1,z_2) \equiv \mu^{loop}_{d, R}(\gamma : z_1, z_2 \in \bar \gamma)$.
Similarly, if 
$\alpha_{\delta, R}(z_1,\neg z_2) = \mu^{loop}_{\delta, R}(\gamma : z_1 \in \bar \gamma,
z_2 \notin \bar \gamma)$ and $\alpha_{\delta, R}(\neg z_1,  z_2)$ is correspondingly
defined, then
\be \nonumber
\left\la e^{i \beta_1 N_{1}} \right\ra_{\delta, R} = e^{-\lambda\alpha_{\delta, R}(z_1,\neg z_2)(1-\cos\beta_1)}
\ee
and
\be \nonumber
\left\la e^{i \beta_2 N_{2}} \right\ra_{\delta, R} = e^{-\lambda\alpha_{\delta, R}(\neg z_1,z_2)(1-\cos\beta_2)} \, .
\ee
Combining the three terms we obtain
\begin{eqnarray*}
\lefteqn{
 \left\la e^{i (\beta_1 N_{\ell}(z_1)+ \beta_2 N_{\ell}(z_2))} \right\ra_{\delta, R} } \\
 & = & e^{-\lambda\alpha_{d, R}(z_1,z_2)(1-\cos(\beta_1 + \beta_2))} \,
 e^{-\lambda\alpha_{\delta, R}(z_1,\neg z_2)(1-\cos\beta_1)} \, e^{-\lambda\alpha_{\delta, R}(\neg z_1,z_2)(1-\cos\beta_2)} \, .
 \end{eqnarray*}
It is easy to see that $\lim_{R \rightarrow \infty}\alpha_{d, R}(z_1,z_2) = \infty$.
(This follows from the scale invariance of $\mu^{loop}$ by considering an increasing---in size---sequence of disjoint,
concentric annuli around $z_{1}$ and $z_{2}$ that are scaled versions of each other.) Hence, in order to remove the
IR cutoff, we must impose \eqref{cc}
and set $\beta_{1}+\beta_{2}= 2 k \pi$, so that  $1-\cos(\beta_1 + \beta_2)=0$.  
 
Assuming that $\beta_{1}+\beta_{2}= 2 k \pi$, we are left with
 \be \label{2pf-charge-cons}
  \left\la e^{i (\beta_1 N_{\ell}(z_1)+ \beta_2 N_{\ell}(z_2))} \right\ra_{\delta, R}
 =e^{
 -\lambda\alpha_{\delta, R}(z_1,\neg z_2)(1-\cos\beta_1)-\lambda\alpha_{\delta, R}(\neg z_1,z_2)(1-\cos\beta_2)} .
 \ee
To remove the infrared cutoff, we use the fact that the Brownian loop soup is \emph{thin}: If $z_1 \neq z_2$,
$\mu^{loop}(\gamma:z_1 \in \bar \gamma, z_2 \notin \bar \gamma, \diam(\gamma) \geq \delta)<\infty$
for any $\delta>0$ (see \cite{nawe}, Lemma 4).
By the obvious monotonicity of $\alpha_{ \delta, R} (z_1,\neg z_2)$ in $R$, this implies that
\begin{equation*}
\lim_{R \rightarrow \infty} \alpha_{\delta, R}(z_1,\neg z_2) = \mu^{loop}\{\gamma \in \eta: z_{1} \in \bar \gamma, 
 z_{2} \notin  \bar \gamma, {\rm diam}(\gamma) \geq \delta \} \equiv \alpha_{\delta}(z_1,\neg z_2) .
\end{equation*}
By scale, rotation and translation invariance of the Brownian loop measure $\mu^{loop}$,
$\alpha_{\delta}(z_1,\neg z_2) $ can only depend on the ratio $x=d/\delta$, so we can introduce the notation
$
\alpha(x) \equiv \alpha_{\delta}(z_1,\neg z_2).
$
The function $\alpha$ has the following properties, which are also immediate consequences of the scale,
rotation and translation invariance of the Brownian loop measure.
 \begin{itemize}
 \item $ \alpha(x) =\alpha_{\delta}(0,\neg z)$ for any $z$ such that $|z|=d$.
 \item For $\sigma \geq 1$, if $\delta < d$, letting $\alpha_{\delta, R}(z) \equiv \alpha_{z, \delta, R}=\mu_{\delta, R}(\gamma: z \in \bar \gamma)$,
 \begin{eqnarray} \label{alpha-sigma}
  \alpha( \sigma x) &=&\alpha_{\delta}(0,\neg \sigma z) \nonumber \\
  &=& \alpha_{\sigma \delta}(0,\neg \sigma z) + \alpha_{\delta, \sigma \delta}(0,\neg \sigma z) \nonumber \\
  &=&   \alpha(  x) +  \alpha_{\delta, \sigma \delta}(0)
   = \alpha(  x) +  \alpha_{1, \sigma }(0).
  \end{eqnarray}
 \end{itemize}

Now let 
 \be \nonumber
 G(x) := \left\la e^{i (\beta_1 N_{\ell}(z_1)+ \beta_2 N_{\ell}(z_2))} \right\ra_{\delta}
 := \lim_{R \rightarrow \infty}  \left\la e^{i (\beta_1 N_{\ell}(z_1)+ \beta_2 N_{\ell}(z_2))} \right\ra_{\delta, R} \, ;
 \ee
using \eqref{2pf-charge-cons} and the definition of the function $\alpha$, we can write 
 \be \nonumber
 G(x)=e^{- \lambda \alpha(x)(2-\cos \beta_{1}- \cos \beta_{2})} \, .
 \ee
Then, for $\sigma \geq 1$,
\be \nonumber
 G(\sigma x)=e^{- \lambda \alpha_{1,\sigma}(0)(2-\cos \beta_{1}- \cos \beta_{2})} G(x).
 \ee
Using Lemma~\ref{FirstLemma}, we have that
 \be \nonumber
 \alpha_{1,\sigma}(0) = \frac{1}{5} \log \sigma .
 \ee
 It then follows that, for $\sigma \geq 1$,
 \be \label{scaling-eq}
 G(\sigma x)= \sigma^{- \frac{\lambda}{5} (2-\cos \beta_{1}- \cos \beta_{2}) } G(x) .
 \ee

For $0<\sigma<1$, \eqref{alpha-sigma} implies
\be \nonumber
\alpha(\sigma x) = \alpha(x) - \alpha_{1,1/\sigma }(0).
\ee
But since
\be \nonumber
\alpha_{1,1/\sigma}(0) = -\frac{1}{5} \log \sigma ,
\ee
equation~\eqref{scaling-eq} is unchanged when $0<\sigma<1$.

The fact that \eqref{scaling-eq} is valid for all $\sigma>0$ immediately implies that
 \be \nonumber
G(x) = C_2 x^{-\frac{\lambda}{5}(2-\cos \beta_{1}- \cos \beta_{2})}
\ee
for some constant $C_2 >0$. \fbox{}

\subsection{Conformal covariance of the $n$-point functions}
We now analyze the $n$-point functions for general $n \geq 1$ and their conformal invariance properties.
In bounded domains $D \subset {\mathbb C}$, we show, for both models, how to remove the UV cutoff $\delta>0$
by dividing by $\delta^{ 2 \sum_{j=1}^n \Delta_{j}}$, with the appropriate $\Delta_j$'s. We also show that this procedure leads
to conformally covariant functions of the domain $D$. The scaling with $\delta$ originates from the fact that loops with diameter
less than $\delta$ can only wind around a single point in the limit $\delta \to 0$, and so for these small loops the $n$-point function
reduces to the product of 1-point functions.

In Section~\ref{n-point-plane}, we deal with the layering model in the full plane, $\mathbb C$, and show that, together
with the UV cutoff $\delta>0$, we can also remove the IR cutoff $R<\infty$, provided we impose the condition
$\sum_{j=1}^n \beta_j \in 2 \pi \Z$ (see \eqref{cc}).  We refer to this condition as  ``charge conservation''
because---apart from the periodicity---it is reminiscent of momentum or charge conservation for the vertex operators
of the free boson.

In the layering model the IR convergence (given ``charge conservation'') is due to the finiteness of the total mass of the loops which
cover some points but not others; this is basically the property that the soup of outer boundaries of a Brownian loop soup is \emph{thin}
in the language of Nacu and Werner \cite{nawe}.

\subsubsection{The layering model in finite domains} \label{n-point-finite-1}

In the theorem below, we let
$\left\langle \prod_{j=1}^n V_{\beta_j}(z_j) \right\rangle_{\delta,D} =
{\mathbb E}_{\lambda,\delta,D}\left(\prod_{j=1}^n e^{i \beta_j N_\ell(z_j)}\right)$ denote
the expectation of the product $\prod_{j=1}^n e^{i \beta_j N_\ell(z_j)}$ with respect to a loop soup in $D$
with intensity $\lambda>0$ containing only loops of diameter at least $\delta>0$, that is, with respect to the
distribution $P_{\lambda,\mu_{\delta,D}} \otimes P_{1/2}$, where $P_{\lambda,\mu_{\delta,D}}$
is the Poisson distribution with intensity measure
$\mu_{\delta,D}
=\mu_D \mathbbm{1}_{\{\diam(\gamma)\geq\delta\}}
=\mu \mathbbm{1}_{\{\gamma \subset D, \diam(\gamma)\geq\delta\}}$
and $P_{1/2}$ is the distribution of a countable sequence of independent Bernoulli random variables with
parameter $1/2$ (remember that each loop belongs to one of two classes with equal probability).
\begin{theorem} \label{thm-bounded-domains}
If $n \in {\mathbb N}$, $D \subset {\mathbb C}$ is bounded and $\boldsymbol{\beta}=(\beta_1,\ldots,\beta_n)$, then 
\begin{equation} \nonumber
\lim_{\delta \to 0} 
\frac{\left\langle \prod_{j=1}^n V_{\beta_j}(z_j)\right\rangle_{\delta, D}}{\delta^{\frac{\lambda}{5}\sum_{j=1}^n (1-\cos\beta_j)}}
=: \phi_D(z_1, \dots, z_n;\boldsymbol{\beta}) 
\end{equation}
exists and is finite and real. Moreover, if $D'$ is another bounded subset of $\mathbb C$ and $f:D \to D'$ is a conformal map such that
$z'_1=f(z_1),\ldots,z'_n=f(z_n)$, then
\begin{equation} \nonumber
\phi_{D'}(z'_1,\ldots,z'_n;\boldsymbol{\beta}) = \prod_{j=1}^n |f'(z_j)|^{-\frac{\lambda}{5}(1-\cos\beta_j)} \phi_D(z_1, \dots, z_n;\boldsymbol{\beta}) \, .
\end{equation}
\end{theorem}

The proof of the theorem will make use of the following lemma, where $B_{\delta}(z)$ denotes the disc of radius $\delta$ centered at $z$,
$\bar\gamma$ denotes the complement of the unique unbounded component of ${\mathbb C} \setminus \gamma$.
\begin{lemma} \label{lemma}
Let $D,D' \subset\mathbb C$ and let $f:D \to D'$ be a conformal map. For $n\geq1$, assume that $z_1,\ldots,z_n \in D$ are distinct
and that $z'_1=f(z_1),\ldots,z'_n=f(z_n)$, and let $s_j=|f'(z_j)|$ for $j=1, \ldots, n$ . Then we have that, for each $j=1,\ldots,n$,
\begin{eqnarray*} \label{eq-lemma}
\lefteqn{\mu_{D'}(\gamma: z'_j \in \bar \gamma, \bar\gamma \not\subset f(B_{\delta}(z_j)), z'_k \notin \bar\gamma \; \forall k \neq j) } \\
& - & \mu_{D'}(\gamma: z'_j \in \bar \gamma, \bar\gamma \not\subset B_{s_j \delta}(z'_j), z'_k \notin \bar\gamma \; \forall k \neq j) = o(1) \text{ as } \delta \to 0 .
\end{eqnarray*}
\end{lemma}

\noindent{\bf Proof.}
Let $B_{in}(z'_j)$ denote the largest (open) disc centered at $z'_j$ contained inside $f(B_{\delta}(z_j)) \cap B_{s_j \delta}(z'_j)$,
and $B_{out}(z'_j)$ denote the smallest disc centered at $z'_j$ containing $f(B_{\delta}(z_j)) \cup B_{s_j \delta}(z'_j)$. A moment
of thought reveals that, for $\delta$ sufficiently small,
\begin{eqnarray} \label{upper-bound}
\lefteqn{|\mu_{D'}(\gamma: z'_j \in \bar \gamma, \bar\gamma \not\subset f(B_{\delta}(z_j)), z'_k \notin \bar\gamma \; \forall k \neq j) } \nonumber \\
& - & \mu_{D'}(\gamma: z'_j \in \bar \gamma, \bar\gamma \not\subset B_{s_j \delta}(z'_j), z'_k \notin \bar\gamma \; \forall k \neq j)| \nonumber \\
& = & \mu_{D'}(\gamma: z'_j \in \bar\gamma, \bar\gamma \not\subset f(B_{\delta}(z_j)) \cap s_j B_{\delta}(z'_j), \bar\gamma \subset f(B_{\delta}(z_j)) \cup s_j B_{\delta}(z'_j), z'_k \notin \bar\gamma \; \forall k \neq j) \nonumber \\
& \leq & \mu_{D'}(\gamma: z'_j \in \bar\gamma, \bar\gamma \not\subset B_{in}(z'_j), \bar\gamma \subset B_{out}(z'_j)) \nonumber \\
& = & \frac{1}{5} \log\frac{\diam(B_{out}(z'_j))}{\diam(B_{in}(z'_j))}\, ,
\end{eqnarray}
where we have used equation \eqref{equation-log-over-5} from Theorem \ref{theorem-conf-restriction} in Appendix \ref{BrownianLoopMeasure}.
Note that, when $D' = {\mathbb C}$, the quantities above involving $\mu_{\mathbb C}$ are bounded because of the fact that
the Brownian loop soup is \emph{thin} (see~\cite{nawe}).

Since $f$ is analytic, for every $w \in \partial B_{\delta}(z_j)$, we have that
\begin{equation} \nonumber
|f(w)-z'_j| = 
s_j \delta + O(\delta^2) \, ,
\end{equation}
which implies that
\begin{equation} \nonumber
\lim_{\delta \to 0} \frac{\diam(B_{out}(z'_j))}{\diam(B_{in}(z'_j))}
= 1 \, .
\end{equation}
In view of \eqref{upper-bound}, this concludes the proof of the lemma. \fbox{} \\

\noindent{\bf Proof of Theorem \ref{thm-bounded-domains}.}
We first show that the limit is finite. We let $\eta$ denote a realization of loops and $\{X_{\gamma}\}_{\gamma \in \eta}$ a collection
of independent, Bernoulli, symmetric random variables taking values in $\{-1,1\}$. Moreover, let $[n] \equiv \{1,\dots,n\}$, let $\mathcal K$
denote the space of assignments of a nonnegative integer to each nonempty subset $S$ of $\{z_1,\ldots,z_n\}$, and for
$S \subset \{z_1,\ldots,z_n\}$, let $I_{S} \subset [n]$ be the set of indices such that $k \in I_S$ if and only if $z_k \in S$. We have that

\begin{eqnarray*}
\left\langle \prod_{j=1}^n V_{\beta_j}(z_j) \right\rangle_{\delta,D}
&=& {\mathbb E}_{\lambda,\delta,D}\left(e^{i \sum_{j=1}^n \beta_j N_{\ell}(z_j)}\right)  \\
&=& \sum_{\boldsymbol{k} \in \mathcal K} {\mathbb E}_{\lambda,\delta,D}\left(e^{i \sum_{j=1}^n \beta_j N_{\ell}(z_j)} | \mathcal L_{\boldsymbol{k} }\right) 
P_{\lambda, \mu_{\delta, D}}(\mathcal L_{\boldsymbol{k} })
\end{eqnarray*}
where  $\mathcal L_{\boldsymbol{k} }= \{\eta: \forall S \subset \{ z_1,\ldots,z_n \}, S \neq \emptyset,
|\{\gamma \in \eta: z_j \in \bar \gamma \quad \forall j \in I_S, z_j \not \in \bar \gamma \quad \forall j \not \in I_S \}|=\boldsymbol{k} (S) \}$.
With probability one with respect to $P_{\lambda, \mu_{\delta, D}}$ we have that, for each $j=1,\ldots,n$,
\be \nonumber
N_{\ell}(z_j)= \sum_{\gamma: z_j \in \bar \gamma, \diam(\gamma) \geq \delta} X_{\gamma}
= \sum_{S \subset \{ z_1,\ldots,z_n \}: z_j \in S} \sum_{\gamma: S \subset \bar\gamma, S^c \subset \bar\gamma^c, \diam(\gamma) \geq \delta} X_{\gamma}.
\ee
With the notation $\sum^{S} := \sum_{\gamma: S \subset \bar \gamma, S^c \subset \bar\gamma^c, \diam(\gamma) \geq \delta}$, and letting
$X$ denote a $(\pm)$-valued, symmetric random variable, we have that
 \begin{eqnarray*}
 {\mathbb E}_{\lambda,\delta, R}\left(e^{i \sum_{j=1}^n \beta_j N_{\ell}(z_j)} | \mathcal L_{\boldsymbol{k} } \right) 
 &=&
 {\mathbb E}_{\lambda,\delta,R}\left(e^{i \sum_{j=1}^n \beta_j \sum_{S \subset \{ z_{1},\ldots,z_{n} \}: z_j \in S} 
 \sum^{S} X_{\gamma} } | \mathcal L_{\boldsymbol{k} } \right) \\
 &=&
 {\mathbb E}_{\lambda,\delta, R}\left(e^{i \sum_{S \subset \{ z_{1},\ldots,z_{n} \}} \sum^{S} (\sum_{j \in I_{S}} \beta_j ) X_{\gamma} }
 | \mathcal L_{\boldsymbol{k} }\right)\\
   &=&
 \prod_{S \subset \{z_1,\ldots,z_n\}, S \neq \emptyset} \left(E\left(e^{i (\sum_{j \in I_{S}} \beta_j ) X } \right)\right)^{\boldsymbol{k}(S)} \\
   &=&
 \prod_{S \subset \{z_1,\ldots,z_n\}, S \neq \emptyset} \Big(\cos\Big(\sum_{j \in I_{S}} \beta_j  \Big)\Big)^{\boldsymbol{k} (S)} .
  \end{eqnarray*}

Next, given $S \subset \{ z_1, \ldots, z_n \}$ with $|S| \geq 2$, let
$\alpha_D(S) := \mu_{D}(\gamma: S \subset \bar \gamma, S^c \subset \bar\gamma^c)$ and
$\alpha_{\delta,D}(z_j) := \mu_{D}(\gamma: \diam(\gamma) \geq \delta, z_j \in \bar \gamma, z_k \notin \bar\gamma \; \forall k \neq j)$.
Furthermore, let $m = \min_{i,j: i \neq j} |z_i -  z_j| \wedge \min_i \dist(z_i,\partial D)$ and note that, when $\delta<m$, we can write

\begin{eqnarray}
\left\langle \prod_{j=1}^n V_{\beta_j}(z_j) \right\rangle_{\delta,D}
&=&
\sum_{\boldsymbol{k} \in \mathcal K}  \prod_{S \subset \{z_1,\ldots,z_n \}, |S|>1} \Big(\cos\Big(\sum_{k \in I_S} \beta_k \Big)\Big)^{\boldsymbol{k} (S)}
\frac{ (\lambda \alpha_D(S))^{\boldsymbol{k}(S)}}{({\boldsymbol{k} (S)})!} e^{- \lambda \alpha_D(S)} \nonumber \\
& & \prod_{j=1}^n \Big( \cos\beta_j \Big)^{\boldsymbol{k}(z_j)} \frac{ (\lambda \alpha_{\delta,D}(z_j))^{\boldsymbol{k}(z_j)}}{({\boldsymbol{k}(z_j)})!} e^{- \lambda \alpha_{\delta,D}(z_j)} \nonumber \\
&=& 
\prod_{S \subset \{z_1,\ldots,z_n \}, |S|>1} \exp{\Big[-\lambda \alpha_D(S) \Big(1-\cos\Big(\sum_{k \in I_S} \beta_k \Big)\Big)\Big]} \nonumber \\
& & \prod_{j=1}^n \exp{\Big[-\lambda \alpha_{\delta,D}(z_j)(1 - \cos \beta_j) \Big]} \, . \label{n-point-func}
\end{eqnarray}

For every $j=1,\ldots,n$, using Lemma~\ref{FirstLemma}, we have that
 \begin{eqnarray*}
  \alpha_{\delta,D}(z_j) &=& \mu_{D}(\gamma: \diam(\gamma) \geq \delta, z_j \in \bar \gamma, z_k \notin \bar\gamma \; \forall k \neq j)  \\
  & = & \mu_{D}(\gamma: m > \diam(\gamma) \geq \delta, z_j \in \bar \gamma) \\
  & + & \mu_{D}(\gamma: \diam(\gamma) \geq m, z_j \in \bar \gamma, z_k \notin \bar\gamma \; \forall k \neq j)\\
  &=&
  \frac{1}{5} \log\frac{m}{\delta} + \alpha_{m,D}(z_j).
 \end{eqnarray*}
Therefore, we obtain
 \begin{eqnarray*}
\lim_{\delta \to 0} 
\frac{\left\langle \prod_{j=1}^n V_{\beta_j}(z_j) \right\rangle_{\delta,D}}{\delta^{\frac{\lambda}{5}\sum_{j=1}^n (1-\cos(\beta_j))}}
& = & \prod_{S \subset \{z_1,\ldots,z_n \}, |S|>1} \exp{\Big[-\lambda \alpha_D(S) \Big(1-\cos\Big(\sum_{k \in I_S} \beta_k \Big)\Big)\Big]} \\
 & & m^{ -\frac{\lambda}{5}\sum_{j=1}^n (1-\cos\beta_j)} e^{-\sum_{j=1}^n \lambda \alpha_{m,D}(z_j)(1-\cos\beta_j)} \\
& = & m^{ -\frac{\lambda}{5}\sum_{j=1}^n (1-\cos\beta_j)}
\exp{\Big[ -\lambda \sum_{j=1}^n \alpha_{m,D}(z_j)(1-\cos\beta_j) \Big]} \\
& & \exp{\Big[ -\lambda \sum_{S \subset \{z_1,\ldots,z_n \}, |S|>1} \alpha_D(S) \Big(1-\cos\Big(\sum_{k \in I_S} \beta_k \Big)\Big) \Big]} \\
& =: & \phi_D(z_1, \dots, z_n; \boldsymbol{\beta}) \, .
\end{eqnarray*}
This concludes the first part of the proof.

To prove the second part of the theorem, using \eqref{n-point-func}, we write
\begin{eqnarray} \label{n-point-funct}
\left\langle \prod_{j=1}^n V_{\beta_j}(z_j) \right\rangle_{\delta,D} & = & \exp{\Big[ -\lambda \sum_{S \subset \{ z_1, \ldots, z_n \}, |S| \geq 2} \alpha_D(S) \Big(1-\cos\Big(\sum_{k \in I_S} \beta_k \Big)\Big) \Big]} \nonumber \\
                                                              &     & \prod_{j=1}^n \exp{\left[-\lambda \alpha_{\delta,D}(z_j)(1 - \cos \beta_j) \right]} \, .
\end{eqnarray}
For each $S \subset \{ z_1, \ldots, z_n \}$ with $|S| \geq 2$, $\alpha_D(S)$ is invariant under conformal transformations, that is,
if $f:D \to D'$ is a conformal map from $D$ to another bounded domain $D'$, and $S' = \{ z'_1,\ldots,z'_n \}$, where $z'_1=f(z_1),\ldots,z'_n=f(z_n)$,
then $\alpha_{D'}(S') = \alpha_D(S)$. Therefore, the first exponential term in \eqref{n-point-funct} is also invariant under conformal transformations.
This implies that, for $\delta$ sufficiently small,
\begin{equation} \nonumber
\frac{\left\langle \prod_{j=1}^n V_{\beta_j}(z_j) \right\rangle_{\delta,D}}{\left\langle \prod_{j=1}^n V_{\beta_j}(z_j') \right\rangle_{\delta,D'}}
= \prod_{j=1}^n \exp{\left\{-\lambda \left[\alpha_{\delta,D}(z_j)-\alpha_{\delta,D'}(z'_j)\right](1 - \cos \beta_j) \right\}} \, . 
\end{equation}

Writing
\begin{eqnarray*}
\alpha_{\delta,D}(z_j) & = & \mu_D(\gamma: \diam(\gamma) \geq \delta, z_j \in \bar \gamma, \bar\gamma \subset B_{\delta}(z_j)) \\
& + & \mu_D(\gamma: z_j \in \bar \gamma, \bar\gamma \not\subset B_{\delta}(z_j), z_k \notin \bar\gamma \; \forall k \neq j)
\end{eqnarray*}
and noticing that $\mu_D(\gamma: \diam(\gamma) \geq \delta, z_j \in \bar \gamma, \bar\gamma \subset B_{\delta}(z_j))
= \mu_{D'}(\gamma: \diam(\gamma) \geq \delta, z'_j \in \bar \gamma, \bar\gamma \subset B_{\delta}(z'_j))$ (where we have assumed,
without loss of generality, that $\delta$ is so small that $B_{\delta}(z_j) \subset D$ and $B_{\delta}(z'_j) \subset D'$), we have that
\begin{eqnarray*}
\alpha_{\delta,D}(z_j) - \alpha_{\delta,D'}(z'_j) & = &
\mu_D(\gamma: z_j \in \bar \gamma, \bar\gamma \not\subset B_{\delta}(z_j), z_k \notin \bar\gamma \; \forall k \neq j) \\
& - & \mu_{D'}(\gamma: z'_j \in \bar \gamma, \bar\gamma \not\subset B_{\delta}(z'_j), z'_k \notin \bar\gamma \; \forall k \neq j) \, .
\end{eqnarray*}
To evaluate this difference, using conformal invariance, we write
\begin{eqnarray*}
\lefteqn{ \mu_D(\gamma: z_j \in \bar \gamma, \bar\gamma \not\subset B_{\delta}(z_j), z_k \notin \bar\gamma \; \forall k \neq j) } \\
& = & \mu_{D'}(\gamma: z'_j \in \bar \gamma, \bar\gamma \not\subset f(B_{\delta}(z_j)), z'_k \notin \bar\gamma \; \forall k \neq j) \, .
\end{eqnarray*}
Thus, letting $s_j = |f'(z_j)|$ and using Lemmas \ref{lemma-measure-equality} and \ref{FirstLemma} from Appendix \ref{BrownianLoopMeasure},
we can write
\begin{eqnarray*}
\alpha_{\delta,D}(z_j) - \alpha_{\delta,D'}(z'_j) & = &
\mu_{D'}(\gamma: z'_j \in \bar \gamma, \bar\gamma \not\subset f(B_{\delta}(z_j)), z'_k \notin \bar\gamma \; \forall k \neq j) \\
& & - \mu_{D'}(\gamma: z'_j \in \bar \gamma, \bar\gamma \not\subset B_{\delta}(z'_j), z'_k \notin \bar\gamma \; \forall k \neq j) \\
& = & \mu_{D'}(\gamma: z'_j \in \bar \gamma, \bar\gamma \not\subset f(B_{\delta}(z_j)), z'_k \notin \bar\gamma \; \forall k \neq j) \\
& & - \mu_{D'}(\gamma: z'_j \in \bar \gamma, \bar\gamma \not\subset B_{s_j \delta}(z'_j), z'_k \notin \bar\gamma \; \forall k \neq j) \\
& & - [ \mu_{D'}(\gamma: z'_j \in \bar \gamma, \bar\gamma \not\subset B_{\delta}(z'_j), z'_k \notin \bar\gamma \; \forall k \neq j) \\
& &  -\mu_{D'}(\gamma: z'_j \in \bar \gamma, \bar\gamma \not\subset B_{s_j \delta}(z'_j), z'_k \notin \bar\gamma \; \forall k \neq j) ] \\
& = & \mu_{D'}(\gamma: z'_j \in \bar \gamma, \bar\gamma \not\subset f(B_{\delta}(z_j)), z'_k \notin \bar\gamma \; \forall k \neq j) \\
& & - \mu_{D'}(\gamma: z'_j \in \bar \gamma, \bar\gamma \not\subset B_{s_j \delta}(z'_j), z'_k \notin \bar\gamma \; \forall k \neq j) \\
& & - \frac{1}{5} \log s_j \, .
\end{eqnarray*}

Using Lemma \ref{lemma}, we obtain that
$\alpha_{\delta,D}(z_j) - \alpha_{\delta,D'}(z'_j) = -\frac{1}{5} \log |f'(z_j)| + o(1)$ as $\delta \to 0$, which gives
\begin{equation} \nonumber
\frac{\left\langle V_{\bf\beta}(z_1,\ldots,z_n) \right\rangle_{\delta,D}}{\left\langle V_{\bf\beta}(z'_1,\ldots,z'_n) \right\rangle_{\delta,D'}}
= e^{-o(1)} \prod_{j=1}^n |f'(z_j)|^{\frac{\lambda}{5}(1-\cos\beta_j)} \text{ as } \delta \to 0.
\end{equation}
Letting $\delta \to 0$ concludes the proof. \fbox{} \\

\subsubsection{The winding model in finite domains} \label{n-point-finite-2}

As above, let $N_{w}(z)$ denote the total number of windings of all loops of a given soup around $z \in {\mathbb C}$.
We have the following theorem.
\begin{theorem} \label{thm-bounded-domains-winding}
If $n \in {\mathbb N}$, $D \subset {\mathbb C}$ is bounded and $\boldsymbol{\beta}=(\beta_1,\ldots,\beta_n)$, then 
\begin{equation} \nonumber
\lim_{\delta \to 0} 
\frac{\left\langle e^{i \beta_1N_{w}(z_1)} \dots e^{i \beta_n N_{w}(z_n)} \right\rangle_{\delta, D}}{\delta^{ \frac{\lambda}{4\pi^2} \sum_{j=1}^n \beta_j(2\pi-\beta_j) }}
=: \psi_D(z_1,\ldots,z_n;\boldsymbol{\beta}) 
\end{equation}
exists and is finite and real. Moreover, if $D'$ is another bounded subset of $\mathbb C$ and $f:D \to D'$ is a conformal map such that
$z'_1=f(z_1),\ldots,z'_n=f(z_n)$, then
\begin{equation} \nonumber
\psi_{D'}(z'_1,\ldots,z'_n;\boldsymbol{\beta}) = \prod_{j=1}^n \left|f'(z_j)\right|^{-\lambda\frac{\beta_j(2\pi-\beta_j)}{4\pi^2}} \psi_D(z_1,\ldots,z_n; \boldsymbol{\beta}) \, ,
\end{equation}
where, in the exponent, the $\beta_j$'s are to be interpreted modulo $2\pi$.
\end{theorem}

\noindent{\bf Proof.}
The proof is analogous to that of Theorem \ref{thm-bounded-domains}. Let
$\alpha_D(S;k_{i_1},\ldots,k_{i_l}) := \mu_D(\gamma: N_{\gamma}(z_{i_j}) = k_{i_j} \text{ for each } z_{i_j} \in S \text{ and } S^c \subset \bar\gamma^c)$,
and $\alpha_{\delta,D}(z_j;k) := \mu_{D}(\gamma: z_j \in \bar \gamma, \theta_{\gamma}(z_j) = k, z_k \notin \bar\gamma \; \forall k \neq j)$.
With this notation we can write
\begin{eqnarray*} \label{n-point-funct-winding}
\lefteqn{\left\langle e^{\beta_1N_{w}(z_1)} \ldots e^{\beta_nN_{w}(z_n)} \right\rangle_{\delta,D} } \\
& = & \exp{\Big[ -\lambda \sum_{l=2}^n \sum_{\stackrel{S \subset \{ z_1, \ldots, z_n \}}{|S|=l}}
\sum_{k_{i_1},\ldots,k_{i_l}=-\infty}^{\infty} \alpha_D(S;k_{i_1},\ldots,k_{i_l}) \left(1-\cos (k_{i_1}\beta_{i_1}+\ldots+k_{i_l}\beta_{i_l}) \right) \Big]} \nonumber \\
                                                              &     & \prod_{j=1}^n \exp{\Big[-\lambda \sum_{k=-\infty}^{\infty} \alpha_{\delta,D}(z_j;k)(1 - \cos(k\beta_j)) \Big]} \, ,
\end{eqnarray*}
For each $S \subset \{ z_1, \ldots, z_n \}$ with $|S| = l \geq 2$, $\alpha_D(S;k_{i_1},\ldots,k_{i_l})$ is invariant under conformal transformations; therefore,
\begin{equation} \nonumber
\frac{\left\langle e^{i \beta_1 N_{w}(z_1)} \ldots e^{i \beta_n N_{w}(z_n)} \right\rangle_{\delta,D}}{\left\langle e^{i \beta_1 N_{w}(z'_1)} \ldots e^{i \beta_n N_{w}(z'_n)} \right\rangle_{\delta,D'}}
= \prod_{j=1}^n \exp{\Big\{-\lambda \sum_{k=-\infty}^{\infty} \left[\alpha_{\delta,D}(z_j;k)-\alpha_{\delta,D'}(z'_j;k)\right](1 - \cos(k\beta_j)) \Big\}} \, . 
\end{equation}
Proceeding as in the proof of Theorem \ref{thm-bounded-domains}, but using Lemma~\ref{SecondLemma}
instead of Lemma~\ref{FirstLemma}, gives
\begin{equation} \nonumber
\alpha_{\delta,D}(z_j;k)-\alpha_{\delta,D'}(z'_j;k) = -c_k \log|f'(z_j)| + o(1) \text{ as } \delta \to 0 \, ,
\end{equation}
where $c_k = \frac{1}{2\pi^2 k^2}$ for $k \in {\mathbb Z} \setminus \{0\}$ and $c_0 = 1/30$.

This, together with the observation, already used at the end of Section~\ref{Sec:WindingOperator}, that
$\sum_{k=-\infty}^{\infty} c_k (1-\cos(k\beta)) = \frac{\beta(2\pi-\beta)}{4\pi^2}$ (where, in the right hand side,
$\beta$ should be interpreted modulo $2\pi$), readily implies the statement of the theorem. \fbox{}

\subsubsection{The layering model in the plane} \label{n-point-plane}

In the theorem below, we let $\left\langle \prod_{j=1}^n V_{\beta_j}(z_j) \right\rangle_{\delta,R}$ denote
the expectation of the product $\prod_{j=1}^n e^{i \beta_j N_\ell(z_j)}$ with respect to a loop soup in $\mathbb C$
with intensity $\lambda>0$ containing only loops $\gamma$ of diameter $0<\delta\leq\diam(\gamma)<R<\infty$.
\begin{theorem} \label{thm-n-point-plane}
If $n \in {\mathbb N}$ and $\boldsymbol{\beta}=(\beta_1,\ldots,\beta_n)$
with $|\boldsymbol{\beta}| = \sum_{j=1}^n \beta_j \in 2\pi {\mathbb Z}$, then 
\begin{equation} \nonumber
\lim_{\delta \to 0, R \to \infty} 
\frac{\left\langle \prod_{j=1}^n V_{\beta_j}(z_j) \right\rangle_{\delta, R}}{\delta^{ \frac{\lambda}{5}\sum_{j=1}^n (1-\cos\beta_j)}}
=: \phi_{\mathbb C}(z_1, \dots, z_n;\boldsymbol{\beta}) 
\end{equation}
exists and is finite and real. Moreover, if $f:{\mathbb C} \to {\mathbb C}$ is a conformal map such that
$z'_1=f(z_1),\ldots,z'_n=f(z_n)$, then
\begin{equation} \nonumber
\phi_{\mathbb C}(z'_1,\ldots,z'_n;\boldsymbol{\beta}) = \prod_{j=1}^n \left|f'(z_j)\right|^{-\frac{\lambda}{5}(1-\cos\beta_j)} \phi_{\mathbb C}(z_1, \dots, z_n;\boldsymbol{\beta}) \, .
\end{equation}
\end{theorem}

\noindent{\bf Proof sketch.} The beginning of the proof proceeds like that of Theorem~\ref{thm-bounded-domains} until equation~\eqref{n-point-func},
leading to the following equation:
\begin{eqnarray*}
\left\langle \prod_{j=1}^n V_{\beta_j}(z_j) \right\rangle_{\delta,R} & = &
\prod_{S \subset \{z_1,\ldots,z_n \}, 1<|S|<n} \exp{\Big[-\lambda \alpha_R(S) \Big(1-\cos\Big(\sum_{k \in I_S} \beta_k \Big)\Big)\Big]} \\
& & \prod_{j=1}^n \exp{\left[-\lambda \alpha_{\delta,R}(z_j)(1 - \cos \beta_j) \right]} \, ,
\end{eqnarray*}
where $\alpha_R(S) := \mu_{\mathbb C}(\gamma: S \subset \bar \gamma, S^c \subset \bar\gamma^c, \diam(\gamma)<R)$,
for $S \subset \{ z_1, \ldots, z_n \}$ with $2\leq|S|<n$, and
$\alpha_{\delta,R}(z_j) := \mu_{\mathbb C}(\gamma: \delta \leq \diam(\gamma) < R, z_j \in \bar \gamma, z_k \notin \bar\gamma \; \forall k \neq j)$,
and where $I_S$ denotes the set of indices such that $k \in I_S$ if and only if $z_k \in S$.

Note, in the equation above, the condition $|S|<n$ in the first product on the right hand side; this condition comes from the fact that
the term $-\lambda \alpha_R(S)$ with $S = \{z_1,\ldots,z_n\}$ is multiplied by $1-\cos(\sum_{k=1}^n \beta_k)=0$, where we have
used the ``charge conservation'' condition $|\boldsymbol{\beta}| = \sum_{j=1}^n \beta_j \in 2\pi {\mathbb Z}$.

For every $j=1,\ldots,n$, using Lemma~\ref{FirstLemma}, we have that
 \begin{eqnarray*}
  \alpha_{\delta,R}(z_j) &=& \mu_{\mathbb C}(\gamma: \delta \leq \diam(\gamma) < R, z_j \in \bar \gamma, z_k \notin \bar\gamma \; \forall k \neq j)  \\
  & = & \mu_{\mathbb C}(\gamma: m > \diam(\gamma) \geq \delta, z_j \in \bar \gamma) \\
  & + & \mu_{D}(\gamma: m\leq\diam(\gamma)<R, z_j \in \bar \gamma, z_k \notin \bar\gamma \; \forall k \neq j)\\
  &=&
  \frac{1}{5} \log\frac{m}{\delta} + \alpha_{m,R}(z_j).
 \end{eqnarray*}
Now note that monotonicity and the fact that the Brownian loop soup is \emph{thin} \cite{nawe} imply that
$\alpha_{m,{\mathbb C}}(z_j) := \lim_{R \to \infty} \alpha_{m,R}(z_j)$ and $\alpha_{\mathbb C}(S) := \lim_{R \to \infty} \alpha_R(S)$,
for $S \subset \{ z_1, \ldots, z_n \}$ with $2\leq|S|<n$, exist and are bounded.
After letting $R\to\infty$, the proof proceeds like that of Theorem~\ref{thm-bounded-domains}, with $D=D'=\mathbb C$. \noindent \fbox{} \\

We have already seen the behavior of the 2-point function in the layering model in Section~\ref{2-point-function-layering-model}; the corollary below
deals with the 3-point function.
\begin{corollary} \label{3-point-function}
Let $z_1, z_2, z_3 \in {\mathbb C}$ be three distinct points, then we have that
\begin{eqnarray*}
\lefteqn{
\phi_{\mathbb C}(z_1,z_2,z_3;\beta_1,\beta_2,\beta_3) = } \\
& & C_{3} \left| \left( {1 \over |z_1 - z_2|} \right)^{\Delta_{l}(\beta_1) + \Delta_{l}(\beta_2) - \Delta_{l}(\beta_3)}  \left( {1 \over |z_1 - z_3|} \right)^{\Delta_{l}(\beta_1) + \Delta_{l}(\beta_3) - \Delta_{l}(\beta_2)} 
 \left( {1 \over |z_2 - z_3|} \right)^{\Delta_{l}(\beta_2) + \Delta_{l}(\beta_3) - \Delta_{l}(\beta_1)} \right|^2
\end{eqnarray*}
for some constant $C_{3}$.
\end{corollary}

\noindent{\bf Proof.}
Theorem~\ref{thm-n-point-plane} implies that the 3-point function in the full plane transforms covariantly under
conformal maps. This immediately implies the corollary following standard argument (see, e.g., \cite{dfms}).
We briefly sketch those arguments below for the reader's convenience.

Scale invariance, rotation invariance, and translation invariance immediately imply that there are constants $C_{abc}$ such that
\be \label{eq-invariance}
\phi_{\mathbb C}(z_1,z_2,z_3;\beta_1,\beta_2,\beta_3) = \sum C_{abc} z_{12}^{-a} z_{13}^{-b} z_{23}^{-c} \, ,
\ee
where $z_{ij}=|z_i-z_j|$ and the sum is over all triplets $a,b,c \geq 0$ satisfying $a+b+c = 2(\Delta_l(\beta_1) + \Delta_l(\beta_2) + \Delta_l(\beta_3))$
(the constraint on the exponents $a,b,c$ follows from Theorem~\ref{thm-n-point-plane} applied to scale transformations).

Now let $f$ be a conformal transformation from $\mathbb C$ to $\mathbb C$; $f$ is then a M\"obius transformation and has the form
$f(z) = \frac{Az+B}{Cz+D}$, with $f'(z) = \frac{AD-BC}{(Cz+D)^2}$. Letting $\gamma_j := |f'(z_j)|^{-1}$, if $\tilde z = f(z)$, it is easy
to check that $\tilde z_{ij} = \gamma_i^{-1/2} \gamma_j^{-1/2} z_{ij}$. Using this fact and Theorem~\ref{thm-n-point-plane}, we have that
\begin{eqnarray*}
\phi_{\mathbb C}(\tilde z_1, \tilde z_2, \tilde z_3;\beta_1,\beta_2,\beta_3) & = &
\left(\gamma_1^{\Delta_l(\beta_1)} \gamma_2^{\Delta_l(\beta_2)} \gamma_3^{\Delta_l(\beta_3)}\right)^2 \sum C_{abc} z_{12}^{-a} z_{13}^{-b} z_{23}^{-c} \\
& = & \left(\gamma_1^{\Delta_l(\beta_1)} \gamma_2^{\Delta_l(\beta_2)} \gamma_3^{\Delta_l(\beta_3)}\right)^2
\sum C_{abc} \frac{\tilde z_{12}^{-a} \tilde z_{13}^{-b} \tilde z_{23}^{-c}}{\gamma_1^{a/2+b/2} \gamma_2^{a/2+c/2} \gamma_3^{b/2+c/2}} \, .
\end{eqnarray*}
For this last expression to be of the correct form~\eqref{eq-invariance}, the $\gamma$'s need to cancel; this immediately leads to the relations
$a = 2(\Delta_l(\beta_2) + \Delta_l(\beta_3) - \Delta_l(\beta_1))$, $b = 2(\Delta_l(\beta_1) + \Delta_l(\beta_3) - \Delta_l(\beta_2))$,
$c = 2(\Delta_l(\beta_1) + \Delta_l(\beta_2) - \Delta_l(\beta_3))$. \fbox{}

\appendix

\section{Appendix: Occupation field and Gaussian free field}\label{appendix}

The theorem used in the proof of Lemma~\ref{coupling3}, stated below, is a version of a recent result
of Le Jan~\cite{lejan1} (see also \cite{lejan2} and Theorem 4.5 of \cite{sznitman-notes}). In this appendix
we give a self-contained proof, both for completeness and because the occupation field $\{L_x\}$ in
Theorem~\ref{gff-of} below is not exactly the same as the occupation field of~\cite{lejan1,lejan2} and
\cite{sznitman-notes}. Indeed, the loop soup that appears in~\cite{lejan1,lejan2} and in Theorem~4.5
of \cite{sznitman-notes} is not the same as the loop soup studied in these lecture notes, although the
two are closely related. In this appendix we use the notation of Sections~\ref{sec:rwls} and \ref{sec:gff}.
\begin{theorem} \label{gff-of}
Let $m$ be a nonnegative function and ${\bf k}$ denote the collection $\{4 (e^{m^2(x)}-1)\}_{x \in D^{\#}}$.
Then the occupation field $\{L_x(\tilde{\cal A}_{1/2,m})\}_{x \in D^{\#}}$ has the same distribution as the field
$\{\frac{1}{2} \phi^2_x\}_{x \in D^{\#}}$, where $\{\phi_x\}_{x \in D^{\#}}$ is the Gaussian free field with
covariance ${\mathbb E}(\phi_x \phi_y) = G^{\bf k}_D(x,y)$.
\end{theorem}

\noindent{\bf Proof.} We will show that the Laplace transform of the occupation field
is the same as that of the field $\{\frac{1}{2} \phi^2_x\}_{x \in D^{\#}}$. For that purpose, it is crucial to notice that
$L_x({\cal A}_{\lambda,m})$ has the gamma distribution with parameters
$\sum_{\tilde\gamma} N_{\tilde\gamma} n(x,\tilde\gamma) +1/2$ and $\frac{1}{k_x+4}$, and consequently,
for any real number $v$ and any collection $\{n_{\tilde\gamma}\}$ of nonnegative numbers,
\begin{equation} \nonumber
{\mathbb E}[\exp(-v L_x) | \{N_{\tilde\gamma}\}=\{n_{\tilde\gamma}\}]
= \left( 1+\frac{v}{k_x+4} \right)^{-1/2}
\prod_{\tilde\gamma} \left( 1+\frac{v}{k_x+4} \right)^{-n_{\tilde\gamma}n(x,\tilde\gamma)} \, ,
\end{equation} 
where the product if over all unrooted lattice loops in $D$.

Let $\{v_x\}_{x \in D^{\#}}$ be a collection of real numbers, the Laplace transform of the
occupation field is given by the following expectation, where the sum $\sum_{\{n_{\tilde\gamma}\}}$
is over all collections of possible multiplicities for the lattice loops $\tilde\gamma$ in $D$,
\begin{eqnarray*}
\lefteqn{ {\mathbb E}\left[\exp\left(-\sum_{x \in D^{\#}} v_x L_x \right)\right] } \\
& = & \sum_{\{n_{\tilde\gamma}\}} {\mathbb E}\left[\exp\left(-\sum_{x \in D^{\#}} v_x L_x \right) \mid \{N_{\tilde\gamma}\}=\{n_{\tilde\gamma}\} \right]
{\mathbb P}(\{N_{\tilde\gamma}\}=\{n_{\tilde\gamma}\}) \, .
\end{eqnarray*}

Recalling~\eqref{rwls}, we can write
\begin{eqnarray*}
\lefteqn{ {\mathbb E}\left[\exp\left(-\sum_{x \in D^{\#}} v_x L_x \right)\right] } \\
& = & \sum_{\{n_{\tilde\gamma}\}} \left[\prod_{x \in D^{\#}} \left( 1+\frac{v_x}{k_x+4} \right)^{-1/2}
\prod_{\tilde\gamma} \left( 1+\frac{v_x}{k_x+4} \right)^{-n_{\tilde\gamma}n(x,\tilde\gamma)}\right] \\
& & \prod_{\tilde\gamma} \exp{\left(-\lambda\nu^{u,{\bf k}}_D(\tilde\gamma)\right)}
\frac{1}{n_{\tilde\gamma}!} \left(\lambda \nu^{u,{\bf k}}_D(\tilde\gamma)\right)^{n_{\tilde\gamma}} \\
& = & \left(\prod_{x \in D^{\#}} (k_x+4)^{1/2}\right) \left(\prod_{\tilde\gamma} \exp{\left(-\lambda\nu^{u,{\bf k}}_D(\tilde\gamma)\right)}\right)
\left(\prod_{x \in D^{\#}} (v_x+k_x+4)^{-1/2}\right) \\
& & \sum_{\{n_{\tilde\gamma}\}} \prod_{\tilde\gamma} \frac{1}{n_{\tilde\gamma}!}
\left[\lambda \nu^{u,{\bf k}}_D(\tilde\gamma) \, \prod_{x \in D^{\#}} \left(\frac{v_x+k_x+4}{k_x+4} \right)^{-n(x,\tilde\gamma)} \right]^{n_{\tilde\gamma}} \, .
\end{eqnarray*}

Now let $\rho_{\tilde\gamma}$ denote the number of rooted loops in the equivalence class $\tilde\gamma$, then we can write
\begin{equation} \nonumber
\nu^{u,{\bf k}}_D(\tilde\gamma) 
= \frac{\rho_{\tilde\gamma}}{|\tilde\gamma|} \prod_{x \in D^{\#}} (k_x+4)^{-n(x,\tilde\gamma)}
\end{equation}
and, letting ${\bf k}+{\bf v}$ denote the collection $\{k_x+v_v\}_{x \in D^{\#}}$, it is easy to see that
\begin{equation} \nonumber
\nu^{u,{\bf k}}_D(\tilde\gamma) \, \prod_{x \in D^{\#}} \left(\frac{v_x+k_x+4}{k_x+4} \right)^{-n(x,\tilde\gamma)}
= \nu^{u,{\bf k}+{\bf v}}_D(\tilde\gamma) \, .
\end{equation}
Using this fact, the fact that~\eqref{rwls} defines a probability distribution, and Eq.~\eqref{partition-function},
one has that
\begin{eqnarray*}
\lefteqn{ \sum_{\{n_{\tilde\gamma}\}} \prod_{\tilde\gamma} \frac{1}{n_{\tilde\gamma}!}
\left[\lambda \nu^{u,{\bf k}}_D(\tilde\gamma) \, \prod_{x \in D^{\#}} \left(\frac{v_x+k_x+4}{k_x+4} \right)^{-n(x,\tilde\gamma)} \right]^{n_{\tilde\gamma}} } \\
& = & \sum_{\{n_{\tilde\gamma}\}} \prod_{\tilde\gamma} \frac{1}{n_{\tilde\gamma}!} \left(\lambda \nu^{u,{\bf k}+{\bf v}}_D(\tilde\gamma) \right)^{n_{\tilde\gamma}} \\
& = & {\cal Z}_{\lambda,{\bf k}+{\bf v}} = \exp{\left(\lambda \sum_{\tilde\gamma} \nu^{u,{\bf k}+{\bf v}}_D(\tilde\gamma)\right)} \, .
\end{eqnarray*}

Then the Laplace transform of the occupation field can be written as
\begin{eqnarray*}
\lefteqn{ {\mathbb E}\left[\exp\left(-\sum_{x \in D^{\#}} v_x L_x \right)\right] } \\
& = & \left(\prod_{x \in D^{\#}} (k_x+4)\right)^{1/2} \exp{\left(-\lambda \sum_{\tilde\gamma} \nu^{u,{\bf k}}_D(\tilde\gamma)\right)} \\
& & \left(\prod_{x \in D^{\#}} (v_x+k_x+4)\right)^{-1/2} \exp{\left(\lambda \sum_{\tilde\gamma} \nu^{u,{\bf k}+{\bf v}}_D(\tilde\gamma)\right)} \\
& = & \left(\prod_{x \in D^{\#}} (k_x+4)\right)^{1/2} \exp{\left(-\lambda \sum_{\tilde\gamma} \frac{\rho_{\tilde\gamma}}{|\tilde\gamma|} \prod_{x \in D^{\#}} (k_x+4)^{-n(x,\tilde\gamma)} \right)} \\
& & \left(\prod_{x \in D^{\#}} (v_x+k_x+4)\right)^{-1/2} \exp{\left(\lambda \sum_{\tilde\gamma} \frac{\rho_{\tilde\gamma}}{|\tilde\gamma|} \prod_{x \in D^{\#}} (v_x+k_x+4)^{-n(x,\tilde\gamma)} \right)} \, .
\end{eqnarray*}
If $\lambda=1/2$, using Lemma~1.2 of~\cite{bfsp} and a standard Gaussian integration formula,
the last expression can be written as
\begin{eqnarray*}
\left( \int_{{\mathbb R}^{D^{\#}}} e^{-H^{\bf k}_D(\varphi)} \prod_{x \in D^{\#}}d\varphi_x \right)^{-1}
\int_{{\mathbb R}^{D^{\#}}} \exp{\left( -H^{\bf k}_D(\varphi) - \frac{1}{2} \sum_{x \in D^{\#}} v_x \varphi^2_x \right)} \prod_{x \in D^{\#}}d\varphi_x \\
= {\mathbb E}^{\bf k}_D \left[ \exp\left( -\sum_{x \in D^{\#}} v_x \frac{\phi^2_x}{2} \right) \right] \, ,
\end{eqnarray*}
which concludes the proof. \fbox{}

\section{Appendix: The Brownian loop measure}\label{BrownianLoopMeasure}

The Brownian loop measure studied in \cite{werner3} plays an important role in these lecture notes as the intensity measure used to define
the Brownian loop soup. Indeed, the conformal invariance properties of the latter are a consequence of those of the Brownian loop measure.
The purpose of this appendix is to discuss those properties, starting with some considerations about Brownian motion ad Brownian bridges.

A two-dimensional, or complex, standard Brownian motion $W \equiv (W_t, t \geq 0)$ is the process $W_t = W^1_t + i W^2_t$, where
$(W^1_t, t \geq 0)$ and $(W^2_t, t \geq 0)$ are two independent, one-dimensional, standard Brownian motions. (\emph{Standard} here means
that it starts at the origin and that the quadratic variation is $1$.)

\begin{lemma} \label{conf-inv-BM} \emph{({\bf Conformal invariance of Brownian motion})}
Let $f: D \to {\mathbb C}$ be a conformal map. Let $W$ be a complex Brownian motion started from $z \in D$ and stopped at
\be \nonumber
\tau_D := \inf \{ t \geq 0: W_t \notin D \} \, .
\ee
Define $S(t) := \int_0^t \left| f'(W_u) \right|^2 du$, $0 \leq t \leq \tau_D$ and let $\sigma(t)$ be such that $\int_0^{\sigma(t)} \left| f'(W_u) \right|^2 du = t$.
Then $Y_s = \left( f(W_{\sigma(s)}), 0 \leq s \leq S(\tau_D) \right)$ is a complex Brownian motion started at $f(z)$ and stopped $\tau_{f(D)}$.
\end{lemma}

\noindent{\bf Proof.} Since $f$ is a conformal map, we can write $f = u + i v$, where $u$ and $v$ are harmonic functions such that $\partial_x u = \partial_y v$
and $\partial_y u = - \partial_x v$. Now note that
\begin{eqnarray*}
d u(W_t) & = & \partial_x u(W_t) d W^1_t + \partial_y u(W_t) d W^2_t + \frac{1}{2} \left( \partial_{xx} u(W_t) + \partial_{yy} u(W_t) \right) dt \\
& = & \partial_x u(W_t) d W^1_t + \partial_y u(W_t) d W^2_t \, .
\end{eqnarray*}
Since the $dt$ term vanishes, $u(W_t)$ is a local martingale; the same applies to $v(W_t)$. Moreover, the quadratic variations are
\begin{eqnarray*}
[u(W)]_t = [v(W)]_t & = & \int_0^t \left[ \left(\partial_x u(W_s) \right)^2 + \left(\partial_y u(W_s) \right)^2 \right] \, ds \\
& = & \int_0^t | f'(W_s) |^2 \, ds = S(t) \, ,
\end{eqnarray*}
implying that $Y$ is a complex Brownian motion. The remaining properties follow immediately from the definition of $Y$. \fbox{} \\

Following Lawler and Werner~\cite{lw,lawler-book}, let ${\cal K}$ be the set of parametrized, continuous, planar curves $\gamma$
defined on a time interval $[0,t_{\gamma}]$, endowed with the metric
\begin{equation*}
d_{\cal K}(\gamma,\gamma') := \inf_{\theta} \left[ \sup_{0 \leq s \leq t_{\gamma}}
|s-\theta(s)| + |\gamma(s) - \gamma'(\theta(s))| \right],
\end{equation*}
where the infimum is over all increasing homeomorphisms
$\theta : [0,t_{\gamma}] \to [0,t_{\gamma}]$. Let ${\cal K}(D)$ be the subset of
${\cal K}$ consisting of those curves that stay in $D$. If $z,w \in {\mathbb C}$,
denote by ${\cal K}_z$ (respectively, ${\cal K}^w$) the set of $\gamma \in {\cal K}$
with $\gamma(0)=z$ (resp., $\gamma(t_{\gamma})=w$). We let
${\cal K}_z^w = {\cal K}_z \cup {\cal K}^w$ and define ${\cal K}_z(D)$,
${\cal K}^w(D)$ and ${\cal K}_z^w(D)$ similarly.

Let $\tilde\mu(z,\cdot;t)$ denote the law of standard complex Brownian motion $(W_s, 0 \leq s \leq t)$ with $W_0=z$. We can write
\begin{equation*}
\tilde\mu(z,\cdot;t) = \int_{\mathbb C} \tilde\mu(z,w;t) d{\bf A}(w),
\end{equation*}
where $\bf A$ denotes area and $\tilde\mu(z,w;t)$ is a measure supported on $\gamma \in {\cal K}_z^w$ with $t_{\gamma}=t$.
The total mass of $\tilde\mu(z,w;t)$ is $|\tilde\mu(z,w;t)| = (2 \pi t)^{-1} \exp(-|z-w|^2/(2t))$. The normalized (probability) measure
$\mu^{br}(z,w;t) := \tilde\mu(z,w;t)/|\tilde\mu(z,w;t)|$ is the law of the \emph{Brownian bridge} from $z$ to $w$ in time $t$.

The $\sigma$-finite measure $\tilde\mu(z,w)$ is defined by
\begin{equation*}
\tilde\mu(z,w) := \int_0^{\infty} \tilde\mu(z,w;t) dt = \int_0^{\infty} \frac{1}{2\pi t} \mu^{br}(z,w;t) e^{-|z-w|^2/2t} dt .
\end{equation*}
Note that the integral explodes at infinity so that the total mass of large loops is infinite. (When $z=w$, the integral also diverges at zero).

If $w=z$, one has $\tilde\mu(z,z;t) = (2 \pi t)^{-1} \mu^{br}_{z,t}$, where $\mu^{br}_{z,t} := \mu^{br}(z,z;t)$ is the law of the Brownian
bridge of time duration $t$ started and ended at $z$. The measure
\begin{equation*}
\tilde\mu(z,z) = \int_0^{\infty} \frac{1}{2 \pi t} \mu^{br}_{z,t} \, dt,
\end{equation*}
is an infinite measure on Brownian loops that start and end at $z$.

If $D \subset {\mathbb C}$ is a domain and $z,w \in D$, we define $\tilde\mu_D(z,w)$ to be $\tilde\mu(z,w)$ restricted to ${\cal K}(D)$.
If $z \neq w$ and $D$ is such that Brownian motion started in $D$ eventually exits $D$, then $|\tilde\mu_D(z,w)|<\infty$. More precisely,
$|\tilde\mu_D(z,w)| = \frac{1}{\pi} G_D(z,w)$, where $G_D(z,w)$ is Green's function normalized so that, for the unit disc ${\mathbb D}$,
$G_{\mathbb D}(0,z) = -\log |z|$.

Suppose that $f: D \to D'$ is a conformal transformation and $\gamma \in {\cal K}$; let
\be \label{eq:time-change}
s(t) := \int_0^{t} \left| f'(\gamma(u)) \right|^2 du .
\ee
If $s(t)<\infty$ for all $t < t_{\gamma}$, define $f \circ \gamma$ by
\be \nonumber
(f \circ \gamma)(s(t)) := f(\gamma(t)) .
\ee
Moreover, if $A \in {\cal K}$, let $f \circ A = \{ \hat\gamma = f \circ \gamma \text{ with } \gamma \in A \}$.
Then
\be \label{conf-inv-equation}
\tilde\mu_{f(D)}(f(z),f(w))(f \circ A) = \tilde\mu_D(z,w)(A) .
\ee
The conformal invariance of $\mu_D(z,w)$ expressed by equation \eqref{conf-inv-equation} is an immediate consequence
of two classical results: the conformal invariance of Brownian motion, which follows from Lemma~\ref{conf-inv-BM}, and
the fact that $G_{f(D)}(f(z),f(w)) = G_{D}(z,w)$. The conformal invariance of $\mu(z,z)$ follows by letting $w \to z$.

Now let $\tilde{\cal K}$ be the set of \emph{loops}, i.e., the set of $\gamma \in {\cal K}$
with $\gamma(0)=\gamma(t_{\gamma})$. Such a $\gamma$ can also be considered as a
function with domain $(-\infty,\infty)$ satisfying $\gamma(s)=\gamma(s+t_{\gamma})$.
For $r \in {\mathbb R}$, define the shift operator $\theta_r : \tilde{\cal K} \to \tilde{\cal K}$ by $t_{\theta_r \gamma}=t_{\gamma}$
and $\theta_r \gamma(s) = \gamma(s+r)$. We say that two loops $\gamma$ and $\gamma'$
are equivalent if for some $r$, $\gamma' = \theta_r \gamma$. We write $[\gamma]$ for
the equivalent class of $\gamma$. Let $\tilde{\cal K}_U$ be the set of \emph{unrooted loops},
i.e., the equivalence classes in $\tilde{\cal K}$. Note that $\tilde{\cal K}_U$ is a metric space under the metric
\begin{equation*}
d_{U}(\gamma,\gamma') := \inf_{r \in [0,t_{\gamma}]} d_{\cal K}(\theta_r \gamma, \gamma').
\end{equation*}
and that a measure supported on $\tilde{\cal K}$ gives a measure on $\tilde{\cal K}_U$ by
``forgetting the root,'' i.e., by considering the map $\gamma \rightarrowtail [\gamma]$.
If $D$ is a domain, define $\tilde{\cal K}_U(D)$ to be the set of unrooted loops that
lie entirely in $D$, i.e., $\gamma[0,t_{\gamma}] \subset D$.

The \emph{Brownian loop measure} $\mu$ on $\tilde{\cal K}_U$ is defined by
\begin{equation} \label{equation-blm}
\mu := \int_{\mathbb C} \frac{1}{t_{\gamma}} \tilde\mu(z,z) d{\bf A}(z)
= \int_{\mathbb C} \int_0^{\infty} \frac{1}{2 \pi t^2} \mu^{br}_{z,t} \, dt \, d{\bf A}(z),
\end{equation}
where $dA$ denotes the Lebesgue measure on $\mathbb C$. Informally, one can say that $\mu$ is
the measure on unrooted loops obtained by averaging $\mu$ over the root (starting point):
\be \label{averaged-measure}
``\mu(\gamma) = \frac{1}{t_{\gamma}} \int_{0}^{t_{\gamma}} \tilde\mu(\gamma(s),\gamma(s)) ds ."
\ee

If $D$ is a domain of the complex plane, we define $\mu_D$ to be $\mu$ restricted to curves in
$\tilde{\cal K}(D)$, i.e., $\mu_D$ is given by~\eqref{equation-blm} with $\mathbb C$ replaced by $D$
and $\tilde\mu(z,z)$ replaced by $\tilde\mu_D(z,z)$. The measure $\mu_D$ is conformally invariant in
the sense that, give a conformal map $f:D \to D'$ and a set $A \in \tilde{\cal K}_U(D)$,
\be \label{eq:conf-inv2}
\mu_{D'}(f \circ A) = \mu_D(A) .
\ee

\begin{eqnarray}
\mu_{D'}(f \circ A) & = & \int_{D'} \frac{1}{t_{\hat\gamma}} \tilde\mu_{D'}(w,w)(f \circ A) d{\bf A}(w) \nonumber \\
& = & \int_{D} \frac{1}{t_{f \circ \gamma}} \tilde\mu_{D}(z,z)(A) |f'(z)|^2 d{\bf A}(z) \label{first-step} \\
& = & \int_{D} \frac{|f'(z)|^2}{t_{f \circ \gamma}} \tilde\mu_{D}(z,z)(A) d{\bf A}(z) \label{second-step} \\
& = & \int_{D} \frac{1}{t_{\gamma}} \tilde\mu_{D}(z,z)(A) d{\bf A}(z) = \mu_{D}(A) , \label{last-step}
\end{eqnarray}
where in~\eqref{first-step} we have used the conformal invariance on $\tilde\mu_{D}$ and in~\eqref{last-step}
we have used the fact that
\be \nonumber
\int_0^{t_{\gamma}} \frac{|f'(\gamma(u))|^2}{t_{f \circ \gamma}} \, du = 1 .
\ee
(To understand the calculation above, it may help to remember \eqref{eq:time-change} and \eqref{averaged-measure}
and to think of~\eqref{second-step} as follows:
\be \nonumber
`` \sum_{\gamma \in A} \int_0^{t_{\gamma}} \frac{|f'(\gamma(u))|^2}{t_{f \circ \gamma}} \, \tilde\mu(\gamma(u),\gamma(u)) \, du . \text{''})
\ee


Because the measure $\mu_D$ is the restriction of the measure $\mu$ to loops that stay in $D$, equation \eqref{eq:conf-inv2}
is called \emph{conformal restriction}. The Brownian loop measure is fully characterized by its conformal restriction property, in a way that we
are now going to make precise. In order to do that, we need the following definitions.

\begin{definition} \label{annular-region}
A bounded open set is an \emph{annular region} in the plane if it is conformally equivalent to some annulus $\{ z: 1<|z|<R \}$.
In other words, $A$ is a bounded open set such that ${\mathbb C} \setminus A$ has two connected components (that are not singletons).
\end{definition}

\begin{definition} \label{non-triviality}
A measure on the set of self-avoiding loops in the plane is \emph{nontrivial} if for some $0<\delta<\Delta<\infty$, the mass of loops
of diameter at least $\delta$ that stay in some disc of radius $\Delta$ is neither zero nor infinite. The bounded connected component
of ${\mathbb C} \setminus A$ will be called the \emph{hole} of $A$.
\end{definition}
For each such annular region $A$, define the set ${\cal U}_A$ of self-avoiding loops that stay in $A$ and that have a non-zero index
around the hole of $A$ (or in other words, that disconnects the two connected components of $\partial A$). We will denote by $\cal G$
the $\sigma$-field generated by this family. This is the $\sigma$-field that we will work with.

The next theorem follows from the restriction property of the Brownian loop measure and \cite{werner3}.
\begin{theorem} \label{theorem-conf-restriction}
Up to a multiplicative constant, the Brownian loop measure $\mu$ is the only nontrivial measure on the set of self-avoiding loops
in the plane satisfying conformal restriction. Furthermore, for any simply connected sets $\tilde D \subset D$ and any $z \in \tilde D$,
\be \label{equation-log-over-5}
\mu(\{ \gamma: \gamma \subset D, \gamma \not\subset \tilde D, \gamma \text{ disconnects } z \text{ from } \partial D \}) = \frac{1}{5} \log f'(z)
\ee
where $f$ is the conformal map from $\tilde D$ to $D$ such that $f(z)=z$ and $f'(z)$ is real and positive.
\end{theorem}

The proof of the theorem contains several parts, so we split it into several lemmas.
\begin{lemma} \label{lemma-uniqueness}
Up to a multiplicative constant, there exists at most one non-trivial measure $\mu^{*}$ on the set of self-avoiding loops
in the plane satisfying conformal restriction. Furthermore, for any simply connected sets $\tilde D \subset D$ and any $z \in \tilde D$,
\be \label{log-equation}
\mu^*(\{ \gamma: \gamma \subset D, \gamma \not\subset \tilde D, \gamma \text{ disconnects } z \text{ from } \partial D \}) = c \log f'(z)
\ee
where $f$ is the conformal map from $\tilde D$ to $D$ such that $f(z)=z$, $f'(z)$ is real and positive, and $0<c<\infty$.
\end{lemma}

\noindent{\bf Sketch of the proof.}
This lemma is essentially Proposition 3 of \cite{werner3}; here we omit some of the technical details of the proof, but discuss the main ideas.
Let $\tilde f$ be the unique conformal map from $D$ to the the unit disc $\mathbb U$ such that $\tilde f(z) = 0 $ and $\tilde f'(z)$
is real and positive. Let $U := \tilde f(\tilde D) \subset {\mathbb U}$ and note that $U$ contains the origin $0$.
In the rest of the proof, we will work with simply connected subsets $U$ of the unit disc $\mathbb U$ containing the origin,
and we will prove that conformal restriction implies a version of \eqref{log-equation} for such sets. The conformal invariance
of $\mu^*$ then implies \eqref{log-equation}.

Let $f_U$ be the unique conformal map from $U$ to $\mathbb U$ such that $f_U(0)=0$ and $f'_U(0)>0$. Since $U \subset {\mathbb U}$,
$f'_U(0) \geq 1$. Since $U=f^{-1}({\mathbb U})$, the map $f_U$ describes $U$ fully. Let
\be \nonumber
A(f_U) = \mu^*(\gamma: 0 \in \bar\gamma, \gamma \subset {\mathbb U}, \gamma \not\subset U) ,
\ee
where $0 \in \bar\gamma$ means that $\gamma$ disconnects the origin from infinity, as in the rest of these lecture notes.
It is easy to convince oneself that
\begin{eqnarray*}
A(f_V \circ f_U) & = & \mu^*(\gamma: 0 \in \bar\gamma, \gamma \text{ intersects } {\mathbb U} \setminus U \text{ or }
f^{-1}_U({\mathbb U} \setminus V)) \\
& = & \mu^*( \{ \gamma: 0 \in \bar\gamma, \gamma \subset {\mathbb U}, \gamma \not\subset U \} \cup
\{ \gamma: 0 \in \bar\gamma, \gamma \subset U, \gamma \text{ int. } f^{-1}_U({\mathbb U} \setminus V) \} ) \\
& = & A(f_U) + A(f_V) .
\end{eqnarray*}

Let's now consider a special class of domains: $U_t := {\mathbb U} \setminus [r_t,1)$, where $r_t$ is the positive real number
such that $f'_{U_t}(0) = e^t$. (We won't need the exact expression for $r_t$.) It is a well-know fact that, if $Z$ is a planar
Brownian motion started from the origin and stopped at its first exit time $T$ of $U_t$, then $\log f'_{U_t}(0) = - E \log(|Z_T|)$.
This implies that $t = E(-\log |Z_T|) \in (0,\infty)$.

Note that $f_{U_t} \circ f_{U_s} = f_{U_{\tau}}$ for some $\tau \in (0,\infty)$; but since $(f_{U_t} \circ f_{U_s})' = e^{t+s}$,
it follows that $f_{U_t} \circ f_{U_s} = f_{U_{t+s}}$, so the family $(f_{U_t})_{t \geq 0}$ is a semi-group. Hence, $t \mapsto A(f_{U_t})$
is a non-decreasing function from $(0,\infty)$ to $(0,\infty)$ such that
\be \nonumber
A(f_{U_{t+s}}) = A(f_{U_t} \circ f_{U_s}) = A(f_{U_t}) + A(f_{U_s}) .
\ee
This implies that $A(f_{U_t}) = c t = c \log f'_{U_t}(0)$ for some constant $0<c<\infty$.

For each $\theta \in [0,2\pi)$ and $t>0$, let $U_{t,\theta} :=  \setminus [r_t e^{i \theta}, e^{i \theta})$;
clearly, $f'_{U_{t,\theta}}(0) = e^t$ and $A(f_{U_{t,\theta}}) = A(f_{U_t}) = c t$. Let $S$ denote the semi-group of conformal maps
generated by the family $(f_{U_{t,\theta}}, t>0, \theta \in [0,2\pi))$ (i.e., the set of finite compositions of those maps).
If we take $f = f_{U_{t,\theta}} \circ f_{U_{s,\alpha}} \in S$, then $f'(0) = f'_{U_{t,\theta}}(0) f'_{U_{s,\alpha}}(0) = e^{t+s}$, that is,
$c \log f'(0) = c (t+s)$. On the other hand, $A(f) = A(f_{U_{t,\theta}}) + A(f_{U_{s,\alpha}}) = c t + c s$, therefore $A(f) = c \log f'(0)$.

The theory developed by Loewner (for the proof of the Bieberbach conjecture in the first highly nontrivial case of the third coefficient)
implies that $S$ is ``dense'' in the class of conformal maps $f_U$ from some simply connected domain $U \subset \mathbb U$ onto
$\mathbb U$ in the sense that, for any $f_U$, one can find a sequence $f_{U_n}$ in $S$ such that $U_n$ is an increasing in $n$
(in the sense of inclusion) family of domains such that $\cup_n U_n = U$. Therefore, $f'_{U_n}(0) \to f'_U(0)$ as $n \to \infty$.
Furthermore, because a loop is a compact subset of the complex plane, it exits $U$ if and only if it exits $U_n$ for every $n$.
Thus, we have that, for any simply connected subset $U$ of $\mathbb U$ containing the origin, we have that
\be \nonumber
A(f_U) = \lim_{n \to \infty} A(f_{U_n}) = c \lim_{n \to \infty} \log f'_{U_n}(0) = c \log f'_U(0) .
\ee

The conformal invariance of $\mu^*$ concludes the proof of \eqref{log-equation}. The uniqueness of $\mu^*$ follows from
rather simple measure-theoretical considerations that we omit. (The interest reader can find the details in \cite{werner3}.)
 \fbox{} \\

In the next lemmas, as in the rest of these lecture notes, we let $\bar \gamma$ denote the ``filling" of the loop $\gamma$,
$B_{z, a}$ be a disk of radius $a$ around $z$ and $\gamma \not\subset B_{z, \delta}$ indicate that the image of $\gamma$
is not fully contained in $B_{z, \delta}$.
\begin{lemma} \label{lemma-measure-equality}
Let $z \in \mathbb C$, then
\be \nonumber
\mu(\gamma: z \in \bar \gamma, \delta \leq \diam(\gamma) < R)
=\mu(\gamma: z \in \bar \gamma,\gamma \not\subset B_{z, \delta}, \gamma \subset B_{z,R}) .
\ee
\end{lemma}

\noindent{\bf Proof.}
Since $\mu(\gamma: \diam(\gamma)=R)=0$, we have that
\begin{eqnarray*}
\lefteqn{\mu(\gamma: z \in \bar \gamma, \delta \leq \diam(\gamma) < R)
- \mu(\gamma: z \in \bar \gamma,\gamma \not\subset B_{z, \delta}, \gamma \subset B_{z,R})} \\
& = & \mu(\gamma: z \in \bar \gamma, \diam(\gamma) \geq \delta, \gamma \subset B_{z,\delta})
- \mu(\gamma: z \in \bar \gamma,\diam(\gamma) \geq R, \gamma \subset B_{z,R}) ,
\end{eqnarray*}
where the last two terms are identical because of the scale invariance of $\mu$. \fbox{}

\begin{lemma} \label{lemma-r-over-5}
For any $z \in {\mathbb C}$ and any $r>0$,
\begin{equation} \nonumber
\mu\left(\gamma: z \in \bar \gamma, 1 \leq \diam(\gamma) < e^r \right)
= \mu\left(\gamma: z \in \bar \gamma, 1 \leq t_{\gamma} < e^{2r}\right) = \frac{r}{5} .
\end{equation}
\end{lemma}

\noindent{\bf Proof.}
Part of the lemma follows from a simple computation:
\begin{eqnarray*}
\mu(\gamma: 0 \in \bar \gamma, 1 \leq t_{\gamma} < e^{2r})
& = & \int_{\mathbb C} \int_1^{e^{2r}} \frac{1}{2 \pi t^2} \,
\mu^{br}_{z,t}(\{ \gamma: 0 \in \bar \gamma \}) \, dt \, d{\bf A}(z) \\
& = & \int_{\mathbb C} \int_1^{e^{2r}} \frac{1}{2 \pi t^2} \,
\mu^{br}_{0,t}(\{ \gamma: z \in \bar \gamma \}) \, dt \, d{\bf A}(z) \\
& = & \int_1^{e^{2r}} \frac{1}{2 \pi t^2} \,
{\mathbb E}^{br}_{0,t}\left( \int_{\mathbb C} \mathbbm{1}_{\{ \gamma: z \in \bar \gamma \}} \, d{\bf A}(z) \right) \, dt \\
& = & \int_1^{e^{2r}} \frac{1}{2 \pi t} \,
{\mathbb E}^{br}_{0,1}\left( \int_{\mathbb C} \mathbbm{1}_{\{ \gamma: z \in \bar \gamma \}} \, d{\bf A}(z) \right) \, dt ,
\end{eqnarray*}
where ${\mathbb E}^{br}_{0,t}$ denotes expectation with respect to a complex Brownian bridge of time length $t$
started at the origin, and where, in the last equality, we have used the fact that
\begin{equation} \nonumber
{\mathbb E}^{br}_{0,t}\left( \int_{\mathbb C} \mathbbm{1}_{\{ \gamma: z \in \bar \gamma \}} \, d{\bf A}(z) \right)
= t \, {\mathbb E}^{br}_{0,1}\left( \int_{\mathbb C} \mathbbm{1}_{\{ \gamma: z \in \bar \gamma \}} \, d{\bf A}(z) \right)
\end{equation}
because of scaling. The expected area of a ``filled-in'' Brownian bridge, computed in \cite{gtf}, is
\begin{equation} \nonumber
{\mathbb E}^{br}_{0,1}\left( \int_{\mathbb C} \mathbbm{1}_{\{ \gamma: z \in \bar \gamma \}} \, d{\bf A}(z) \right)
= \frac{\pi}{5} ,
\end{equation}
so that
\begin{equation} \label{frac-equation}
\mu(\gamma: z \in \bar \gamma, 1 \leq t_{\gamma} < e^{2r}) = \frac{r}{5} .
\end{equation}

Combined with  equation \eqref{log-equation}, equation \eqref{frac-equation} implies that
\begin{equation} \nonumber
\frac{\mu(\gamma: z \in \bar \gamma, 1 \leq \diam(\gamma) < e^r)}{\mu(\gamma: z \in \bar \gamma, 1 \leq t_{\gamma} < e^{2r})} = 5c \, ,
\end{equation}
independently of $r$. Therefore, to conclude the proof of the lemma, it suffices to show that
\begin{equation} \nonumber
\lim_{r \to \infty}\frac{\mu(\gamma: z \in \bar \gamma, 1 \leq \diam(\gamma) < e^r)}{\mu(\gamma: z \in \bar \gamma, 1 \leq t_{\gamma} < e^{2r})} = 1 .
\end{equation}

Let $A := \{ \gamma: z \in \bar \gamma, 1 \leq \diam(\gamma) < e^r \}$,
$B \equiv \{ \gamma: z \in \bar \gamma, 1 \leq t_{\gamma} < e^{2r} \}$, and define the disjoint sets
\begin{itemize}
\item $B_1 := \{ \gamma: z \in \bar \gamma, \diam(\gamma) < e^{-r}, 1 \leq t_{\gamma} < e^{2r} \} \subset B$,
\item $B_2 := \{ \gamma: z \in \bar \gamma, e^{-r} \leq \diam(\gamma) < 1 , 1 \leq t_{\gamma} < e^{2r} \} \subset B$,
\item $B_3 := \{ \gamma: z \in \bar \gamma, 1 \leq \diam(\gamma) < e^r , 1 \leq t_{\gamma} < e^{2r} \}
= A \cap B$,
\item $B_4 := \{ \gamma: z \in \bar \gamma, e^{r} \leq \diam(\gamma) < e^{2r} , 1 \leq t_{\gamma} < e^{2r} \} \subset B$,
\item $B_5 := \{ \gamma: z \in \bar \gamma, \diam(\gamma) \geq e^{2r} , 1 \leq t_{\gamma} < e^{2r} \} \subset B$,
\end{itemize}
and the disjoint sets
\begin{itemize}
\item $A_1 := \{ \gamma: z \in \bar \gamma, 1 \leq \diam(\gamma) < e^r , t_{\gamma} < e^{-2r} \} \subset A$,
\item $A_2 := \{ \gamma: z \in \bar \gamma, 1 \leq \diam(\gamma) < e^r , e^{-2r} \leq t_{\gamma} < 1 \}
\subset A$,
\item $A_3 := \{ \gamma: z \in \bar \gamma, 1 \leq \diam(\gamma) < e^r , 1 \leq t_{\gamma} < e^{2r} \}
= A \cap B$,
\item $A_4 := \{ \gamma: z \in \bar \gamma, 1 \leq \diam(\gamma) < e^r , e^{2r} \leq t_{\gamma} < e^{4r} \}
\subset A$,
\item $A_5 := \{ \gamma: z \in \bar \gamma, 1 \leq \diam(\gamma) < e^r , t_{\gamma} \geq e^{4r} \}
\subset A$ .
\end{itemize}
We clearly have that $\mu(A) = \sum_{i=1}^5 \mu(A_i)$ and
$\mu(B) = \sum_{i=1}^5 \mu(B_i)$.

Note that $A_4$ can be obtained from $B_2$ by scaling each loop $\gamma$ in $B_2$ by a factor of $e^r$ and scaling
$t_{\gamma}$ by a factor $e^{2r}$. Similarly, $A_2$ can be obtained from $B_4$ by scaling each loop $\gamma$ in
$B_4$ by a factor of $e^{-r}$ and scaling $t_{\gamma}$ by a factor $e^{-2r}$. Because of the scaling properties of
Brownian motion, those transformations are $\mu$-measure preserving; thus
$\mu(A_4) = \mu(B_2)$ and $\mu(A_2) = \mu(B_4)$.

To conclude the proof, we show that the $\mu$-measures of $B_1, B_5, A_1, A_5$ go to zero as $r \to \infty$.
For $\mu(B_5)$, we have the upper bound
\begin{eqnarray*}
\mu(B_5) & = & \int_{\mathbb C} \int_1^{e^{2r}} \frac{1}{2\pi t^2} \, \mu^{br}_{z,t}(\gamma: 0 \in \bar\gamma, \diam(\gamma) \geq e^{2r}) \, dt \, d{\bf A}(z) \\
& \leq & \pi e^{4r} \int_1^{e^{2r}} \frac{1}{2\pi t^2} \,
\mu^{br}_{0,t}\left(\gamma: \diam(\gamma) \geq e^{2r} \right) \, dt \\
& + & \sum_{n=1}^{\infty} \pi \left[\left(e^{2r}+n \right)^2 - \left(e^{2r}+n-1 \right)^2 \right] \\
& & \int_1^{e^{2r}} \frac{1}{2\pi t^2} \,
\mu^{br}_{0,t}\left(\gamma: \diam(\gamma) \geq e^{2r}+n-1 \right) \, dt ,
\end{eqnarray*}
where the first term in the upper bound comes from rooted loops with root within the disc of radius $e^{2r}$
centered at the origin, and the other terms come from rooted loops with roots inside one of a family of concentric
annuli around that disc.

Plugging \eqref{brownian-bound} from the proof of Lemma~\ref{no-loop} into the upper bound for $\mu(B_5)$ gives
\begin{eqnarray*}
\mu(B_5) & \leq & 2e^{4r} \int_1^{e^{2r}} t^{-2} \, \exp\left( -\frac{e^{4r}}{288 t} \right) \, dt \\
& + & 2 \sum_{n=1}^{\infty} (2e^{2r} + 2n -1) \int_1^{e^{2r}} t^{-2}
\, \exp\left( -\frac{\left(e^{2r}+n-1 \right)^2}{288 t} \right) \, dt , \\
& \leq & \text{const} \, \exp\left( -\frac{e^{2r}}{288} \right) \to 0 \; \text{ as } r \to \infty .
\end{eqnarray*}

If we now let $W_s$ and $B_s := W_s - s W_t$, $s \in [0,t]$, denote standard two-dimensional Brownian motion and
Brownian bridge, respectively, we have that $\max_{s \in [0,t]} |B_s| \geq \left| W_{1/2} - \frac{1}{2} W_t \right|$.
Combining this observation with the Markov property of Brownian motion and the Gaussian distribution of its increments leads
to the following upper bound for $\mu(B_1)$:
\begin{eqnarray*}
\mu(B_1) & = & \int_{\mathbb C} \int_1^{e^{2r}} \frac{1}{2\pi t^2} \, \mu^{br}_{z,t}\left(\gamma: 0 \in \bar\gamma, \diam(\gamma) < e^{-r} \right) \, dt \, d{\bf A}(z) \\
& \leq & \pi e^{-2r} \int_1^{e^{2r}} \Pr\left( \left| W_t - 2W_{1/2} \right| \leq e^{-r} \right) \, dt \\
& \leq & \text{const} \, e^{-2r} \to 0 \; \text{ as } r \to \infty .
\end{eqnarray*}

Finally, let $B'_5 \equiv \{ \gamma: z \in \bar \gamma, e^{2r} \leq \diam(\gamma) < e^{3r} , t_{\gamma} < e^{2r} \} \subset B_5$
and $B'_1 \equiv \{ \gamma: z \in \bar \gamma, e^{-2r} \leq \diam(\gamma) < e^{-r}, t_{\gamma} \geq 1 \} \subset B_1$.
Note that $B'_5$ can be obtained from $A_1$ by scaling each loop $\gamma$ in $A_1$ by a factor of $e^{2r}$ and scaling
$t_{\gamma}$ by a factor $e^{4r}$. Similarly, $B'_1$ can be obtained from $A_5$ by scaling each loop $\gamma$ in $A_5$
by a factor of $e^{-2r}$ and scaling $t_{\gamma}$ by a factor $e^{-4r}$. Because of the scaling properties of Brownian motion,
those transformations are $\mu$-measure preserving; thus $\mu(A_1) = \mu(B'_5) \leq \mu(B_5)$
and $\mu(A_5) = \mu(B'_1) \leq \mu(B_1)$, concluding the proof of the lemma. \fbox{} \\

\noindent{\bf Proof of Theorem \ref{theorem-conf-restriction}.}
We have already seen that the Brownian loop measure $\mu$ satisfies conformal restriction; Lemma \ref{lemma-uniqueness} then
implies that $\mu$ is unique, up to a multiplicative constant, and satisfies \eqref{log-equation} for some constant $0<c<\infty$.
Combining \eqref{log-equation} with Lemmas \ref{lemma-measure-equality} and \ref{lemma-r-over-5} gives
\begin{eqnarray*}
c r & = & \mu(\gamma: z \in \bar \gamma,\gamma \not\subset B_{z,1}, \gamma \subset B_{z,e^r}) \\
& = & \mu\left(\gamma: z \in \bar \gamma, 1 \leq \diam(\gamma) < e^r \right) \\
& = & \mu\left(\gamma: z \in \bar \gamma, 1 \leq t_{\gamma} < e^{2r}\right)  = \frac{r}{5} . \; \fbox{}
\end{eqnarray*}

We conclude this appendix with two simple but useful lemmas.
\begin{lemma} \label{FirstLemma}
Let $z \in \mathbb C$, then
\be \nonumber
\mu(\gamma: z \in \bar \gamma, \delta \leq \diam(\gamma) < R)
= \frac{1}{5} \log\frac{R}{\delta} .
\ee
\end{lemma}

\noindent{\bf Proof.}
This is an immediate consequence of Lemma \ref{lemma-measure-equality} and Theorem \ref{theorem-conf-restriction}. \fbox{}

\begin{lemma} \label{SecondLemma}
Let $z \in \mathbb C$ and $k \in {\mathbb Z} \setminus \{0\}$, then
\be \nonumber
\mu^{loop}(\gamma: \gamma \text{ has winding number } k \text{ around } z, \delta \leq \diam(\gamma) < R)
= \frac{1}{2\pi^2 k^2} \log\frac{R}{\delta} .
\ee
\end{lemma}

\noindent{\bf Proof.} It is easy to check that the measure on loops surrounding the origin induced by $\mu$,
but restricted to loops that wind $k$ times around a given $z \in {\mathbb C}$, satisfies the conformal restriction
property. Therefore, following the proofs of lemmas \ref{lemma-measure-equality} and \ref{lemma-uniqueness},
we obtain
\begin{eqnarray*}
\lefteqn{\mu(\gamma: \gamma \text{ has winding number } k \text{ around } z, \delta \leq \diam(\gamma) < R)} \\
& = & \mu(\gamma: \gamma \text{ has winding number } k \text{ around } z, \gamma \not\subset B_{z, \delta}, \gamma \subset B_{z,R}) \nonumber \\
& = & c_k \log\frac{R}{\delta} ,
\end{eqnarray*}
for some positive constant $c_k<\infty$. In order to find the constants $c_k$, we can proceed as in the proof
of Theorem~\ref{theorem-conf-restriction}, using the fact that the expected area of a ``filled-in'' Brownian loop
winding $k$ times around the origin was computed in \cite{gtf} and is equal to $1/2\pi k^2$
for $k \in {\mathbb Z} \setminus \{0\}$ (and $\pi/30$ for $k=0$). \fbox{}


\end{document}